\pdfoutput=1

\documentclass[12pt,a4paper]{article}

\usepackage{amsfonts,amssymb,amsthm,amsmath,latexsym}
\usepackage{mathtools}
\usepackage{subeqnarray,eqsection,indent,url,cite,bm}
\usepackage{multirow}  
\usepackage{tensor}  
\usepackage{graphicx}
\usepackage{float}  

\date{May 10, 2021 \\[2mm] revised March 3, 2022}   

\begin{document}

\title{\vspace*{-1cm}Total positivity of some polynomial matrices \\[0.5mm]
       that enumerate labeled trees and forests \\[5mm]
       I.~Forests of rooted labeled trees
      }

\author{ \\
      \hspace*{-1cm}
      {\large Alan D.~Sokal${}^{1,2}$}
   \\[5mm]
     \hspace*{-2cm}
      \normalsize
           ${}^1$Department of Mathematics, University College London,
                    Gower Street, London WC1E 6BT, UK   \\[1mm]
     \hspace*{-2.4cm}
      \normalsize
           ${}^2$Department of Physics, New York University,
                    726 Broadway, New York, NY 10003, USA
       \\[3mm]
      \hspace*{-1.2cm}
      \normalsize {\tt sokal@nyu.edu} \\[1cm]
}

\maketitle
\thispagestyle{empty}   

\begin{abstract}
We consider the lower-triangular matrix of generating polynomials
that enumerate $k$-component forests of rooted trees on the vertex set $[n]$
according to the number of improper edges
(generalizations of the Ramanujan polynomials).
We show that this matrix is coefficientwise totally positive
and that the sequence of its row-generating polynomials
is coefficientwise Hankel-totally positive.
More generally, we define the generic rooted-forest polynomials
by introducing also a weight $m! \, \phi_m$
for each vertex with $m$~proper children.
We show that if the weight sequence~$\bm{\phi}$ is Toeplitz-totally positive,
then the two foregoing total-positivity results continue to hold.
Our proofs use production matrices and exponential Riordan arrays.
\end{abstract}


\bigskip
\bigskip
\noindent
{\bf Mathematics Subject Classification (MSC 2010) codes:}
05A15 (Primary);
05A19, 05A20, 05C05, 05C30, 15B05, 15B36, 15B48, 30E05, 44A60 (Secondary).

\clearpage

\newtheorem{theorem}{Theorem}[section]
\newtheorem{proposition}[theorem]{Proposition}
\newtheorem{lemma}[theorem]{Lemma}
\newtheorem{corollary}[theorem]{Corollary}
\newtheorem{definition}[theorem]{Definition}
\newtheorem{conjecture}[theorem]{Conjecture}
\newtheorem{question}[theorem]{Question}
\newtheorem{problem}[theorem]{Problem}
\newtheorem{example}[theorem]{Example}
\newtheorem{remark}[theorem]{Remark}

\renewcommand{\theenumi}{\alph{enumi}}
\renewcommand{\labelenumi}{(\theenumi)}
\def\eop{\hbox{\kern1pt\vrule height6pt width4pt
depth1pt\kern1pt}\medskip}
\def\prf{\par\noindent{\bf Proof.\enspace}\rm}
\def\rmk{\par\medskip\noindent{\bf Remark\enspace}\rm}

\newcommand{\textbfit}[1]{\textbf{\textit{#1}}}

\newcommand{\bigdash}{%
\smallskip\begin{center} \rule{5cm}{0.1mm} \end{center}\smallskip}

\newcommand{\safepar}{ {\protect\hfill\protect\break\hspace*{5mm}} }

\newcommand{\be}{\begin{equation}}
\newcommand{\ee}{\end{equation}}
\newcommand{\<}{\langle}
\renewcommand{\>}{\rangle}
\newcommand{\widebar}{\overline}
\def\reff#1{(\protect\ref{#1})}
\def\spose#1{\hbox to 0pt{#1\hss}}
\def\ltapprox{\mathrel{\spose{\lower 3pt\hbox{$\mathchar"218$}}
    \raise 2.0pt\hbox{$\mathchar"13C$}}}
\def\gtapprox{\mathrel{\spose{\lower 3pt\hbox{$\mathchar"218$}}
    \raise 2.0pt\hbox{$\mathchar"13E$}}}
\def\textprime{${}^\prime$}
\def\proof{\par\medskip\noindent{\sc Proof.\ }}
\def\firstproof{\par\medskip\noindent{\sc First Proof.\ }}
\def\secondproof{\par\medskip\noindent{\sc Second Proof.\ }}
\def\alternateproof{\par\medskip\noindent{\sc Alternate Proof.\ }}
\def\algebraicproof{\par\medskip\noindent{\sc Algebraic Proof.\ }}
\def\combinatorialproof{\par\medskip\noindent{\sc Combinatorial Proof.\ }}
\def\proofof#1{\bigskip\noindent{\sc Proof of #1.\ }}
\def\firstproofof#1{\bigskip\noindent{\sc First Proof of #1.\ }}
\def\secondproofof#1{\bigskip\noindent{\sc Second Proof of #1.\ }}
\def\thirdproofof#1{\bigskip\noindent{\sc Third Proof of #1.\ }}
\def\algebraicproofof#1{\bigskip\noindent{\sc Algebraic Proof of #1.\ }}
\def\combinatorialproofof#1{\bigskip\noindent{\sc Combinatorial Proof of #1.\ }}
\def\sketchofproof{\par\medskip\noindent{\sc Sketch of proof.\ }}
\def\sketchofproofof#1{\par\medskip\noindent{\sc Sketch of proof of #1.\ }}
\renewcommand{\qed}{ $\square$ \bigskip}
\newcommand{\myendremark}{ $\blacksquare$ \bigskip}
\def\half{ {1 \over 2} }
\def\third{ {1 \over 3} }
\def\twothird{ {2 \over 3} }
\def\smfrac#1#2{{\textstyle{#1\over #2}}}
\def\smhalf{ {\smfrac{1}{2}} }
\newcommand{\real}{\mathop{\rm Re}\nolimits}
\renewcommand{\Re}{\mathop{\rm Re}\nolimits}
\newcommand{\imag}{\mathop{\rm Im}\nolimits}
\renewcommand{\Im}{\mathop{\rm Im}\nolimits}
\newcommand{\sgn}{\mathop{\rm sgn}\nolimits}
\newcommand{\tr}{\mathop{\rm tr}\nolimits}
\newcommand{\supp}{\mathop{\rm supp}\nolimits}
\newcommand{\disc}{\mathop{\rm disc}\nolimits}
\newcommand{\diag}{\mathop{\rm diag}\nolimits}
\newcommand{\tridiag}{\mathop{\rm tridiag}\nolimits}
\newcommand{\AZ}{\mathop{\rm AZ}\nolimits}
\newcommand{\perm}{\mathop{\rm perm}\nolimits}
\def\hboxscript#1{ {\hbox{\scriptsize\em #1}} }
\renewcommand{\emptyset}{\varnothing}
\newcommand{\eqdef}{\stackrel{\rm def}{=}}

\newcommand{\restrict}{\upharpoonright}

\newcommand{\compinv}{{\langle -1 \rangle}}   

\newcommand{\scra}{{\mathcal{A}}}
\newcommand{\scrb}{{\mathcal{B}}}
\newcommand{\scrc}{{\mathcal{C}}}
\newcommand{\scrd}{{\mathcal{D}}}
\newcommand{\scre}{{\mathcal{E}}}
\newcommand{\scrf}{{\mathcal{F}}}
\newcommand{\scrg}{{\mathcal{G}}}
\newcommand{\scrh}{{\mathcal{H}}}
\newcommand{\scri}{{\mathcal{I}}}
\newcommand{\scrj}{{\mathcal{J}}}
\newcommand{\scrk}{{\mathcal{K}}}
\newcommand{\scrl}{{\mathcal{L}}}
\newcommand{\scrm}{{\mathcal{M}}}
\newcommand{\scrn}{{\mathcal{N}}}
\newcommand{\scro}{{\mathcal{O}}}
\newcommand\scroo{
  \mathchoice
    {{\scriptstyle\mathcal{O}}}
    {{\scriptstyle\mathcal{O}}}
    {{\scriptscriptstyle\mathcal{O}}}
    {\scalebox{0.6}{$\scriptscriptstyle\mathcal{O}$}}
  }
\newcommand{\scrp}{{\mathcal{P}}}
\newcommand{\scrq}{{\mathcal{Q}}}
\newcommand{\scrr}{{\mathcal{R}}}
\newcommand{\scrs}{{\mathcal{S}}}
\newcommand{\scrt}{{\mathcal{T}}}
\newcommand{\scrv}{{\mathcal{V}}}
\newcommand{\scrw}{{\mathcal{W}}}
\newcommand{\scrz}{{\mathcal{Z}}}

\newcommand{\bfa}{{\mathbf{a}}}
\newcommand{\bfb}{{\mathbf{b}}}
\newcommand{\bfc}{{\mathbf{c}}}
\newcommand{\bfd}{{\mathbf{d}}}
\newcommand{\bfe}{{\mathbf{e}}}
\newcommand{\bfF}{{\mathbf{F}}}
\newcommand{\bfh}{{\mathbf{h}}}
\newcommand{\bfj}{{\mathbf{j}}}
\newcommand{\bfi}{{\mathbf{i}}}
\newcommand{\bfk}{{\mathbf{k}}}
\newcommand{\bfl}{{\mathbf{l}}}
\newcommand{\bfm}{{\mathbf{m}}}
\newcommand{\bfn}{{\mathbf{n}}}
\newcommand{\bfp}{{\mathbf{p}}}
\newcommand{\bfx}{{\mathbf{x}}}
\newcommand{\bfy}{{\mathbf{y}}}
\renewcommand{\k}{{\mathbf{k}}}
\newcommand{\n}{{\mathbf{n}}}
\newcommand{\vv}{{\mathbf{v}}}
\newcommand{\bv}{{\mathbf{v}}}
\newcommand{\w}{{\mathbf{w}}}
\newcommand{\x}{{\mathbf{x}}}
\newcommand{\y}{{\mathbf{y}}}
\newcommand{\cc}{{\mathbf{c}}}
\newcommand{\zero}{{\mathbf{0}}}
\newcommand{\one}{{\mathbf{1}}}
\newcommand{\bmm}{{\mathbf{m}}}

\newcommand{\ahat}{{\widehat{a}}}
\newcommand{\Zhat}{{\widehat{Z}}}

\newcommand{\C}{{\mathbb C}}
\newcommand{\D}{{\mathbb D}}
\newcommand{\Z}{{\mathbb Z}}
\newcommand{\N}{{\mathbb N}}
\newcommand{\Q}{{\mathbb Q}}
\newcommand{\PP}{{\mathbb P}}
\newcommand{\R}{{\mathbb R}}
\newcommand{\RR}{{\mathbb R}}
\newcommand{\E}{{\mathbb E}}

\newcommand{\Sym}{{\mathfrak{S}}}
\newcommand{\SymB}{{\mathfrak{B}}}
\newcommand{\Alt}{{\mathrm{Alt}}}

\newcommand{\germanA}{{\mathfrak{A}}}
\newcommand{\germanB}{{\mathfrak{B}}}
\newcommand{\germanQ}{{\mathfrak{Q}}}
\newcommand{\germanh}{{\mathfrak{h}}}

\newcommand{\myle}{\preceq}
\newcommand{\myge}{\succeq}
\newcommand{\mygt}{\succ}

\newcommand{\B}{{\sf B}}
\newcommand{\OB}{B^{\rm ord}}
\newcommand{\OS}{{\sf OS}}
\newcommand{\OO}{{\sf O}}
\newcommand{\SP}{{\sf SP}}
\newcommand{\OSP}{{\sf OSP}}
\newcommand{\Eu}{{\sf Eu}}
\newcommand{\ERR}{{\sf ERR}}
\newcommand{\sfB}{{\sf B}}
\newcommand{\sfD}{{\sf D}}
\newcommand{\sfE}{{\sf E}}
\newcommand{\sfF}{{\sf F}}
\newcommand{\sfG}{{\sf G}}
\newcommand{\sfJ}{{\sf J}}
\newcommand{\sfL}{{\sf L}}
\newcommand{\sfP}{{\sf P}}
\newcommand{\sfQ}{{\sf Q}}
\newcommand{\sfS}{{\sf S}}
\newcommand{\sfT}{{\sf T}}
\newcommand{\sfW}{{\sf W}}
\newcommand{\sfMV}{{\sf MV}}
\newcommand{\AMV}{{\sf AMV}}
\newcommand{\BM}{{\sf BM}}
\newcommand{\NC}{{\sf NC}}

\newcommand{\emIB}{B^{\rm irr}}
\newcommand{\emIP}{P^{\rm irr}}
\newcommand{\emOB}{B^{\rm ord}}
\newcommand{\emCB}{B^{\rm cyc}}
\newcommand{\emSC}{P^{\rm cyc}}

\newcommand{\propv}{{\rm propv}}
\newcommand{\stat}{{\rm stat}}
\newcommand{\cyc}{{\rm cyc}}
\newcommand{\Asc}{{\rm Asc}}
\newcommand{\asc}{{\rm asc}}
\newcommand{\Des}{{\rm Des}}
\newcommand{\des}{{\rm des}}
\newcommand{\Exc}{{\rm Exc}}
\newcommand{\exc}{{\rm exc}}
\newcommand{\Wex}{{\rm Wex}}
\newcommand{\wex}{{\rm wex}}
\newcommand{\Fix}{{\rm Fix}}
\newcommand{\fix}{{\rm fix}}
\newcommand{\lrmax}{{\rm lrmax}}
\newcommand{\rlmax}{{\rm rlmax}}
\newcommand{\Rec}{{\rm Rec}}
\newcommand{\rec}{{\rm rec}}
\newcommand{\Arec}{{\rm Arec}}
\newcommand{\arec}{{\rm arec}}
\newcommand{\ERec}{{\rm ERec}}
\newcommand{\erec}{{\rm erec}}
\newcommand{\EArec}{{\rm EArec}}
\newcommand{\earec}{{\rm earec}}
\newcommand{\recarec}{{\rm recarec}}
\newcommand{\nonrec}{{\rm nonrec}}
\newcommand{\Cpeak}{{\rm Cpeak}}
\newcommand{\cpeak}{{\rm cpeak}}
\newcommand{\Cval}{{\rm Cval}}
\newcommand{\cval}{{\rm cval}}
\newcommand{\Cdasc}{{\rm Cdasc}}
\newcommand{\cdasc}{{\rm cdasc}}
\newcommand{\Cddes}{{\rm Cddes}}
\newcommand{\cddes}{{\rm cddes}}
\newcommand{\cdrise}{{\rm cdrise}}
\newcommand{\cdfall}{{\rm cdfall}}
\newcommand{\Peak}{{\rm Peak}}
\newcommand{\peak}{{\rm peak}}
\newcommand{\Val}{{\rm Val}}
\newcommand{\val}{{\rm val}}
\newcommand{\Dasc}{{\rm Dasc}}
\newcommand{\dasc}{{\rm dasc}}
\newcommand{\Ddes}{{\rm Ddes}}
\newcommand{\ddes}{{\rm ddes}}
\newcommand{\inv}{{\rm inv}}
\newcommand{\maj}{{\rm maj}}
\newcommand{\rs}{{\rm rs}}
\newcommand{\cross}{{\rm cr}}
\newcommand{\crosshat}{{\widehat{\rm cr}}}
\newcommand{\nest}{{\rm ne}}
\newcommand{\rodd}{{\rm rodd}}
\newcommand{\reven}{{\rm reven}}
\newcommand{\lodd}{{\rm lodd}}
\newcommand{\leven}{{\rm leven}}
\newcommand{\sg}{{\rm sg}}
\newcommand{\bl}{{\rm bl}}
\newcommand{\tran}{{\rm tr}}
\newcommand{\area}{{\rm area}}
\newcommand{\ret}{{\rm ret}}
\newcommand{\peaks}{{\rm peaks}}
\newcommand{\hl}{{\rm hl}}
\newcommand{\sll}{{\rm sl}}
\newcommand{\negg}{{\rm neg}}
\newcommand{\imp}{{\rm imp}}
\newcommand{\pc}{{\rm pc}}
\newcommand{\osg}{{\rm osg}}
\newcommand{\ons}{{\rm ons}}
\newcommand{\isg}{{\rm isg}}
\newcommand{\ins}{{\rm ins}}
\newcommand{\LL}{{\rm LL}}
\newcommand{\height}{{\rm ht}}
\newcommand{\as}{{\rm as}}

\newcommand{\ba}{{\bm{a}}}
\newcommand{\bahat}{{\widehat{\bm{a}}}}
\newcommand{\sfa}{{{\sf a}}}
\newcommand{\bb}{{\bm{b}}}
\newcommand{\bc}{{\bm{c}}}
\newcommand{\bchat}{{\widehat{\bm{c}}}}
\newcommand{\bd}{{\bm{d}}}
\newcommand{\bee}{{\bm{e}}}
\newcommand{\beh}{{\bm{eh}}}
\newcommand{\bff}{{\bm{f}}}
\newcommand{\bg}{{\bm{g}}}
\newcommand{\bh}{{\bm{h}}}
\newcommand{\bll}{{\bm{\ell}}}
\newcommand{\bp}{{\bm{p}}}
\newcommand{\br}{{\bm{r}}}
\newcommand{\bs}{{\bm{s}}}
\newcommand{\bu}{{\bm{u}}}
\newcommand{\bw}{{\bm{w}}}
\newcommand{\bx}{{\bm{x}}}
\newcommand{\by}{{\bm{y}}}
\newcommand{\bz}{{\bm{z}}}
\newcommand{\bA}{{\bm{A}}}
\newcommand{\bB}{{\bm{B}}}
\newcommand{\bC}{{\bm{C}}}
\newcommand{\bE}{{\bm{E}}}
\newcommand{\bF}{{\bm{F}}}
\newcommand{\bG}{{\bm{G}}}
\newcommand{\bH}{{\bm{H}}}
\newcommand{\bI}{{\bm{I}}}
\newcommand{\bJ}{{\bm{J}}}
\newcommand{\bM}{{\bm{M}}}
\newcommand{\bN}{{\bm{N}}}
\newcommand{\bP}{{\bm{P}}}
\newcommand{\bQ}{{\bm{Q}}}
\newcommand{\bR}{{\bm{R}}}
\newcommand{\bS}{{\bm{S}}}
\newcommand{\bT}{{\bm{T}}}
\newcommand{\bW}{{\bm{W}}}
\newcommand{\bX}{{\bm{X}}}
\newcommand{\bY}{{\bm{Y}}}
\newcommand{\bIB}{{\bm{B}^{\rm irr}}}
\newcommand{\bOB}{{\bm{B}^{\rm ord}}}
\newcommand{\bOS}{{\bm{OS}}}
\newcommand{\bERR}{{\bm{ERR}}}
\newcommand{\bSP}{{\bm{SP}}}
\newcommand{\bMV}{{\bm{MV}}}
\newcommand{\bBM}{{\bm{BM}}}
\newcommand{\balpha}{{\bm{\alpha}}}
\newcommand{\bbeta}{{\bm{\beta}}}
\newcommand{\bgamma}{{\bm{\gamma}}}
\newcommand{\bdelta}{{\bm{\delta}}}
\newcommand{\bkappa}{{\bm{\kappa}}}
\newcommand{\bmu}{{\bm{\mu}}}
\newcommand{\bomega}{{\bm{\omega}}}
\newcommand{\bphi}{{\bm{\phi}}}
\newcommand{\bsigma}{{\bm{\sigma}}}
\newcommand{\btau}{{\bm{\tau}}}
\newcommand{\bpsi}{{\bm{\psi}}}
\newcommand{\bxi}{{\bm{\xi}}}
\newcommand{\bzeta}{{\bm{\zeta}}}
\newcommand{\bone}{{\bm{1}}}
\newcommand{\bzero}{{\bm{0}}}

\newcommand{\phihat}{{\widehat{{\rule{0ex}{1.45ex}\mkern-2mu\smash{\phi}}}}}

\newcommand{\Cbar}{{\overline{C}}}
\newcommand{\Dbar}{{\overline{D}}}
\newcommand{\dbar}{{\overline{d}}}
\def\Ctilde{{\widetilde{C}}}
\def\Ftilde{{\widetilde{F}}}
\def\Gtilde{{\widetilde{G}}}
\def\Htilde{{\widetilde{H}}}
\def\Ptilde{{\widetilde{P}}}
\def\Chat{{\widehat{C}}}
\def\ctilde{{\widetilde{c}}}
\def\zbar{{\overline{Z}}}
\def\pitilde{{\widetilde{\pi}}}

\newcommand{\sech}{{\rm sech}}

%
%
\newcommand{\sn}{{\rm sn}}
\newcommand{\cn}{{\rm cn}}
\newcommand{\dn}{{\rm dn}}
\newcommand{\sm}{{\rm sm}}
\newcommand{\cm}{{\rm cm}}

%
%
\newcommand{\zfz}{ {{}_0 \! F_0} }
\newcommand{\zfo}{ {{}_0  F_1} }
\newcommand{\ofz}{ {{}_1 \! F_0} }
\newcommand{\ofo}{ {{}_1 \! F_1} }
\newcommand{\oft}{ {{}_1 \! F_2} }

%
%
\newcommand{\FHyper}[2]{ {\tensor[_{#1 \!}]{F}{_{#2}}\!} }
\newcommand{\FHYPER}[5]{ {\FHyper{#1}{#2} \!\biggl(
   \!\!\begin{array}{c} #3 \\[1mm] #4 \end{array}\! \bigg|\, #5 \! \biggr)} }
\newcommand{\tfo}{ {\FHyper{2}{1}} }
\newcommand{\tfz}{ {\FHyper{2}{0}} }
\newcommand{\threefz}{ {\FHyper{3}{0}} }
\newcommand{\FHYPERbottomzero}[3]{ {\FHyper{#1}{0} \hspace*{-0mm}\biggl(
   \!\!\begin{array}{c} #2 \\[1mm] \hbox{---} \end{array}\! \bigg|\, #3 \! \biggr)} }
\newcommand{\FHYPERtopzero}[3]{ {\FHyper{0}{#1} \hspace*{-0mm}\biggl(
   \!\!\begin{array}{c} \hbox{---} \\[1mm] #2 \end{array}\! \bigg|\, #3 \! \biggr)} }

\newcommand{\phiHyper}[2]{ {\tensor[_{#1}]{\phi}{_{#2}}} }
\newcommand{\psiHyper}[2]{ {\tensor[_{#1}]{\psi}{_{#2}}} }
\newcommand{\PhiHyper}[2]{ {\tensor[_{#1}]{\Phi}{_{#2}}} }
\newcommand{\PsiHyper}[2]{ {\tensor[_{#1}]{\Psi}{_{#2}}} }
\newcommand{\phiHYPER}[6]{ {\phiHyper{#1}{#2} \!\left(
   \!\!\begin{array}{c} #3 \\ #4 \end{array}\! ;\, #5, \, #6 \! \right)\!} }
\newcommand{\psiHYPER}[6]{ {\psiHyper{#1}{#2} \!\left(
   \!\!\begin{array}{c} #3 \\ #4 \end{array}\! ;\, #5, \, #6 \! \right)} }
\newcommand{\PhiHYPER}[5]{ {\PhiHyper{#1}{#2} \!\left(
   \!\!\begin{array}{c} #3 \\ #4 \end{array}\! ;\, #5 \! \right)\!} }
\newcommand{\PsiHYPER}[5]{ {\PsiHyper{#1}{#2} \!\left(
   \!\!\begin{array}{c} #3 \\ #4 \end{array}\! ;\, #5 \! \right)\!} }
\newcommand{\zerophizero}{ {\phiHyper{0}{0}} }
\newcommand{\ophizero}{ {\phiHyper{1}{0}} }
\newcommand{\zphio}{ {\phiHyper{0}{1}} }
\newcommand{\ophio}{ {\phiHyper{1}{1}} }
\newcommand{\tphio}{ {\phiHyper{2}{1}} }
\newcommand{\tphiz}{ {\phiHyper{2}{0}} }
\newcommand{\tPhio}{ {\PhiHyper{2}{1}} }
\newcommand{\opsio}{ {\psiHyper{1}{1}} }

%
%
\newcommand{\stirlingsubset}[2]{\genfrac{\{}{\}}{0pt}{}{#1}{#2}}
\newcommand{\stirlingcycleold}[2]{\genfrac{[}{]}{0pt}{}{#1}{#2}}
\newcommand{\stirlingcycle}[2]{\left[\! \stirlingcycleold{#1}{#2} \!\right]}
\newcommand{\assocstirlingsubset}[3]{{\genfrac{\{}{\}}{0pt}{}{#1}{#2}}_{\! \ge #3}}
\newcommand{\genstirlingsubset}[4]{{\genfrac{\{}{\}}{0pt}{}{#1}{#2}}_{\! #3,#4}}
\newcommand{\irredstirlingsubset}[2]{{\genfrac{\{}{\}}{0pt}{}{#1}{#2}}^{\!\rm irr}}
\newcommand{\euler}[2]{\genfrac{\langle}{\rangle}{0pt}{}{#1}{#2}}
\newcommand{\eulergen}[3]{{\genfrac{\langle}{\rangle}{0pt}{}{#1}{#2}}_{\! #3}}
\newcommand{\eulersecond}[2]{\left\langle\!\! \euler{#1}{#2} \!\!\right\rangle}
\newcommand{\eulersecondgen}[3]{{\left\langle\!\! \euler{#1}{#2} \!\!\right\rangle}_{\! #3}}
\newcommand{\binomvert}[2]{\genfrac{\vert}{\vert}{0pt}{}{#1}{#2}}
\newcommand{\binomsquare}[2]{\genfrac{[}{]}{0pt}{}{#1}{#2}}
\newcommand{\qbinom}[3]{\genfrac{(}{)}{0pt}{}{#1}{#2}_{\!\! #3}}


\newenvironment{sarray}{
             \textfont0=\scriptfont0
             \scriptfont0=\scriptscriptfont0
             \textfont1=\scriptfont1
             \scriptfont1=\scriptscriptfont1
             \textfont2=\scriptfont2
             \scriptfont2=\scriptscriptfont2
             \textfont3=\scriptfont3
             \scriptfont3=\scriptscriptfont3
           \renewcommand{\arraystretch}{0.7}
           \begin{array}{l}}{\end{array}}

\newenvironment{scarray}{
             \textfont0=\scriptfont0
             \scriptfont0=\scriptscriptfont0
             \textfont1=\scriptfont1
             \scriptfont1=\scriptscriptfont1
             \textfont2=\scriptfont2
             \scriptfont2=\scriptscriptfont2
             \textfont3=\scriptfont3
             \scriptfont3=\scriptscriptfont3
           \renewcommand{\arraystretch}{0.7}
           \begin{array}{c}}{\end{array}}


\newcommand*\circled[1]{\tikz[baseline=(char.base)]{
  \node[shape=circle,draw,inner sep=1pt] (char) {#1};}}
\newcommand{\ostar}{{\circledast}}
\newcommand{\ostarN}{{\,\circledast_{\vphantom{\dot{N}}N}\,}}
\newcommand{\ostarPsi}{{\,\circledast_{\vphantom{\dot{\Psi}}\Psi}\,}}
\newcommand{\starN}{{\,\ast_{\vphantom{\dot{N}}N}\,}}
\newcommand{\starpsi}{{\,\ast_{\vphantom{\dot{\bpsi}}\!\bpsi}\,}}
\newcommand{\starone}{{\,\ast_{\vphantom{\dot{1}}1}\,}}
\newcommand{\startwo}{{\,\ast_{\vphantom{\dot{2}}2}\,}}
\newcommand{\starinfty}{{\,\ast_{\vphantom{\dot{\infty}}\infty}\,}}
\newcommand{\starT}{{\,\ast_{\vphantom{\dot{T}}T}\,}}

\newcommand*{\Scale}[2][4]{\scalebox{#1}{$#2$}}

\newcommand*{\Scaletext}[2][4]{\scalebox{#1}{#2}} 

\clearpage

\tableofcontents

\clearpage

\section{Introduction and statement of results}

It is well known\footnote{
   See e.g.\ \cite{Clarke_58}, \cite[pp.~26--27]{Moon_70},
   \cite[p.~70]{Comtet_74}, \cite[pp.~25--28]{Stanley_99}
   or \cite{Avron_16}.
   See also \cite{Riordan_68b,Francon_74,Sagan_83,Pitman_02,Guo_17}
   and \cite[pp.~235--240]{Aigner_18}
   for related information.
}
that the number of forests of rooted trees on $n$ labeled vertices
is $f_n = (n+1)^{n-1}$,
and that the number of forests of rooted trees on $n$ labeled vertices
having $k$ components (i.e.\ $k$ trees) is
\be
   f_{n,k} \;=\;  \binom{n-1}{k-1} \, n^{n-k} 
           \;=\;  \binom{n}{k} \, k \, n^{n-k-1}
 \label{def.fnk}
\ee
(to be interpreted as $\delta_{k0}$ when $n=0$).
In particular, the number of rooted trees on $n$ labeled vertices
is $f_{n,1} = n^{n-1}$.
The first few $f_{n,k}$ and $f_n$ are
%
%
\vspace*{-5mm}
\begin{table}[H]
\centering
\small
\begin{tabular}{c|rrrrrrrrr|r}
$n \setminus k$ & 0 & 1 & 2 & 3 & 4 & 5 & 6 & 7 & 8 & $f_n$ \\
\hline
0 & 1 &  &  &  &  &  &  &  &  & 1  \\
1 & 0 & 1 &  &  &  &  &  &  &  & 1  \\
2 & 0 & 2 & 1 &  &  &  &  &  &  & 3  \\
3 & 0 & 9 & 6 & 1 &  &  &  &  &  & 16  \\
4 & 0 & 64 & 48 & 12 & 1 &  &  &  &  & 125  \\
5 & 0 & 625 & 500 & 150 & 20 & 1 &  &  &  & 1296  \\
6 & 0 & 7776 & 6480 & 2160 & 360 & 30 & 1 &  &  & 16807  \\
7 & 0 & 117649 & 100842 & 36015 & 6860 & 735 & 42 & 1 &  & 262144  \\
8 & 0 & 2097152 & 1835008 & 688128 & 143360 & 17920 & 1344 & 56 & 1 & 4782969 \\
\end{tabular}
\end{table}
\vspace*{-5mm}

\noindent
\!\!\cite[A061356/A137452 and A000272]{OEIS}.
By adding a new vertex 0 and connecting it to the roots of all the trees,
we see that $f_n$ is also
the number of (unrooted) trees on $n+1$ labeled vertices, and that
$f_{n,k}$ is the number of (unrooted) trees on $n+1$ labeled vertices
in which some specified vertex (here vertex 0) has degree $k$.

The unit-lower-triangular matrix $(f_{n,k})_{n,k \ge 0}$
has the exponential generating function
\be
   \sum_{n=0}^\infty \sum_{k=0}^n f_{n,k} \, {t^n \over n!} \, x^k
   \;=\;
   e^{x T(t)}
   \;,
 \label{eq.fnk.egf}
\ee
where
\be
   T(t)  \;\eqdef\; \sum_{n=1}^\infty n^{n-1} \, {t^n \over n!}
 \label{def.treefn}
\ee
is the {\em tree function}\/ \cite{Corless_96}.\footnote{
   In the analysis literature, expressions involving the tree function
   are often written in terms of the {\em Lambert $W$ function}\/
   $W(t) = -T(-t)$, which is the inverse function to $w \mapsto w e^w$
   \cite{Corless_96,Kalugin_12b}.
}
An equivalent statement is that
the unit-lower-triangular matrix $(f_{n,k})_{n,k \ge 0}$
is \cite{Barry_10_OEIS}
the exponential Riordan array \cite{Deutsch_04,Deutsch_09,Barry_16}
$\scrr[F,G]$ with $F(t) = 1$ and $G(t) = T(t)$;
we will discuss this connection in Section~\ref{subsec.EGF.1}.

The principal purpose of this paper is to prove the total positivity
of some matrices related to (and generalizing) $f_n$ and $f_{n,k}$.
Recall first that a finite or infinite matrix of real numbers is called
{\em totally positive}\/ (TP) if all its minors are nonnegative,
and {\em strictly totally positive}\/ (STP)
if all its minors are strictly positive.\footnote{
   {\bf Warning:}  Many authors
   (e.g.\ \cite{Gantmakher_37,Gantmacher_02,Fomin_00,Fallat_11})
   use the terms ``totally nonnegative'' and ``totally positive''
   for what we have termed ``totally positive'' and
   ``strictly totally positive'', respectively.
   So it is very important, when seeing any claim about
   ``totally positive'' matrices, to ascertain which sense of
   ``totally positive'' is being used!
   (This is especially important because many theorems in this subject
    require {\em strict}\/ total positivity for their validity.)
}
Background information on totally positive matrices can be found
in \cite{Karlin_68,Gantmacher_02,Pinkus_10,Fallat_11};
they have applications to many areas of pure and applied mathematics.\footnote{
   Including combinatorics
   \cite{Brenti_89,Brenti_95,Brenti_96,Fomin_00,Skandera_03},
   stochastic processes \cite{Karlin_59,Karlin_68},
   statistics \cite{Karlin_68},
   the mechanics of oscillatory systems \cite{Gantmakher_37,Gantmacher_02},
   the zeros of polynomials and entire functions
   \cite{Karlin_68,Asner_70,Kemperman_82,Holtz_03,Pinkus_10,Dyachenko_14},
   spline interpolation \cite{Schoenberg_53,Karlin_68,Gasca_96},
   Lie theory \cite{Lusztig_94,Lusztig_98,Fomin_99,Lusztig_08}
   and cluster algebras \cite{Fomin_10,Fomin_forthcoming},
   the representation theory of the infinite symmetric group
   \cite{Thoma_64,Borodin_17},
   the theory of immanants \cite{Stembridge_91},
   planar discrete potential theory \cite{Curtis_98,Fomin_01}
   and the planar Ising model \cite{Lis_17},
   and several other areas \cite{Gasca_96}.
}

Our first result is the following:

\begin{theorem}
   \label{thm1.1}
\quad\hfill\vspace*{-1mm}
\begin{itemize}
   \item[(a)]  The unit-lower-triangular matrix $F = (f_{n,k})_{n,k \ge 0}$
       is totally positive.\footnote{
   I trust that there will be no confusion between my use of the letter $F$
   for the matrix $(f_{n,k})_{n,k \ge 0}$ or its generalizations,
   and also for the power series $F(t)$ in an exponential Riordan array
   $\scrr[F,G]$.
   The meaning should be unambiguous from the context.
}
   \item[(b)]  The Hankel matrix
       $(f_{n+n'+1,1})_{n,n' \ge 0}$
       is totally positive.
\end{itemize}
\end{theorem}

It is known \cite{Gantmakher_37,Pinkus_10}
that a Hankel matrix of real numbers is totally positive
if and only if the underlying sequence
is a Stieltjes moment sequence,
i.e.\ the moments of a positive measure on $[0,\infty)$.
And it is also known
that $(f_{n+1,1})_{n \ge 0} = ((n+1)^n)_{n \ge 0}$
is a Stieltjes moment sequence.\footnote{
   The integral representation
   \cite{Bouwkamp_86} \cite[Corollary~2.4]{Kalugin_12b}
   $$
      {(n+1)^n \over n!}
      \;=\;
      {1 \over \pi}
      \int\limits_0^\pi
          \biggl( {\sin\nu \over \nu} \, e^{\nu \cot\nu} \biggr) ^{\! n+1}
                 \: d\nu
   $$
   shows that $(n+1)^n/n!$ is a Stieltjes moment sequence.
   Moreover, $n! = \int_0^\infty x^n \, e^{-x} \, dx$
   is a Stieltjes moment sequence.
   Since the entrywise product of two Stieltjes moment sequences
   is easily seen to be a Stieltjes moment sequence,
   it follows that $(n+1)^n$ is a Stieltjes moment sequence.
   But I do not know any simple formula
   (i.e.\ one involving only a single integral over a real variable)
   for its Stieltjes integral representation.
 \label{footnote.stieltjes.fn}
}
So Theorem~\ref{thm1.1}(b) is equivalent to this known result.
But our proof here is combinatorial and linear-algebraic, not analytic.

However, this is only the beginning of the story,
because our main interest \cite{Sokal_flajolet,Sokal_OPSFA,Sokal_totalpos}
is not with sequences and matrices of real numbers,
but rather with sequences and matrices of {\em polynomials}\/
(with integer or real coefficients) in one or more indeterminates $\bfx$:
in applications they will typically be generating polynomials that enumerate
some combinatorial objects with respect to one or more statistics.
We equip the polynomial ring $\R[\bfx]$ with the coefficientwise
partial order:  that is, we say that $P$ is nonnegative
(and write $P \myge 0$)
in case $P$ is a polynomial with nonnegative coefficients.
We then say that a matrix with entries in $\R[\bfx]$ is
\textbfit{coefficientwise totally positive}
if all its minors are polynomials with nonnegative coefficients;
and we say that a sequence $\ba = (a_n)_{n \ge 0}$ with entries in $\R[\bfx]$
is \textbfit{coefficientwise Hankel-totally positive}
if its associated infinite Hankel matrix
$H_\infty(\ba) = (a_{n+n'})_{n,n' \ge 0}$
is coefficientwise totally positive.
Most generally, we can consider sequences and matrices
with entries in an arbitrary partially ordered commutative ring;
total positivity and Hankel-total positivity
are then defined in the obvious way
(see Section~\ref{subsec.totalpos.prelim}).
Coefficientwise Hankel-total positivity of a sequence of polynomials
$(P_n(\bfx))_{n \ge 0}$ {\em implies}\/ the pointwise Hankel-total positivity
(i.e.\ the Stieltjes moment property) for all $\bfx \ge 0$,
but it is vastly stronger.

Returning now to the matrix $F = (f_{n,k})_{n,k \ge 0}$,
let us define its \textbfit{row-generating polynomials} in the usual way:
\be
   F_n(x)  \;=\;  \sum_{\ell=0}^n f_{n,\ell} \, x^\ell
   \;.
 \label{def.Fn}
\ee
More generally, let us define its
\textbfit{binomial partial row-generating polynomials}
\be
   F_{n,k}(x)
   \;=\;
   \sum_{\ell=k}^n f_{n,\ell} \, \binom{\ell}{k} \, x^{\ell-k}
   \;.
 \label{def.Fnk}
\ee
Thus, $F_n(x)$ is the generating polynomial for
forests of rooted trees on $n$ labeled vertices,
with a weight $x$ for each component;
and $F_{n,k}(x)$ is the generating polynomial for
forests of rooted trees on $n$ labeled vertices
with $k$ distinguished components,
with a weight $x$ for each undistinguished component.
In fact, we have the explicit formulae
\begin{eqnarray}
   F_n(x)      & = &  x \, (x+n)^{n-1}
      \label{eq.Fn.explicit}   \\[2mm]
   F_{n,k}(x)  & = &  \binom{n}{k} \, (x+k) \, (x+n)^{n-k-1}
      \label{eq.Fnk.explicit}
\end{eqnarray}
as can easily be verified by expanding the right-hand sides.
The $F_n(x)$ are a specialization of
the celebrated {\em Abel polynomials}\/ $A_n(x;a) = x (x-an)^{n-1}$
\cite{Mullin_70,Francon_74,Sagan_83,Roman_84}
to $a=-1$,
while the $F_{n,k}(x)$ can be found in \cite{Pitman_02,Wang_11}.

{}From \reff{def.Fnk} we see that $F_{n,k}(x)$ is a polynomial of degree $n-k$
with nonnegative integer coefficients,
with leading coefficient $\binom{n}{k}$;
in particular, $F_{n,n}(x) = 1$.
Moreover, $F_{n,0}(x) = F_n(x)$ [because $\binom{\ell}{0} = 1$]
and $F_{n,k}(0) = f_{n,k}$ [because $\binom{k}{k} = 1$].
So the matrix $F(x) = \bigl( F_{n,k}(x) \bigr)_{n,k \ge 0}$
is a unit-lower-triangular matrix, with entries in $\Z[x]$,
that has the row-generating polynomials $F_n(x)$ in its zeroth column
and that reduces to $F = (f_{n,k})_{n,k \ge 0}$ when $x=0$.
Because of the presence of the binomial coefficients $\binom{\ell}{k}$
in \reff{def.Fnk}, we call $F(x)$
the \textbfit{binomial row-generating matrix} of the matrix $F$.\footnote{
   Let us remark that the ordinary row-generating matrix of
   a lower-triangular matrix ---
   that is, \reff{def.Fnk} without the factors $\binom{\ell}{k}$ ---
   has been introduced recently by several authors
   \cite{Chang_16,Mu_17a,Zhu_17a}.
   I do not know whether the binomial row-generating matrix
   has been used previously, but I suspect that it has been.
}
Please note that the definition \reff{def.Fnk}
can be written as a matrix factorization
\be
   F(x)  \;=\;  F B_x
   \;,
 \label{eq.FBx}
\ee
where $B_x$ is the weighted binomial matrix
\be
   (B_x)_{ij}  \;=\;  \binom{i}{j} \, x^{i-j}
 \label{def.Bx}
\ee
(note that it too is unit-lower-triangular);
this factorization will play a central role in our proofs.
Our second result is then:

\begin{theorem}
   \label{thm1.2}
\quad\hfill\vspace*{-1mm}
\begin{itemize}
   \item[(a)]  The unit-lower-triangular polynomial matrix
      $F(x) = \bigl( F_{n,k}(x) \bigr)_{n,k \ge 0}$
      is coefficientwise totally positive.
   \item[(b)]  The polynomial sequence
       $\bF = \bigl( F_{n}(x) \bigr)_{n \ge 0}$
       is coefficientwise Hankel-totally positive.
       [That is, the Hankel matrix
        $H_\infty(\bF) = \bigl( F_{n+n'}(x) \bigr)_{n,n' \ge 0}$
        is coefficientwise totally positive.]
\end{itemize}
\end{theorem}

It is not difficult to see (see Lemma~\ref{lemma.binomialmatrix.TP} below)
that the matrix $B_x$ is coefficientwise totally positive;
and it is an immediate consequence of the Cauchy--Binet formula
that the product of two (coefficientwise) totally positive matrices
is (coefficientwise) totally positive.
So Theorem~\ref{thm1.2}(a) is actually an immediate consequence
of Theorem~\ref{thm1.1}(a) together with \reff{eq.FBx}
and Lemma~\ref{lemma.binomialmatrix.TP}.
But Theorem~\ref{thm1.2}(b) will take more work.

But this is still not the end of the story, because we want to
generalize these polynomials further by adding further variables.
First let us agree that the vertices of our forest $\scrf$ of rooted trees
will henceforth be labeled by the totally ordered set $[n] = \{1,2,\ldots,n\}$.
Given a rooted tree $T \in \scrf$ and two vertices $i,j$ of $T$,
we say that $j$ is a {\em descendant}\/ of $i$
if the unique path from the root of $T$ to $j$ passes through $i$.
(Note in particular that every vertex is a descendant of itself.)
Now let $e = ij$ be an edge of $T$, ordered so that $j$ is a descendant of $i$;
then $i$ is the {\em parent}\/ of $j$, and $j$ is a {\em child}\/ of $i$.
We say that the edge $e = ij$ is \textbfit{improper}
if there exists a descendant of $j$ (possibly $j$ itself)
that is lower-numbered than $i$;
otherwise we say that $e = ij$ is \textbfit{proper}.

Now let $f_{n,k,m}$
be the number of forests of rooted trees on the vertex set $[n]$
that have $k$ components and $m$ improper edges
(note that $0 \le m \le n-k$ since a forest with $k$ components
 has $n-k$ edges).
And introduce the generating polynomial
that gives a weight $y$ for each improper edge
and a weight $z$ for each proper edge:
\be
   f_{n,k}(y,z)
   \;=\;
   \sum_{m=0}^{n-k} f_{n,k,m} \, y^m z^{n-k-m}
   \;.
 \label{def.fnkyz}
\ee
The first few $f_{n,k}(y,z)$ are
%
%
\vspace*{-5mm}
\begin{table}[H]
\centering
\small
\begin{tabular}{c|lllll}
$n \setminus k$ & 0 & 1 & 2 & 3 & 4 \\
\hline
0 & 1 &  &  &  &  \\
1 & 0 & 1 &  &  &  \\
2 & 0 & $z + y$ & 1 &  &  \\
3 & 0 & $2z^2 + 4zy + 3y^2$ & $3z + 3y$ & 1 &  \\
4 & 0 & $6z^3 + 18z^2y + 25zy^2 + 15y^3$ & $11z^2 + 22zy + 15y^2$ & $6z + 6y$ & 1  \\
\end{tabular}
\end{table}
\vspace{-5mm}

\noindent
Clearly $f_{n,k}(y,z)$ is a homogeneous polynomial of degree $n-k$
with nonnegative integer coefficients;
it is a polynomial refinement of $f_{n,k}$
in the sense that $f_{n,k}(1,1) = f_{n,k}$.
(Of course, it was redundant to introduce the two variables $y$ and $z$
 instead of just one of them; we did it because it makes the formulae
 more symmetric.)
In particular, the polynomials $f_{n,1}(y,z)$ enumerate rooted trees
according to the number of improper edges;
they are homogenized versions of the celebrated {\em Ramanujan polynomials}\/
\cite{Shor_95,Dumont_96,Zeng_99,Guo_07,Lin_14,Josuat-Verges_15,Guo_18,Chen_21,Randazzo_21}
\cite[A054589]{OEIS}.\footnote{
   Also, Drake \cite[Example~1.7.3]{Drake_08} shows that $f_{n,1,n-1-\ell}$
   is the number of {\em ordered}\/ rooted trees
   (i.e.\ trees in which the children of each vertex are linearly ordered)
   on the vertex set $[n]$ with $\ell$ ascents
   (i.e.\ edges $ij$ in which the child $j$ is higher-numbered than the
   parent $i$)
   in which the child vertex of every ascent is a leaf.
}

The unit-lower-triangular matrix $(f_{n,k}(y,z))_{n,k \ge 0}$
is also the exponential Riordan array $\scrr[F,G]$ with $F(t) = 1$ and
\be
   G(t)
   \;=\;
   {1 \over z}
   \biggl[ T\Big( \Big(1 - {z \over y} + {z^2 \over y} t \Big) \:
                  e^{- \, \displaystyle \big(1 - {z \over y} \big)}
            \Big)
            \:-\: \Big(1 - {z \over y} \Big)
   \biggr]
   \;,
\ee
where $T(t)$ is the tree function \reff{def.treefn};
we will show this in Section~\ref{subsec.EGF.2}.

\bigskip

{\bf Remark.}
Let us write the homogenized Ramanujan polynomials as
$f_{n+1,1}(y,z) = \sum\limits_{m=0}^n r(n,m) \, y^m z^{n-m}$,
so that $r(n,m) = f_{n+1,1,m}$ is the number of rooted trees
on the vertex set $[n+1]$ with $m$ improper edges.
Then Shor \cite{Shor_95} and Dumont--Ramamonjisoa \cite{Dumont_96} showed that
\be
   r(n,m)  \;=\;  n \, r(n-1,m) \:+\: (n+m-1) \, r(n-1,m-1)
   \;.
\ee
That is, the Ramanujan polynomials are the row-generating polynomials
for the Graham--Knuth--Patashnik (GKP) recurrence
\cite[Problem~6.94, pp.~319 and~564]{Graham_94}
\cite{Theoret_94,Neuwirth_01,Spivey_11,Barbero_14,GKP_recurrence}
\begin{equation}
  T(n,m)
  \;=\;
  (\alpha n + \beta m + \gamma)    \, T(n-1,m)
  \:+\:
  (\alpha' n + \beta' m + \gamma') \, T(n-1,m-1)
  \quad\hbox{for $n \ge 1$}
\end{equation}
with initial condition $T(0,m) = \delta_{m0}$,
specialized to
$(\alpha,\beta,\gamma, \alpha',\beta',\gamma') = (1,0,0,1,1,-1)$.
It is an interesting open problem to extend the results presented here
to other cases of the GKP recurrence,
such as the generalized Ramanujan polynomials considered in
\cite{Guo_07,Lin_14,Randazzo_21}.
Concerning total positivity for some special cases of the GKP recurrence,
see \cite[especially Conjecture~1.3 and Theorem~1.5]{FPSAC2021};
and concerning coefficientwise Hankel-total positivity for
the row-generating polynomials of the GKP recurrence,
see \cite[Conjectures~6.9 and 6.10]{GKP_recurrence}.

It should, however, be remarked that the coefficient matrix
of the Ramanujan polynomials, $R = (r(n,m))_{n,m \ge 0}$,
is {\em not}\/ totally positive:
the lower-left $7 \times 7$ minor of the leading $9 \times 9$ matrix
is $-3709251874944000$.
\myendremark

\medskip

Now we can again introduce row-generating polynomials
and binomial partial row-generating polynomials:
we generalize \reff{def.Fn} and \reff{def.Fnk} by defining
\be
   F_n(x,y,z)  \;=\;  \sum_{\ell=0}^n f_{n,\ell}(y,z) \, x^\ell
 \label{def.Fnyz}
\ee
and
\be
   F_{n,k}(x,y,z)
   \;=\;
   \sum_{\ell=k}^n f_{n,\ell}(y,z) \, \binom{\ell}{k} \, x^{\ell-k}
   \;.
 \label{def.Fnkyz}
\ee
Thus, $F_n(x,y,z)$ is the generating polynomial for
forests of rooted trees on the vertex set $[n]$,
with a weight $x$ for each component
and a weight $y$ (resp.~$z$) for each improper (resp.\ proper) edge;
and $F_{n,k}(x,y,z)$ is the generating polynomial for
forests of rooted trees on the vertex set $[n]$
with $k$ distinguished components,
with a weight $x$ for each undistinguished component
and a weight $y$ (resp.~$z$) for each improper (resp.\ proper) edge.
Note that $F_n(x,y,z)$ [resp.\ $F_{n,k}(x,y,z)$]
is a homogeneous polynomial of degree $n$ (resp.~$n-k$) in $x,y,z$.
Our third result is then:

\begin{theorem}
   \label{thm1.3}
\quad\hfill\vspace*{-1mm}
\begin{itemize}
   \item[(a)]  The unit-lower-triangular polynomial matrix
      $F(x,y,z) = \bigl( F_{n,k}(x,y,z) \bigr)_{n,k \ge 0}$
      is coefficientwise totally positive (jointly in $x,y,z$).
   \item[(b)]  The polynomial sequence
       $\bF = \bigl( F_{n}(x,y,z) \bigr)_{n \ge 0}$
       is coefficientwise Hankel-totally positive (jointly in $x,y,z$).
   \item[(c)] The polynomial sequence
       $\bF^\triangle = \bigl( f_{n+1,1}(y,z) \bigr)_{n \ge 0}$
       is coefficientwise Hankel-totally positive (jointly in $y,z$).
\end{itemize}
\end{theorem}

\noindent
Here part~(c) is an easy consequence of part~(b),
obtained by restricting to $n \ge 1$, dividing by~$x$, and taking $x \to 0$.


We remark that Chen {\em et al.}\/ \cite[Corollary~3.3]{Chen_11}
have proven that the sequence $\bigl( f_{n+1,1}(y,z) \bigr)_{n \ge 0}$
of Ramanujan polynomials is coefficientwise strongly log-convex
(i.e.\ coefficientwise Hankel-totally positive of order~2).\footnote{
   See also Lin and Zeng \cite[Theorem~1.5]{Lin_14}
   for a generalization to some polynomials with additional indeterminates.
}
Theorem~\ref{thm1.3}(c) is thus an extension of this result
to prove coefficientwise Hankel-total positivity of all orders.

But this is {\em still}\/ not the end of the story,
because we can add even more variables --- in fact, an infinite set.
Given a rooted tree $T$ on a totally ordered vertex set
and vertices $i,j \in T$ such that $j$ is a child of $i$,
we say that $j$ is a \textbfit{proper child} of $i$
if the edge $e = ij$ is proper
(that~is, $j$ and all its descendants are higher-numbered than~$i$).
Now let $\bphi = (\phi_m)_{m \ge 0}$ be indeterminates,
and let $f_{n,k}(y,\bphi)$ be the generating polynomial for
$k$-component forests of rooted trees on the vertex set $[n]$
with a weight $\phihat_m \eqdef m! \, \phi_m$
for each vertex with $m$ proper children
and a weight $y$ for each improper edge.
(We will see later why it is convenient to introduce the factors~$m!$
in this definition.
Observe also that the variables~$z$ are now redundant,
because they would simply scale $\phi_m \to z^m \phi_m$.)
We call the polynomials $f_{n,k}(y,\bphi)$
the \textbfit{generic rooted-forest polynomials}.
Here $\bphi = (\phi_m)_{m \ge 0}$ are in the first instance indeterminates,
so that $f_{n,k}(y,\bphi)$ belongs to the polynomial ring $\Z[y,\bphi]$;
but we can then, if we wish, substitute specific values for $\bphi$
in any commutative ring $R$, leading to values $f_{n,k}(y,\bphi) \in R[y]$.
(Similar substitutions can of course also be made for $y$.)
When doing this we will use the same notation $f_{n,k}(y,\bphi)$,
as the desired interpretation for $\bphi$ should be clear from the context.

The polynomial $f_{n,k}(y,\bphi)$ is quasi-homogeneous of degree $n-k$
when $\phi_m$ is assigned weight $m$ and $y$ is assigned weight 1.
It follows from this quasi-homogeneity that
the variable $y$ is now in principle redundant,
since it can be absorbed into $\bphi$:
namely, if we define a rescaled $\bphi$ by
\be
   (\bphi^c)_m  \;\eqdef\; c^m \, \phi_m
   \;,
\ee
then
\be
   f_{n,k}(y,\bphi)   \;=\;  y^{n-k} \, f_{n,k}(1,\bphi^{1/y})
   \;.
 \label{eq.fnk.rescaling}
\ee
However, we prefer to retain the redundant variable $y$,
in order to avoid the division by $y$ inherent in \reff{eq.fnk.rescaling};
in particular, this facilitates the study of the limiting case $y=0$.

The lower-triangular matrix $(f_{n,k}(y,\bphi))_{n,k \ge 0}$
is also an exponential Riordan array $\scrr[F,G]$ with $F(t) = 1$,
as we will show in Section~\ref{subsec.EGF.3};
but this time the function $G(t)$ is rather more complicated.

Now define the row-generating polynomials
\be
   F_n(x,y,\bphi)  \;=\;  \sum_{\ell=0}^n f_{n,\ell}(y,\bphi) \, x^\ell
 \label{def.Fnyz.phi}
\ee
and
\be
   F_{n,k}(x,y,\bphi)
   \;=\;
   \sum_{\ell=k}^n f_{n,\ell}(y,\bphi) \, \binom{\ell}{k} \, x^{\ell-k}
   \;.
 \label{def.Fnkyz.phi}
\ee
Thus, $F_n(x,y,\bphi)$ is the generating polynomial for
forests of rooted trees on the vertex set $[n]$,
with a weight $x$ for each component,
$y$ for each improper edge,
and $m! \, \phi_m$ for each vertex with $m$~proper children;
and $F_{n,k}(x,y,\bphi)$ is the generating polynomial for
forests of rooted trees on the vertex set $[n]$
with $k$ distinguished components,
with a weight $x$ for each undistinguished component,
$y$ for each improper edge,
and $m! \, \phi_m$ for each vertex with $m$~proper children.
Our fundamental result is then the following:

\begin{theorem}
   \label{thm1.4}
Fix $1 \le r \le \infty$.
Let $R$ be a partially ordered commutative ring,
and let $\bphi = (\phi_m)_{m \ge 0}$ be a sequence in $R$
that is Toeplitz-totally positive of order~$r$.
Then:
\begin{itemize}
      \item[(a)]  The lower-triangular polynomial matrix
      $F(x,y,\bphi) = \bigl( F_{n,k}(x,y,\bphi) \bigr)_{n,k \ge 0}$
      is coefficientwise totally positive of order~$r$ (jointly in $x,y$).
   \item[(b)]  The polynomial sequence
       $\bF = \bigl( F_{n}(x,y,\bphi) \bigr)_{n \ge 0}$
       is coefficientwise Hankel-totally positive of order~$r$
       (jointly in $x,y$).
   \item[(c)]  The polynomial sequence
       $\bF^\triangle = \bigl( f_{n+1,1}(y,\bphi) \bigr)_{n \ge 0}$
       is coefficientwise Hankel-totally positive of order~$r$
       (in $y$).
\end{itemize}
\end{theorem}

\noindent
(The concept of Toeplitz-total positivity
 in a partially ordered commutative ring
 will be explained in detail in Section~\ref{subsec.totalpos.prelim}.
 Total positivity of order~$r$ means that the minors of size $\le r$
 are nonnegative.)
Here (a) and (b) are once again the key results;
(c)~is an easy consequence of~(b),
obtained by restricting to $n \ge 1$, dividing by~$x$, and taking $x \to 0$.
Specializing Theorem~\ref{thm1.4} to $r = \infty$, $R = \Q$ and $\phi_m = z^m/m!$
(which is indeed Toeplitz-totally positive: see \reff{eq.thm.aissen} below),
we recover Theorem~\ref{thm1.3}.

Theorem~\ref{thm1.4} generalizes the main result of our recent paper
\cite{latpath_lah} on the generic Lah polynomials,
to which it reduces when $y=0$;
we will explain this connection in Section~\ref{sec.lah}.


The main tool in our proofs is the theory of production matrices
\cite{Deutsch_05,Deutsch_09}
as applied to total positivity \cite{Sokal_totalpos},
combined with the theory of exponential Riordan arrays
\cite{Deutsch_04,Deutsch_09,Barry_16}.
Therefore, in Section~\ref{sec.prelim} we review some facts
about total positivity, production matrices, and exponential Riordan arrays
that will play a central role in our arguments.
This development culminates in Theorem~\ref{thm.hankelTP.expriordan};
it is the fundamental theoretical result that underlies all our proofs.
In Section~\ref{sec.EGF} we show that the matrices
$(f_{n,k})_{n,k \ge 0}$, $(f_{n,k}(y,z))_{n,k \ge 0}$
and $(f_{n,k}(y,\bphi))_{n,k \ge 0}$
are exponential Riordan arrays $\scrr[F,G]$ with $F=1$,
and we compute their generating functions $G$.
In Section~\ref{sec.proofs} we prove Theorems~\ref{thm1.1}--\ref{thm1.4},
by exhibiting the production matrices
for $F$, $F(x)$, $F(x,y,z)$ and $F(x,y,\bphi)$
and proving that these production matrices are
coefficientwise totally positive.
In Section~\ref{sec.lah} we discuss the connection
with the generic Lah polynomials that were introduced in \cite{latpath_lah}.
Finally, in Section~\ref{sec.open} we pose some open problems.

A sequel devoted to a different (but closely related)
class of polynomials enumerating rooted labeled trees,
written in collaboration with Xi Chen,
will appear elsewhere \cite{Chen-Sokal_trees_totalpos}.

{\bf Note added:}
Some related ideas concerning total positivity and exponential Riordan arrays
can be found in a recent paper of Zhu \cite{Zhu_21}.

\section{Preliminaries}   \label{sec.prelim}

Here we review some definitions and results from
\cite{Sokal_totalpos,latpath_lah}
that will be needed in the sequel.
We also include a brief review of
exponential Riordan arrays \cite{Deutsch_04,Deutsch_09,Barry_16}
and Lagrange inversion \cite{Gessel_16}.
The key result in this section
--- obtained by straightforward combination of the others ---
is Theorem~\ref{thm.hankelTP.expriordan}.

\subsection{Partially ordered commutative rings and total positivity}
   \label{subsec.totalpos.prelim}

In this paper all rings will be assumed to have an identity element 1
and to be nontrivial ($1 \ne 0$).

A \textbfit{partially ordered commutative ring} is a pair $(R,\scrp)$ where
$R$ is a commutative ring and $\scrp$ is a subset of $R$ satisfying
\begin{itemize}
   \item[(a)]  $0,1 \in \scrp$.
   \item[(b)]  If $a,b \in \scrp$, then $a+b \in \scrp$ and $ab \in \scrp$.
   \item[(c)]  $\scrp \cap (-\scrp) = \{0\}$.
\end{itemize}
We call $\scrp$ the {\em nonnegative elements}\/ of $R$,
and we define a partial order on $R$ (compatible with the ring structure)
by writing $a \le b$ as a synonym for $b-a \in \scrp$.
Please note that, unlike the practice in real algebraic geometry
\cite{Brumfiel_79,Lam_84,Prestel_01,Marshall_08},
we do {\em not}\/ assume here that squares are nonnegative;
indeed, this property fails completely for our prototypical example,
the ring of polynomials with the coefficientwise order,
since $(1-x)^2 = 1-2x+x^2 \not\myge 0$.

Now let $(R,\scrp)$ be a partially ordered commutative ring
and let $\bfx = \{x_i\}_{i \in I}$ be a collection of indeterminates.
In the polynomial ring $R[\bfx]$ and the formal-power-series ring $R[[\bfx]]$,
let $\scrp[\bfx]$ and $\scrp[[\bfx]]$ be the subsets
consisting of polynomials (resp.\ series) with nonnegative coefficients.
Then $(R[\bfx],\scrp[\bfx])$ and $(R[[\bfx]],\scrp[[\bfx]])$
are partially ordered commutative rings;
we refer to this as the \textbfit{coefficientwise order}
on $R[\bfx]$ and $R[[\bfx]]$.

A (finite or infinite) matrix with entries in a
partially ordered commutative ring
is called \textbfit{totally positive} (TP) if all its minors are nonnegative;
it is called \textbfit{totally positive of order~$\bm{r}$} (TP${}_r$)
if all its minors of size $\le r$ are nonnegative.
It follows immediately from the Cauchy--Binet formula that
the product of two TP (resp.\ TP${}_r$) matrices is TP
(resp.\ TP${}_r$).\footnote{
   For infinite matrices, we need some condition to ensure that
   the product is well-defined.
   For instance, the product $AB$ is well-defined whenever
   $A$ is row-finite (i.e.\ has only finitely many nonzero entries in each row)
   or $B$ is column-finite.
}
This fact is so fundamental to the theory of total positivity
that we shall henceforth use it without comment.

We say that a sequence $\ba = (a_n)_{n \ge 0}$
with entries in a partially ordered commutative ring
is \textbfit{Hankel-totally positive} 
(resp.\ \textbfit{Hankel-totally positive of order~$\bm{r}$})
if its associated infinite Hankel matrix
$H_\infty(\ba) = (a_{i+j})_{i,j \ge 0}$
is TP (resp.\ TP${}_r$).
We say that $\ba$
is \textbfit{Toeplitz-totally positive} 
(resp.\ \textbfit{Toeplitz-totally positive of order~$\bm{r}$})
if its associated infinite Toeplitz matrix
$T_\infty(\ba) = (a_{i-j})_{i,j \ge 0}$
(where $a_n \eqdef 0$ for $n < 0$)
is TP (resp.\ TP${}_r$).\footnote{
   When $R = \R$, Toeplitz-totally positive sequences are traditionally called
   {\em P\'olya frequency sequences}\/ (PF),
   and Toeplitz-totally positive sequences of order $r$
   are called {\em P\'olya frequency sequences of order $r$}\/ (PF${}_r$).
   See \cite[chapter~8]{Karlin_68} for a detailed treatment.
}

When $R = \R$, Hankel- and Toeplitz-total positivity have simple
analytic characterizations.
A sequence $(a_n)_{n \ge 0}$ of real numbers
is Hankel-totally positive if and only if it is a Stieltjes moment sequence
\cite[Th\'eor\`eme~9]{Gantmakher_37} \cite[section~4.6]{Pinkus_10}.
And a sequence $(a_n)_{n \ge 0}$ of real numbers 
is Toeplitz-totally positive if and only if its ordinary generating function
can be written as
\be
   \sum_{n=0}^\infty a_n t^n
   \;=\;
   C e^{\gamma t} t^m \prod_{i=1}^\infty {1 + \alpha_i t  \over  1 - \beta_i t}
 \label{eq.thm.aissen}
\ee
with $m \in \N$, $C,\gamma,\alpha_i,\beta_i \ge 0$,
$\sum \alpha_i < \infty$ and $\sum \beta_i < \infty$:
this is the celebrated Aissen--Schoenberg--Whitney--Edrei theorem
\cite[Theorem~5.3, p.~412]{Karlin_68}.
However, in a general partially ordered commutative ring $R$,
the concepts of Hankel- and Toeplitz-total positivity are more subtle.

We will need a few easy facts about the total positivity of special matrices:

\begin{lemma}[Bidiagonal matrices]
  \label{lemma.bidiagonal}
Let $A$ be a matrix with entries in a partially ordered commutative ring,
with the property that all its nonzero entries belong to two consecutive
diagonals.
Then $A$ is totally positive if and only if all its entries are nonnegative.
\end{lemma}

\proof
The nonnegativity of the entries (i.e.\ TP${}_1$)
is obviously a necessary condition for TP.
Conversely, for a matrix of this type it is easy to see that
every nonzero minor is simply a product of some entries.
\qed

\begin{lemma}[Toeplitz matrix of powers]
   \label{lemma.toeplitz.power}
Let $R$ be a partially ordered commutative ring, let $x \in R$,
and consider the infinite Toeplitz matrix
\be
   T_x
   \;\eqdef\;
   T_\infty(x^\N)
   \;=\;
   \begin{bmatrix}
      1   &     &     &     &         \\
      x   &  1  &     &     &            \\
      x^2 &  x  &  1  &     &            \\
      x^3 &  x^2  &  x  &  1    &            \\
      \vdots  & \vdots   & \vdots  & \vdots   & \ddots
   \end{bmatrix}
   \;.
 \label{def.Tx}
\ee
Then every minor of $T_x$ is either zero or else a power of $x$.
Hence $T_x$ is TP $\iff$ $T_x$ is TP${}_1$ $\iff$ $x \ge 0$.

In particular, if $x$ is an indeterminate, then $T_x$ 
is totally positive in the ring $\Z[x]$ equipped with the
coefficientwise order.
\end{lemma}

\proof
Consider a submatrix $A = (T_x)_{IJ}$
with rows $I = \{i_1 < \ldots < i_k \}$
and columns $J = \{j_1 < \ldots < j_k \}$.
We will prove by induction on $k$ that
$\det A$ is either zero or a power of $x$.
It is trivial if $k=0$ or 1.
If $A_{12} = A_{22} = 0$,
then $A_{1s} = A_{2s} = 0$ for all $s \ge 2$ by definition of $T_x$,
and $\det A = 0$.
If $A_{12}$ and $A_{22}$ are both nonzero,
then the first column of $A$ is $x^{j_2 - j_1}$ times the second column,
and again $\det A = 0$.
Finally, if $A_{12} = 0$ and $A_{22} \ne 0$
(by definition of $T_x$ this is the only other possibility),
then $A_{1s} = 0$ for all $s \ge 2$;
we then replace the first column of $A$ by
the first column minus $x^{j_2 - j_1}$ times the second column,
so that the new first column has $x^{i_1-j_1}$ in its first entry
(or zero if $i_1 < j_1$) and zeroes elsewhere.
Then $\det A$ equals $x^{i_1-j_1}$ (or zero if $i_1 < j_1$)
times the determinant of its last $k-1$ rows and columns,
so the claim follows from the inductive hypothesis.
\qed

\noindent
See also Example~\ref{exam.toeplitz.power.TP} below
for a second proof of the total positivity of $T_x$,
using production matrices.

\begin{lemma}[Binomial matrix]
   \label{lemma.binomialmatrix.TP}
In the ring $\Z$, the binomial matrix
$B = \bigl( {\textstyle \binom{n}{k}} \bigr) _{n,k \ge 0}$
is totally positive.
More generally, the weighted binomial matrix
$B_{x,y}  = \bigl( x^{n-k} y^k {\textstyle \binom{n}{k}} \bigr) _{\! n,k \ge 0}$
is totally positive in the ring $\Z[x,y]$ equipped with the
coefficientwise order.
\end{lemma}

\proof
It is well known that the binomial matrix $B$ is totally positive,
and this can be proven by a variety of methods:
e.g.\ using production matrices
\cite[pp.~136--137, Example~6.1]{Karlin_68}
\cite[pp.~108--109]{Pinkus_10},
by diagonal similarity to a totally positive Toeplitz matrix
\cite[p.~109]{Pinkus_10},
by exponentiation of a nonnegative lower-subdiagonal matrix
\cite[p.~63]{Fallat_11},
or by an application of the Lindstr\"om--Gessel--Viennot lemma
\cite[p.~24]{Fomin_00}.

Then $B_{x,y} = D B D'$ where $D = \diag\bigl( (x^n)_{n \ge 0} \bigr)$
and $D' = \diag\bigl( (x^{-k} y^k)_{k \ge 0} \bigr)$.
By Cauchy--Binet, $B_{x,y}$ is totally positive in the ring $\Z[x,x^{-1},y]$
equipped with the coefficientwise order.
But because $B$ is lower-triangular,
the elements of $B_{x,y}$ actually lie in the subring $\Z[x,y]$.
\qed

\noindent
See also Example~\ref{exam.binomial.matrix.TP} below
for an {\em ab initio}\/ proof of Lemma~\ref{lemma.binomialmatrix.TP}
using production matrices.

Finally, let us show that the sufficiency half
of the Aissen--Schoenberg--Whitney--Edrei theorem holds
(with a slight modification to avoid infinite products)
in a general partially ordered commutative ring.
We give two versions, depending on whether or~not
it is assumed that the ring $R$ contains the rationals:

\begin{lemma}[Sufficient condition for Toeplitz-total positivity]
   \label{lemma.toeplitz.1}
Let $R$ be a partially ordered commutative ring,
let $N$ be a nonnegative integer,
and let $\alpha_1,\ldots,\alpha_N$, $\beta_1,\ldots,\beta_N$ and $C$
be nonnegative elements in $R$.
Define the sequence $\ba = (a_n)_{n \ge 0}$ in $R$ by
\be
   \sum_{n=0}^\infty a_n t^n
   \;=\;
   C \, \prod_{i=1}^N {1 + \alpha_i t  \over  1 - \beta_i t}
   \;.
 \label{eq.thm.toeplitz.1}
\ee
Then the Toeplitz matrix $T_\infty(\ba)$ is totally positive.
\end{lemma}

\noindent
Of course, it is no loss of generality to have the same number $N$
of alphas and betas,
since some of the $\alpha_i$ or $\beta_i$ could be zero.

\begin{lemma}[Sufficient condition for Toeplitz-total positivity, with rationals]
   \label{lemma.toeplitz.1a}
Let $R$ be a partially ordered commutative ring containing the rationals,
let $N$ be a nonnegative integer,
and let $\alpha_1,\ldots,\alpha_N$, $\beta_1,\ldots,\beta_N$, $\gamma$ and $C$
be nonnegative elements in $R$.
Define the sequence $\ba = (a_n)_{n \ge 0}$ in $R$ by
\be
   \sum_{n=0}^\infty a_n t^n
   \;=\;
   C \, e^{\gamma t} \, \prod_{i=1}^N {1 + \alpha_i t  \over  1 - \beta_i t}
   \;.
 \label{eq.thm.toeplitz.1a}
\ee
Then the Toeplitz matrix $T_\infty(\ba)$ is totally positive.
\end{lemma}

\proofof{Lemma~\ref{lemma.toeplitz.1}}
We make a series of elementary observations:

1) The sequence $\ba = (1,\alpha,0,0,0,\ldots)$, corresponding
to the generating function $A(t) = 1 + \alpha t$,
is Toeplitz-totally positive if and only if $\alpha \ge 0$.
The ``only if'' is trivial,
and the ``if'' follows from Lemma~\ref{lemma.bidiagonal}
because the Toeplitz matrix $T_\infty(\ba)$ is bidiagonal.

2) The sequence $\ba = (1,\beta,\beta^2,\beta^3,\ldots)$,
corresponding to the generating function $A(t) = 1/(1 - \beta t)$,
is Toeplitz-totally positive if and only if $\beta \ge 0$.
The ``only if'' is again trivial,
and the ``if'' follows from Lemma~\ref{lemma.toeplitz.power}.

3) If $\ba$ and $\bb$ are sequences
with ordinary generating functions $A(t)$ and $B(t)$,
then the convolution $\bc = \ba * \bb$,
defined by $c_n = \sum_{k=0}^n a_k b_{n-k}$,
has ordinary generating function $C(t) = A(t) \, B(t)$;
moreover, the Toeplitz matrix $T_\infty(\bc)$ is simply
the matrix product $T_\infty(\ba) \, T_\infty(\bb)$.
It thus follows from the Cauchy--Binet formula that if
$\ba$ and~$\bb$ are Toeplitz-totally positive, then so is $\bc$.

4) A Toeplitz-totally positive sequence can be multiplied
by a nonnegative constant $C$, and it is still Toeplitz-totally positive.

Combining these observations proves the lemma.
\qed

\proofof{Lemma~\ref{lemma.toeplitz.1a}}
We add to the proof of Lemma~\ref{lemma.toeplitz.1}
the following additional observation:

5) The sequence $\ba = (\gamma^n/n!)_{n \ge 0}$,
corresponding to the generating function $A(t) = e^{\gamma t}$,
is Toeplitz-totally positive if and only if $\gamma \ge 0$.
The ``only if'' is again trivial,
and the ``if'' follows from Lemma~\ref{lemma.binomialmatrix.TP}
because $\gamma^{n-k}/(n-k)! = \binom{n}{k} \gamma^{n-k} \times k!/n!$
and hence $T_\infty(\ba) = D^{-1} B_{\gamma,1} D$
where $D = \diag(\, (n!)_{n \ge 0})$.
\qed

\subsection{Production matrices}   \label{subsec.production}

The method of production matrices \cite{Deutsch_05,Deutsch_09}
has become in recent years an important tool in enumerative combinatorics.
In the special case of a tridiagonal production matrix,
this construction goes back to Stieltjes' \cite{Stieltjes_1889,Stieltjes_1894}
work on continued fractions:
the production matrix of a classical S-fraction or J-fraction is tridiagonal.
In~the present paper, by contrast,
we shall need production matrices that are lower-Hessenberg
(i.e.\ vanish above the first superdiagonal)
but are not in general tridiagonal.
We therefore begin by reviewing briefly
the basic theory of production matrices.
The important connection of production matrices with total positivity
will be treated in the next subsection.

Let $P = (p_{ij})_{i,j \ge 0}$ be an infinite matrix
with entries in a commutative ring $R$.
In~order that powers of $P$ be well-defined,
we shall assume that $P$ is either row-finite
(i.e.\ has only finitely many nonzero entries in each row)
or column-finite.

Let us now define an infinite matrix $A = (a_{nk})_{n,k \ge 0}$ by
\be
   a_{nk}  \;=\;  (P^n)_{0k}
 \label{def.iteration}
\ee
(in particular, $a_{0k} = \delta_{0k}$).
Writing out the matrix multiplications explicitly, we have
\be
   a_{nk}
   \;=\;
   \sum_{i_1,\ldots,i_{n-1}}
      p_{0 i_1} \, p_{i_1 i_2} \, p_{i_2 i_3} \,\cdots\,
        p_{i_{n-2} i_{n-1}} \, p_{i_{n-1} k}
   \;,
 \label{def.iteration.walk}
\ee
so that $a_{nk}$ is the total weight for all $n$-step walks in $\N$
from $i_0 = 0$ to $i_n = k$, in~which the weight of a walk is the
product of the weights of its steps, and a step from $i$ to $j$
gets a weight $p_{ij}$.
Yet another equivalent formulation is to define the entries $a_{nk}$
by the recurrence
\be
   a_{nk}  \;=\;  \sum_{i=0}^\infty a_{n-1,i} \, p_{ik}
   \qquad\hbox{for $n \ge 1$}
 \label{def.iteration.bis}
\ee
with the initial condition $a_{0k} = \delta_{0k}$.

We call $P$ the \textbfit{production matrix}
and $A$ the \textbfit{output matrix},
and we write $A = \scro(P)$.
Note that if $P$ is row-finite, then so is $\scro(P)$;
if $P$ is lower-Hessenberg, then $\scro(P)$ is lower-triangular;
if $P$ is lower-Hessenberg with invertible superdiagonal entries,
then $\scro(P)$ is lower-triangular with invertible diagonal entries;
and if $P$ is unit-lower-Hessenberg
(i.e.\ lower-Hessenberg with entries 1 on the superdiagonal),
then $\scro(P)$ is unit-lower-triangular.
In all the applications in this paper, $P$ will be lower-Hessenberg.

The matrix $P$ can also be interpreted as the adjacency matrix
for a weighted directed graph on the vertex set $\N$
(where the edge $ij$ is omitted whenever $p_{ij}  = 0$).
Then $P$ is row-finite (resp.\ column-finite)
if and only if every vertex has finite out-degree (resp.\ finite in-degree).

This iteration process can be given a compact matrix formulation.
Let us define the \textbfit{augmented production matrix}
\be
   \widetilde{P}
   \;\eqdef\;
   \left[
    \begin{array}{c}
         1 \;\; 0 \;\; 0 \;\; 0 \;\; \cdots \; \vphantom{\Sigma} \\
         \hline
         P
    \end{array}
    \right]
   \;.
 \label{def.prodmat_augmented}
\ee
Then the recurrence \reff{def.iteration.bis}
together with the initial condition $a_{0k} = \delta_{0k}$ can be written as
\be
   A
   \;=\;
   \left[
    \begin{array}{c}
         1 \;\; 0 \;\; 0 \;\; 0 \;\; \cdots \; \vphantom{\Sigma} \\
         \hline
         AP
    \end{array}
    \right]
   \;=\;
   \left[
    \begin{array}{c|c}
         1 & \bzero  \\
         \hline
         \bzero & A
    \end{array}
    \right]
   \left[
    \begin{array}{c}
         1 \;\; 0 \;\; 0 \;\; 0 \;\; \cdots \; \vphantom{\Sigma} \\
         \hline
         P
    \end{array}
    \right]
   \;=\;
   \left[
    \begin{array}{c|c}
         1 & \bzero  \\
         \hline
         \bzero & A
    \end{array}
    \right]
   \widetilde{P}
   \;.
 \label{eq.prodmat.u}
\ee
This identity can be iterated to give the factorization
\be
   A
   \;=\;
   \cdots\,
   \left[
    \begin{array}{c|c}
         I_3 & \bzero  \\
         \hline
                &   \\[-4mm]
         \bzero & \widetilde{P}
    \end{array}
    \right]
   \left[
    \begin{array}{c|c}
         I_2 & \bzero  \\
         \hline
                &   \\[-4mm]
         \bzero & \widetilde{P}
    \end{array}
    \right]
   \left[
    \begin{array}{c|c}
         I_1 & \bzero  \\
         \hline
                &   \\[-4mm]
         \bzero & \widetilde{P}
    \end{array}
    \right]
    \widetilde{P}
 \label{eq.prodmat.u.iterated}
\ee
where $I_k$ is the $k \times k$ identity matrix;
and conversely, \reff{eq.prodmat.u.iterated} implies \reff{eq.prodmat.u}.

Now let $\Delta = (\delta_{i+1,j})_{i,j \ge 0}$
be the matrix with 1 on the superdiagonal and 0 elsewhere.
Then for any matrix $M$ with rows indexed by $\N$,
the product $\Delta M$ is simply $M$ with its zeroth row removed
and all other rows shifted upwards.
(Some authors use the notation $\overline{M} \eqdef \Delta M$.)
The recurrence \reff{def.iteration.bis} can then be written as
\be
   \Delta \, \scro(P)  \;=\;  \scro(P) \, P
   \;.
 \label{def.iteration.bis.matrixform}
\ee
It follows that if $A$ is a row-finite matrix
that has a row-finite inverse $A^{-1}$
and has first row $a_{0k} = \delta_{0k}$,
then $P = A^{-1} \Delta A$ is the unique matrix such that $A = \scro(P)$.
This holds, in particular, if $A$ is lower-triangular with
invertible diagonal entries and $a_{00} = 1$;
then $A^{-1}$ is lower-triangular
and $P = A^{-1} \Delta A$ is lower-Hessenberg.
And if $A$ is unit-lower-triangular,
then $P = A^{-1} \Delta A$ is unit-lower-Hessenberg.

We shall repeatedly use the following easy facts:

\begin{lemma}[Production matrix of a product]
   \label{lemma.production.AB}
Let $P = (p_{ij})_{i,j \ge 0}$ be a row-finite matrix
(with entries in a commutative ring $R$),
with output matrix $A = \scro(P)$;
and let $B = (b_{ij})_{i,j \ge 0}$
be a lower-triangular matrix with invertible (in $R$) diagonal entries.
Then
\be
   AB \;=\;  b_{00} \, \scro(B^{-1} P B)
   \;.
\ee
That is, up to a factor $b_{00}$,
the matrix $AB$ has production matrix $B^{-1} P B$.
\end{lemma}

\proof
Since $P$ is row-finite, so is $A = \scro(P)$;
then the matrix products $AB$ and $B^{-1} P B$
arising in the lemma are well-defined.  Now
\be
   a_{nk}
   \;=\;
   \sum_{i_1,\ldots,i_{n-1}}
      p_{0 i_1} \, p_{i_1 i_2} \, p_{i_2 i_3} \,\cdots\,
        p_{i_{n-2} i_{n-1}} \, p_{i_{n-1} k}
   \;,
\ee
while
\be
   \scro(B^{-1} P B)_{nk}
   \;=\;
   \sum_{j,i_1,\ldots,i_{n-1},i_n}
      (B^{-1})_{0j} \,
      p_{j i_1} \, p_{i_1 i_2} \, p_{i_2 i_3} \,\cdots\,
        p_{i_{n-2} i_{n-1}} \, p_{i_{n-1} i_n} \, b_{i_n k}
   \;.
\ee
But $B$ is lower-triangular with invertible diagonal entries,
so $B$ is invertible and $B^{-1}$ is lower-triangular,
with $(B^{-1})_{0j} = b_{00}^{-1} \delta_{j0}$.
It follows that $AB = b_{00} \, \scro(B^{-1} P B)$.
\qed

\begin{lemma}[Production matrix of a down-shifted matrix]
   \label{lemma.down-shifted}
Let $P = (p_{ij})_{i,j \ge 0}$ be a row-finite or column-finite matrix
(with entries in a commutative ring $R$),
with output matrix $A = \scro(P)$;
and let $c$ be an element of $R$.
Now define
\be
   Q
   \;=\;
   \left[
   \begin{array}{c|c@{\hspace*{2mm}}c@{\hspace*{2mm}}c}
        0      &   c & 0 & \cdots \\
        \hline
        0      &     &   &        \\
        0      &     & P &        \\[-1mm]
        \vdots &     &   &        \\
   \end{array}
   \right]
   \;=\;
   c {\bf e}_{01} \,+\, \Delta^{\rm T} P \Delta
\ee
and
\be
   B
   \;=\;
   \left[
   \begin{array}{c|c@{\hspace*{2mm}}c@{\hspace*{2mm}}c}
        1      &   0 & 0 & \cdots \\
        \hline
        0      &     &   &        \\
        0      &     & cA &        \\[-1mm]
        \vdots &     &   &        \\
   \end{array}
   \right]
   \;=\;
   {\bf e}_{00} \,+\, c \Delta^{\rm T} A \Delta
   \;.
\ee
Then $B = \scro(Q)$.
\end{lemma}

\proof
We use \reff{def.iteration.walk} and its analogue for $Q$:
\be
   \scro(Q)_{nk}
   \;=\;
   \sum_{i_1,\ldots,i_{n-1}}
      q_{0 i_1} \, q_{i_1 i_2} \, q_{i_2 i_3} \,\cdots\,
        q_{i_{n-2} i_{n-1}} \, q_{i_{n-1} k}
   \;.
 \label{def.iteration.walk.BQ}
\ee
In \reff{def.iteration.walk.BQ}, the only nonzero contributions come from
$i_1 = 1$, with $q_{01} = c$;
and then we must also have $i_2,i_3,\ldots \ge 1$ and $k \ge 1$,
with $q_{ij} = p_{i-1,j-1}$.
Hence $\scro(Q)_{nk} = c a_{n-1,k-1}$ for $n \ge 1$.
\qed


\subsection{Production matrices and total positivity}
   \label{subsec.totalpos.prodmat}

Let $P = (p_{ij})_{i,j \ge 0}$ be a matrix with entries in a
partially ordered commutative ring $R$.
We will use $P$ as a production matrix;
let $A = \scro(P)$ be the corresponding output matrix.
As before, we assume that $P$ is either row-finite or column-finite.

When $P$ is totally positive, it turns out \cite{Sokal_totalpos}
that the output matrix $\scro(P)$ has {\em two}\/ total-positivity properties:
firstly, it is totally positive;
and secondly, its zeroth column is Hankel-totally positive.
Since \cite{Sokal_totalpos} is not yet publicly available,
we shall present briefly here (with proof) the main results
that will be needed in the sequel.

The fundamental fact that drives the whole theory is the following:

\begin{proposition}[Minors of the output matrix]
   \label{prop.iteration.homo}
Every $k \times k$ minor of the output matrix $A = \scro(P)$
can be written as a sum of products of minors of size $\le k$
of the production matrix $P$.
\end{proposition}

In this proposition the matrix elements $\bfp = \{p_{ij}\}_{i,j \ge 0}$
should be interpreted in the first instance as indeterminates:
for instance, we can fix a row-finite or column-finite set
$S \subseteq \N \times \N$
and define the matrix $P^S = (p^S_{ij})_{i,j \in \N}$ with entries
\be
   p^S_{ij}
   \;=\;
   \begin{cases}
       p_{ij}  & \textrm{if $(i,j) \in S$} \\[1mm]
       0       & \textrm{if $(i,j) \notin S$}
   \end{cases}
\ee
Then the entries (and hence also the minors) of both $P$ and $A$
belong to the polynomial ring $\Z[\bfp]$,
and the assertion of Proposition~\ref{prop.iteration.homo} makes sense.
Of course, we can subsequently specialize the indeterminates $\bfp$
to values in any commutative ring $R$.

\proofof{Proposition~\ref{prop.iteration.homo}}
%
For any infinite matrix $X = (x_{ij})_{i,j \ge 0}$,
let us write $X_N = (x_{ij})_{0 \le i \le N-1 ,\, j \ge 0}$
for the submatrix consisting of the first $N$ rows
(and {\em all}\/ the columns) of $X$.
Every $k \times k$ minor of $A$ is of course
a $k \times k$ minor of $A_N$ for some $N$,
so it suffices to prove that the claim about minors holds for all the $A_N$.
But this is easy: the fundamental identity \reff{eq.prodmat.u} implies
\be
   A_N
   \;=\;
   \left[
    \begin{array}{c|c}
         1 & \bzero  \\
         \hline
         \bzero & A_{N-1}
    \end{array}
    \right]
   \,
   \left[
    \begin{array}{c}
         1 \;\; 0 \;\; 0 \;\; 0 \;\; \cdots \; \vphantom{\Sigma} \\
         \hline
         P
    \end{array}
    \right]
   \;.
 \label{eq.proof.prop.iteration.homo}
\ee
So the result follows by induction on $N$, using the Cauchy--Binet formula.
\qed

If we now specialize the indeterminates $\bfp$
to values in some partially ordered commutative ring $R$,
we can immediately conclude:

\begin{theorem}[Total positivity of the output matrix]
   \label{thm.iteration.homo}
Let $P$ be an infinite matrix that is either row-finite or column-finite,
with entries in a partially ordered commutative ring $R$.
If $P$ is totally positive of order~$r$, then so is $A = \scro(P)$.
\end{theorem}

\medskip

{\bf Remarks.}
1.  In the case $R = \R$, Theorem~\ref{thm.iteration.homo}
is due to Karlin \cite[pp.~132--134]{Karlin_68};
see also \cite[Theorem~1.11]{Pinkus_10}.
Karlin's proof is different from ours.

2.  Our quick inductive proof of Proposition~\ref{prop.iteration.homo}
follows an idea of Zhu \cite[proof of Theorem~2.1]{Zhu_13},
which was in turn inspired in part by Aigner \cite[pp.~45--46]{Aigner_99}.
The same idea recurs in recent work of several authors
\cite[Theorem~2.1]{Zhu_14}
\cite[Theorem~2.1(i)]{Chen_15a}
\cite[Theorem~2.3(i)]{Chen_15b}
\cite[Theorem~2.1]{Liang_16}
\cite[Theorems~2.1 and 2.3]{Chen_19}
\cite{Gao_non-triangular_transforms}.
However, all of these results concerned only special cases:
\cite{Aigner_99,Zhu_13,Chen_15b,Liang_16}
treated the case in which the production matrix $P$ is tridiagonal;
\cite{Zhu_14} treated a (special) case in which $P$ is upper bidiagonal;
\cite{Chen_15a} treated the case in which
$P$ is the production matrix of a Riordan array;
\cite{Chen_19,Gao_non-triangular_transforms}
treated (implicitly) the case in which $P$ is upper-triangular and Toeplitz.
But the argument is in fact completely general, as we have just seen;
there is no need to assume any special form for the matrix $P$.

3. A slightly different version of this proof
was presented in \cite{latpath_SRTR,latpath_lah}.
The simplified reformulation
given here,
using the augmented production matrix,
is due to Mu and Wang \cite{Mu_20}.
\myendremark

\medskip

\begin{example}[Toeplitz matrix of powers]
   \label{exam.toeplitz.power.TP}
\rm
Let $P = x {\bf e}_{00} + y \Delta$,
where $x$ and $y$ are indeterminates
(here ${\bf e}_{ij}$ denotes the matrix with an entry~1 in position~$ij$
 and 0 elsewhere).
By Lemma~\ref{lemma.bidiagonal}, $P$ is TP
in the ring $\Z[x,y]$ equipped with the coefficientwise order.
An easy computation shows that
$\scro(x {\bf e}_{00} + y\Delta)_{nk} = x^{n-k} y^k \, {\rm I}[k \le n]$.
(Here ${\rm I}[\hbox{\sl proposition}] = 1$ if {\sl proposition} is true,
and 0 if it is false.)
When $y=1$, this is the Toeplitz matrix of powers \reff{def.Tx}.
So Theorem~\ref{thm.iteration.homo} implies that $T_x$ is TP
in the ring $\Z[x]$ equipped with the coefficientwise order.
This gives a second proof of the total positivity stated in
Lemma~\ref{lemma.toeplitz.power}.
\myendremark
\end{example}
\vspace*{-8mm}

\begin{example}[Binomial matrix]
   \label{exam.binomial.matrix.TP}
\rm
Let $P$ be the upper-bidiagonal Toeplitz matrix
$xI + y\Delta$, where $x$ and $y$ are indeterminates.
By Lemma~\ref{lemma.bidiagonal}, $P$ is TP
in the ring $\Z[x,y]$ equipped with the coefficientwise order.
An easy computation shows that $\scro(xI + y\Delta) = B_{x,y}$,
the weighted binomial matrix
with entries $(B_{x,y})_{nk} = x^{n-k} y^k \binom{n}{k}$.
So Theorem~\ref{thm.iteration.homo} implies that $B_{x,y}$ is TP
in the ring $\Z[x,y]$ equipped with the coefficientwise order.
This gives an {\em ab initio}\/ proof of Lemma~\ref{lemma.binomialmatrix.TP}.
\myendremark
\end{example}

\bigskip

Now define 
$\scroo_0(P)$ to be the zeroth-column sequence of $\scro(P)$, i.e.
\be
   \scroo_0(P)_n  \;\eqdef\;  \scro(P)_{n0}  \;\eqdef\;  (P^n)_{00}
   \;.
 \label{def.scroo0}
\ee
Then the Hankel matrix of $\scroo_0(P)$ has matrix elements
\begin{eqnarray}
   & &
   \!\!\!\!\!\!\!
   H_\infty(\scroo_0(P))_{nn'}
   \;=\;
   \scroo_0(P)_{n+n'}
   \;=\;
   (P^{n+n'})_{00}
   \;=\;
   \sum_{k=0}^\infty (P^n)_{0k} \, (P^{n'})_{k0}
   \;=\;
          \nonumber \\
   & &
   \sum_{k=0}^\infty (P^n)_{0k} \, ((P^{\rm T})^{n'})_{0k}
   \;=\;
   \sum_{k=0}^\infty \scro(P)_{nk} \, \scro(P^{\rm T})_{n'k}
   \;=\;
   \big[ \scro(P) \, {\scro(P^{\rm T})}^{\rm T} \big]_{nn'}
   \;.
   \qquad
\end{eqnarray}
(Note that the sum over $k$ has only finitely many nonzero terms:
 if $P$ is row-finite, then there are finitely many nonzero $(P^n)_{0k}$,
 while if $P$ is column-finite,
 there are finitely many nonzero $(P^{n'})_{k0}$.)
We have therefore proven:

\begin{lemma}[Identity for Hankel matrix of the zeroth column]
   \label{lemma.hankel.karlin}
Let $P$ be a row-finite or column-finite matrix
with entries in a commutative ring $R$.
Then
\be
   H_\infty(\scroo_0(P))
   \;=\;
   \scro(P) \, {\scro(P^{\rm T})}^{\rm T}
   \;.
\ee
\end{lemma}

{\bf Remark.}
If $P$ is row-finite, then $\scro(P)$ is row-finite;
$\scro(P^{\rm T})$ need not be row- or column-finite,
but the product $\scro(P) \, {\scro(P^{\rm T})}^{\rm T}$
is anyway well-defined.
Similarly, if $P$ is column-finite, then ${\scro(P^{\rm T})}^{\rm T}$
is column-finite;
$\scro(P)$ need not be row- or column-finite,
but the product $\scro(P) \, {\scro(P^{\rm T})}^{\rm T}$
is again well-defined.
\myendremark

\medskip

Combining Proposition~\ref{prop.iteration.homo}
with Lemma~\ref{lemma.hankel.karlin} and the Cauchy--Binet formula,
we obtain:

\begin{corollary}[Hankel minors of the zeroth column]
   \label{cor.iteration2}
Every $k \times k$ minor of the infinite Hankel matrix
$H_\infty(\scroo_0(P)) = ((P^{n+n'})_{00})_{n,n' \ge 0}$
can be written as a sum of products
of the minors of size $\le k$ of the production matrix $P$.
\end{corollary}

And specializing the indeterminates $\bfp$
to nonnegative elements in a partially ordered commutative ring,
in such a way that $P$ is row-finite or column-finite,
we deduce:

\begin{theorem}[Hankel-total positivity of the zeroth column]
   \label{thm.iteration2bis}
Let $P = (p_{ij})_{i,j \ge 0}$ be an infinite row-finite or column-finite
matrix with entries in a partially ordered commutative ring $R$,
and define the infinite Hankel matrix
$H_\infty(\scroo_0(P)) = ((P^{n+n'})_{00})_{n,n' \ge 0}$.
If $P$ is totally positive of order~$r$, then so is $H_\infty(\scroo_0(P))$.
\end{theorem}

One might hope that Theorem~\ref{thm.iteration2bis}
could be strengthened to show not only Hankel-TP of the zeroth column
of the output matrix $A = \scro(P)$,
but in fact Hankel-TP of the row-generating polynomials $A_n(x)$
for all $x \ge 0$ (at least when $R = \R$) ---
or even more strongly,
coefficientwise Hankel-TP of the row-generating polynomials.
Alas, this hope is vain, for these properties do not hold {\em in general}\/:

\begin{example}[Failure of Hankel-TP of the row-generating polynomials]
   \label{exam.hankel.karlin.rowgen}
\rm
Let $P = {\bf e}_{00} + \Delta$
be the upper-bidiagonal matrix with 1 on the superdiagonal
and $1,0,0,0,\ldots$ on the diagonal;
by Lemma~\ref{lemma.bidiagonal} it is TP.
Then $A = \scro(P)$ is the lower-triangular matrix will all entries 1
(see Example~\ref{exam.toeplitz.power.TP}),
so that $A_n(x) = \sum_{k=0}^n x^k$.
Since $A_0(x) \, A_2(x) - A_1(x)^2 = -x$,
the sequence $(A_n(x))_{n \ge 0}$ is not even log-convex
(i.e.\ Hankel-TP${}_2$) for any real number $x > 0$.
%
\myendremark
\end{example}

Nevertheless, in one important special case ---
which includes all the matrices arising in the present paper ---
the total positivity of the production matrix {\em does}\/
imply the coefficientwise Hankel-TP of the row-generating polynomials
of the output matrix:
see Theorem~\ref{thm.hankelTP.expriordan} below.

\subsection[An identity for $B_x^{-1} P B_x$]{An identity for $\bm{B_x^{-1} P B_x}$}

An important role will be played later in this paper
by a simple but remarkable identity \cite[Lemma~3.6]{latpath_lah}
for $B_x^{-1} P B_x$,
where $B_x$ is the $x$-binomial matrix
and $P$ is a particular diagonal similarity transform (by factorials)
of a lower-Hessenberg Toeplitz matrix:

\begin{lemma}[Identity for $B_x^{-1} P B_x$]
   \label{lemma.BxinvPBx}
Let $\bphi = (\phi_i)_{i \ge 0}$ and $x$ be indeterminates,
and work in the ring $\Z[\bphi,x]$.
Define the lower-Hessenberg matrix $P = (p_{ij})_{i,j \ge 0}$ by
\be
   p_{ij}
   \;=\; 
   \begin{cases}
      0   &  \textrm{if $j=0$ or $j > i+1$}  \\[1mm]
      {\displaystyle {i! \over (j-1)!}} \, \phi_{i-j+1}
          &  \textrm{if $1 \le j \le i+1$}
   \end{cases}
 \label{eq.lemma.BxinvPBx.1}
\ee
and the unit-lower-triangular $x$-binomial matrix $B_x$ by
\be
   (B_x)_{nk}  \;=\;  \binom{n}{k} \, x^{n-k}
   \;.
\ee
Let $\Delta = (\delta_{i+1,j})_{i,j \ge 0}$
be the matrix with 1 on the superdiagonal and 0 elsewhere.
Then
\be
   B_x^{-1} P B_x   \;=\;   P (I + x \Delta^{\rm T})
   \;.
 \label{eq.lemma.BxinvPBx}
\ee
\end{lemma}

In \cite{latpath_lah} we proved \reff{eq.lemma.BxinvPBx}
by a computation using a binomial sum.
Here is a simpler proof:

\proofof{Lemma~\ref{lemma.BxinvPBx}}
We have $P = D T_\infty(\bphi) D^{-1} \Delta$
and $B_x = D \, T_\infty\big( (x^n/n!)_{n \ge 0} \big) \, D^{-1}$,
where $D = \diag\big( (n!)_{n \ge 0} \big)$.
Now two Toeplitz matrices always commute:
$T_\infty(\ba) \, T_\infty(\bb) =
 T_\infty(\ba * \bb) =
 T_\infty(\bb) \, T_\infty(\ba)$.
It follows that $D T_\infty(\bphi) D^{-1}$ and $B_x$ commute.
On the other hand, the classic recurrence for binomial coefficients implies
\be
   \Delta B_x  \;=\; B_x \, (xI + \Delta)
\ee
(cf.\ Example~\ref{exam.binomial.matrix.TP}).  Therefore
\begin{subeqnarray}
   B_x^{-1} P B_x
   & = &
   B_x^{-1} \, D T_\infty(\bphi) D^{-1} \, \Delta B_x
      \\[2mm]
   & = &
   B_x^{-1} \, D T_\infty(\bphi) D^{-1} \, B_x \, (xI + \Delta)
      \\[2mm]
   & = &
   D T_\infty(\bphi) D^{-1} \, (xI + \Delta)
      \\[2mm]
   & = &
   D T_\infty(\bphi) D^{-1} \, \Delta \, (I + x\Delta^{\rm T})
\end{subeqnarray}
since $\Delta \Delta^{\rm T} = I$.
\qed

\subsection{A lemma on diagonal scaling}   \label{subsec.diagonal_scaling}

Given a lower-triangular matrix $A = (a_{nk})_{n,k \ge 0}$ 
with entries in a commutative ring $R$,
let us define the matrix $A^\sharp = (a^\sharp_{nk})_{n,k \ge 0}$ by
\be
   a^\sharp_{nk}  \;=\;  {n! \over k!} \: a_{nk}
   \;;
\ee
this is well-defined since $a_{nk} \neq 0$ only when $n \ge k$,
in which case $n!/k!$ is an integer.

If $R$ contains the rationals, we can of course write
$A^\sharp = D A D^{-1}$ where $D = \diag\big( (n!)_{n \ge 0} \big)$.
And if $R$ is a partially ordered commutative ring
that contains the rationals and $A$ is TP${}_r$,
then we deduce immediately from $A^\sharp = D A D^{-1}$
that also $A^\sharp$ is TP${}_r$.
The following simple lemma \cite[Lemma~3.7]{latpath_lah}
shows that this conclusion holds even when $R$ does not contain the rationals:

\begin{lemma}
   \label{lemma.diagmult.TP}
Let $A = (a_{ij})_{i,j \ge 0}$ be a lower-triangular matrix
with entries in a partially ordered commutative ring $R$,
and let $\bd = (d_i)_{i \ge 1}$.
Define the lower-triangular matrix
$A^{\sharp\bd} = (a^{\sharp\bd}_{ij})_{i,j \ge 0}$ by
\be
   a^{\sharp\bd}_{ij}  \;=\;  d_{j+1} d_{j+2} \cdots d_i \, a_{ij}
   \;.
\ee
Then:
\begin{itemize}
   \item[(a)] If $A$ is TP${}_r$ and $\bd$ are indeterminates,
      then $A^{\sharp\bd}$ is TP${}_r$ in the ring $R[\bd]$ equipped with
      the coefficientwise order.
   \item[(b)] If $A$ is TP${}_r$ and $\bd$ are nonnegative elements of $R$,
      then $A^{\sharp\bd}$ is TP${}_r$ in the ring $R$.
\end{itemize}
\end{lemma}

\proof
(a) Let $\bd = (d_i)_{i \ge 1}$ be commuting indeterminates,
and let us work in the ring $R[\bd,\bd^{-1}]$
equipped with the coefficientwise order.
Let $D = \diag(1,\, d_1,\, d_1 d_2,\, \ldots)$.
Then $D$ is invertible, and both $D$ and
$D^{-1} = \diag(1,\, d_1^{-1},\, d_1^{-1} d_2^{-1},\, \ldots)$
have nonnegative elements.
It follows that $A^{\sharp\bd} = D A D^{-1}$ is TP${}_r$
in the ring $R[\bd,\bd^{-1}]$ equipped with the coefficientwise order.
But the matrix elements $a^{\sharp\bd}_{ij}$
actually belong to the subring $R[\bd] \subseteq R[\bd,\bd^{-1}]$.
So $A^{\sharp\bd}$ is TP${}_r$ in the ring $R[\bd]$
equipped with the coefficientwise order.

(b) follows from (a) by specializing indeterminates.
\qed

\noindent
The special case $A^{\sharp\bd} = A^\sharp$ corresponds to taking $d_i = i$.

\subsection{Exponential Riordan arrays}

Let $R$ be a commutative ring containing the rationals,
and let $F(t) = \sum_{n=0}^\infty f_n t^n/n!$
and $G(t) = \sum_{n=1}^\infty g_n t^n/n!$ be formal power series
with coefficients in $R$; we set $g_0 = 0$.
Then the \textbfit{exponential Riordan array}
\cite{Deutsch_04,Deutsch_09,Barry_16}
associated to the pair $(F,G)$
is the infinite lower-triangular matrix
$\scrr[F,G] = (\scrr[F,G]_{nk})_{n,k \ge 0}$ defined by
\be
   \scrr[F,G]_{nk}
   \;=\;
   {n! \over k!} \:
   [t^n] \, F(t) G(t)^k
   \;.
 \label{def.RFG}
\ee
That is, the $k$th column of $\scrr[F,G]$
has exponential generating function $F(t) G(t)^k/k!$.
The bivariate egf is
\be
   \sum_{n=0}^\infty \sum_{k=0}^n \scrr[F,G]_{nk} \, {t^n \over n!} \, x^k
   \;=\;
   F(t) \, e^{x G(t)}
   \;.
\ee
Please note that the diagonal elements of $\scrr[F,G]$
are $\scrr[F,G]_{nn} = f_0 g_1^n$,
so the matrix $\scrr[F,G]$ is invertible
in the ring $R^{\N \times \N}_{\rm lt}$ of lower-triangular matrices
if and only if $f_0$ and $g_1$ are invertible in $R$.

We shall use an easy but important result that is sometimes called
the \emph{fundamental theorem of exponential Riordan arrays} (FTERA):

\begin{lemma}[Fundamental theorem of exponential Riordan arrays]
   \label{lemma.FETRA}
Let $\bb = (b_n)_{n \ge 0}$ be a sequence with
exponential generating function $B(t) = \sum_{n=0}^\infty b_n t^n/n!$.
Considering $\bb$ as a column vector and letting $\scrr[F,G]$
act on it by matrix multiplication, we obtain a sequence $\scrr[F,G] \bb$
whose exponential generating function is $F(t) \, B(G(t))$.
\end{lemma}

\proof
We compute
\begin{subeqnarray}
   \sum_{k=0}^n \scrr[F,G]_{nk} \, b_k
   & = &
   \sum_{k=0}^\infty {n! \over k!} \, [t^n] \, F(t) G(t)^k \, b_k
            \\[2mm]
   & = &
   n! \: [t^n] \: F(t) \sum_{k=0}^\infty b_k \, {G(t)^k \over k!}
            \\[2mm]
   & = &
   n! \: [t^n] \: F(t) \, B(G(t))
   \;.
\end{subeqnarray}
\qed

We can now determine the production matrix of an exponential
Riordan array $\scrr[F,G]$:

\begin{theorem}[Production matrices of exponential Riordan arrays]
   \label{thm.riordan.exponential.production}
Let $L$ be a lower-triangular matrix
(with entries in a commutative ring $R$ containing the rationals)
with invertible diagonal entries and $L_{00} = 1$,
and let $P = L^{-1} \Delta L$ be its production matrix.
Then $L$ is an exponential Riordan array
if and only~if $P = (p_{nk})_{n,k \ge 0}$ has the form
\be
   p_{nk}
   \;=\;
   {n! \over k!} \: (z_{n-k} \,+\, k \, a_{n-k+1})
 \label{eq.thm.riordan.exponential.production}
\ee
for some sequences $\ba = (a_n)_{n \ge 0}$ and $\bz = (z_n)_{n \ge 0}$
in $R$.

More precisely, $L = \scrr[F,G]$ if and only~if $P$
is of the form \reff{eq.thm.riordan.exponential.production}
where the ordinary generating functions
$A(s) = \sum_{n=0}^\infty a_n s^n$ and $Z(s) = \sum_{n=0}^\infty z_n s^n$
are connected to $F(t)$ and $G(t)$ by
\be
   G'(t) \;=\; A(G(t))  \;,\qquad
   {F'(t) \over F(t)} \;=\; Z(G(t))
 \label{eq.prop.riordan.exponential.production.1}
\ee
or equivalently
\be
   A(s)  \;=\;  G'(\bar{G}(s))  \;,\qquad
   Z(s)  \;=\;  {F'(\bar{G}(s)) \over F(\bar{G}(s))}
 \label{eq.prop.riordan.exponential.production.2}
\ee
where $\bar{G}(s)$ is the compositional inverse of $G(t)$.
\end{theorem}

\par\bigskip\noindent{\sc Proof}
(mostly contained in \cite[pp.~217--218]{Barry_16}).
Suppose that $L = \scrr[F,G]$.
The hypotheses on $L$ imply that $f_0 = 1$
and that $g_1$ is invertible in $R$;
so $G(t)$ has a compositional inverse.
Now let $P = (p_{nk})_{n,k \ge 0}$ be a matrix;
its column exponential generating functions are, by definition,
$P_k(t) = \sum_{n=0}^\infty p_{nk} \, t^n/n!$.
Applying the FTERA to each column of $P$,
we see that $\scrr[F,G] P$ is a matrix
whose column exponential generating functions
are $\big( F(t) \, P_k(G(t)) \big)_{k \ge 0}$.
On~the other hand, $\Delta \, \scrr[F,G]$
is the matrix $\scrr[F,G]$ with its zeroth row removed
and all other rows shifted upwards,
so it has column exponential generating functions
\be
   {d \over dt} \, \big( F(t) \, G(t)^k/k! \big)
   \;=\;
   {1 \over k!} \: \Big[ F'(t) \, G(t)^k
                         \:+\: k \, F(t) \, G(t)^{k-1} \, G'(t) \Big]
   \;.
\ee
Comparing these two results, we see that
$\Delta \, \scrr[F,G] = \scrr[F,G] \, P$
if and only~if
\be
   P_k(G(t))
   \;=\;
   {1 \over k!} \:
   {F'(t) \, G(t)^k \:+\: k \, F(t) \, G(t)^{k-1} \, G'(t)
    \over
    F(t)}
   \;,
\ee
or in other words
\be
   P_k(t)
   \;=\;
   {1 \over k!}  \:
      \biggl[ {F'(\bar{G}(t)) \over F(\bar{G}(t))} \, t^k
              \:+\: k \, t^{k-1} \, G'(\bar{G}(t))
      \biggr]
   \;.
\ee
Therefore
\begin{subeqnarray}
   p_{nk}
   & = &
   {n! \over k!} \: [t^n] \,
      \biggl[ {F'(\bar{G}(t)) \over F(\bar{G}(t))} \, t^k
              \:+\: k \, t^{k-1} \, G'(\bar{G}(t))
      \biggr]
    \\[2mm]
   & = &
   {n! \over k!} \:
      \biggl[ [t^{n-k}] \: {F'(\bar{G}(t)) \over F(\bar{G}(t))}
              \:+\: k \, [t^{n-k+1}] \: G'(\bar{G}(t))
      \biggr]
    \\[2mm]
   & = &
   {n! \over k!} \: (z_{n-k} \,+\, k \, a_{n-k+1})
\end{subeqnarray}
where $\ba = (a_n)_{n \ge 0}$ and $\bz = (z_n)_{n \ge 0}$
are given by \reff{eq.prop.riordan.exponential.production.2}.

Conversely, suppose that $P = (p_{nk})_{n,k \ge 0}$ has the form
\reff{eq.thm.riordan.exponential.production}.
Define $F(t)$ and $G(t)$
as the unique solutions (in the formal-power-series ring $R[[t]]$)
of the differential equations \reff{eq.prop.riordan.exponential.production.1}
with initial conditions $F(0) = 1$ and $G(0) = 0$.
Then running the foregoing computation backwards
shows that $\Delta \, \scrr[F,G] = \scrr[F,G] \, P$.
\qed

The exponential Riordan arrays arising in the present paper
will all have $F(t) = 1$:
these are said to belong to the {\em associated subgroup}\/
(or {\em Lagrange subgroup}\/).
Such matrices (sometimes with the zeroth row and column removed)
are also known as {\em Jabotinsky matrices}\/ \cite{Jabotinsky_47}
or {\em convolution matrices}\/ \cite{Knuth_92}.
Their entries are also identical to the
{\em partial Bell polynomials}\/ \cite[pp.~133--137]{Comtet_74}
${\bf B}_{n,k}(g_1,g_2,\ldots)$
where $G(t) = \sum_{n=1}^\infty g_n t^n/n!$.

Let us also observe that the matrices $P$
occurring in Lemma~\ref{lemma.BxinvPBx}
are precisely the production matrices
\reff{eq.thm.riordan.exponential.production}
with $\bz = 0$ (and $\ba = \bphi$):
that is, they are the production matrices of exponential Riordan arrays
$\scrr[F,G]$ with $F(t) = 1$.
This observation allows us to improve Theorem~\ref{thm.iteration2bis}
--- from Hankel-total positivity of the zeroth column to
coefficientwise Hankel-total positivity of the row-generating polynomials ---
for the special case of exponential Riordan arrays
$\scrr[F,G]$ with $F(t) = 1$:

\begin{theorem}[Hankel-TP for row-generating polynomials of exponential
   Riordan array]
   \label{thm.hankelTP.expriordan}
Let $R$ be a partially ordered commutative ring containing the rationals;
let $A = (a_{nk})_{n,k \ge 0} = \scrr[1,G]$ be an exponential Riordan array
of the associated subgroup,
with entries in $R$ and with invertible diagonal elements;
let $A_n(x) = \sum_{k=0}^n a_{nk} \, x^k$
be its row-generating polynomials;
and let $P = A^{-1} \Delta A$ be its production matrix.

If $P$ is totally positive of order~$r$ in the ring $R$,
then the sequence $(A_n(x))_{n \ge 0}$ of row-generating polynomials
is Hankel-totally positive of order~$r$
in the ring $R[x]$ equipped with the coefficientwise order.
\end{theorem}

\proof
The row-generating polynomials $A_n(x)$ form the zeroth column
of the binomial row-generating matrix $A B_x$.
By Lemma~\ref{lemma.production.AB},
the production matrix of $A B_x$ is $B_x^{-1} P B_x$.
By Theorem~\ref{thm.riordan.exponential.production},
the production matrix $P = (p_{nk})_{n,k \ge 0}$ has the form
\be
   p_{nk}
   \;=\;
   {n! \over (k-1)!} \: a_{n-k+1}
\ee
for some sequence $\ba = (a_n)_{n \ge 0}$ in $R$.
By Lemma~\ref{lemma.BxinvPBx}, we have
\be
   B_x^{-1} P B_x   \;=\;   P (I + x \Delta^{\rm T})
   \;.
\ee
By Lemma~\ref{lemma.bidiagonal}, the matrix $I + x \Delta^{\rm T}$
is totally positive
in the ring $\Z[x]$ equipped with the coefficientwise order;
and by hypothesis,
the matrix $P$ is totally positive of order~$r$ in the ring~$R$.
It follows that $B_x^{-1} P B_x$ is totally positive of order~$r$
in the ring $R[x]$ equipped with the coefficientwise order.
Theorem~\ref{thm.iteration2bis} then implies that
the sequence $(A_n(x))_{n \ge 0}$ of row-generating polynomials
is Hankel-totally positive of order~$r$
in the ring $R[x]$ equipped with the coefficientwise order.
\qed

\subsection{Lagrange inversion}

We will use Lagrange inversion in the following form \cite{Gessel_16}:
If $\Phi(u)$ is a formal power series
with coefficients in a commutative ring $R$ containing the rationals,
then there exists a unique formal power series $f(t)$
with zero constant term satisfying
\be
   f(t)  \;=\;  t \, \Phi(f(t))
   \;,
\ee
and it is given by
\be
   [t^n] \, f(t)  \;=\;  {1 \over n} \, [u^{n-1}] \, \Phi(u)^n
     \quad\hbox{for $n \ge 1$}
   \;;
\ee
and more generally, if $H(u)$ is any formal power series, then
\be
   [t^n] \, H(f(t))  \;=\;  {1 \over n} \, [u^{n-1}] \, H'(u) \, \Phi(u)^n
     \quad\hbox{for $n \ge 1$}
   \;.
 \label{eq.lagrange.H}
\ee
In particular, taking $H(u) = u^k$ with integer $k \ge 0$, we have
\be
   [t^n] \, f(t)^k  \;=\;  {k \over n} \, [u^{n-k}] \, \Phi(u)^n
     \quad\hbox{for $n \ge 1$}
   \;.
 \label{eq.lagrange.k}
\ee

\section[The matrices $\bm{(f_{n,k})_{n,k \ge 0}}$,
            $\bm{(f_{n,k}(y,z))_{n,k \ge 0}}$ and
            $\bm{(f_{n,k}(y,\bphi))_{n,k \ge 0}}$
            as exponential Riordan arrays]{The matrices $\bm{(f_{n,k})_{n,k \ge 0}}$,
            $\bm{(f_{n,k}(y,z))_{n,k \ge 0}}$ and \hfill\break
            $\bm{(f_{n,k}(y,\bphi))_{n,k \ge 0}}$
            as exponential Riordan arrays}   \label{sec.EGF}

In this section we show that the matrices
$(f_{n,k})_{n,k \ge 0}$, $(f_{n,k}(y,z))_{n,k \ge 0}$ and
\linebreak
$(f_{n,k}(y,\bphi))_{n,k \ge 0}$
are exponential Riordan arrays $\scrr[F,G]$ with $F=1$,
and we compute their generating functions $G$.
Much of the contents of the first two subsections is known
\cite{Barry_10_OEIS,Dumont_96},
but we think it useful to bring it all together in one place;
it will motivate our generalization in Section~\ref{subsec.EGF.3}
and will play a key role in the remainder of the paper.

\subsection[The matrix $(f_{n,k})_{n,k \ge 0}$]{The matrix $\bm{(f_{n,k})_{n,k \ge 0}}$}  \label{subsec.EGF.1}

We recall that $f_{n,k}$ is defined combinatorially
as the number of $k$-component forests of rooted trees
on a total of $n$ labeled vertices.
Such a forest can be constructed as follows:
partition the vertex set $V$ into subsets $V_1,\ldots,V_k$
of cardinalities $n_i = |V_i| \ge 1$;
construct a rooted tree on each subset $V_i$;
and finally divide by $k!$ because the trees are distinguishable
(since they are labeled) and any permutation of them gives rise
to the same forest.  It follows that
\be
   f_{n,k}
   \;=\;
   {1 \over k!}
   \!\!\!
   \sum\limits_{\begin{scarray}
                   n_1,\ldots,n_k \ge 1 \\
                   n_1 +\ldots+ n_k = n
                \end{scarray}}
   \!\!\!\!\!
   \binom{n}{n_1,\ldots,n_k} \, f_{n_1,1} \,\cdots\, f_{n_k,1}
   \;.
 \label{eq.EGF.star0}
\ee
In terms of the column exponential generating functions
\be
   \scrf_k(t)  \;\eqdef\;  \sum_{n=0}^\infty f_{n,k} \, {t^n \over n!}
   \;,
\ee
we have
\be
   \scrf_k(t)  \;=\;  {\scrf_1(t)^k \over k!}
   \;.
 \label{eq.EGF.star1}
\ee
It follows from \reff{eq.EGF.star1} and \reff{def.RFG}
that the matrix $(f_{n,k})_{n,k \ge 0}$
is the exponential Riordan array $\scrr[F,G]$
with $F(t) = 1$ and $G(t) = \scrf_1(t)$.

On the other hand, a rooted tree on $n$ labeled vertices
can be obtained by choosing a root and then forming a forest
of rooted trees on the remaining $n-1$ labeled vertices: thus
\be
   f_{n,1}  \;=\; n \, \sum_{k=0}^\infty f_{n-1,k}
   \;.
\ee
Multiplying by $t^n/n!$ and summing over $n \ge 1$, we get
\begin{subeqnarray}
   \scrf_1(t)
   & = &
   t \, \sum_{k=0}^\infty \scrf_k(t)
       \\[2mm]
   & = &
   t \, e^{\scrf_1(t)} \quad\hbox{by \reff{eq.EGF.star1}}
   \;.
     \slabel{eq.EGF.star2.b}
\end{subeqnarray}
This is the well-known functional equation for the
exponential generating function of rooted trees.

We can now (as is also well known\footnote{
   See e.g.\ \cite[Example~5.4.4]{Stanley_99}.
})
apply Lagrange inversion to the functional equation \reff{eq.EGF.star2.b}
to compute $f_{n,k}$.  Using \reff{eq.lagrange.k}, we have
\be
   [t^n] \, \scrf_1(t)^k
   \;=\;
   {k \over n} \, [u^{n-k}] \, (e^u)^n
   \;=\;
   {k \over n} \: {n^{n-k} \over (n-k)!}
\ee
and hence, using \reff{eq.EGF.star1},
\be
   f_{n,k}  \;=\; n! \, [t^n] \, \scrf_k(t)
            \;=\; {n! \over k!} \, [t^n] \, \scrf_1(t)^k
            \;=\;  \binom{n}{k} \, k \, n^{n-k-1}
   \;,
\ee
in agreement with \reff{def.fnk}.
This is, of course, one of the many classic proofs of \reff{def.fnk}.
In~particular, for $k=1$ we have
$\scrf_1(t) = \sum_{n=1}^\infty n^{n-1} \, t^n/n!$,
which is the celebrated {\em tree function}\/ $T(t)$ \cite{Corless_96}.

All this is, of course, extremely well known
(except possibly for the interpretation as an exponential Riordan array,
 which is known \cite{Barry_10_OEIS} but perhaps not as well known
 as it should be).
It is, however, a useful warm-up for the generalization
in which we introduce the variables $y$ and $z$,
to which we now turn.

\subsection[The matrix $(f_{n,k}(y,z))_{n,k \ge 0}$]{The matrix $\bm{(f_{n,k}(y,z))_{n,k \ge 0}}$}  \label{subsec.EGF.2}

Recall that $f_{n,k}(y,z)$ is defined combinatorially
as the generating polynomial for $k$-component forests of rooted trees
on the vertex set $[n]$,
in which each improper edge gets a weight $y$
and each proper edge gets a weight $z$.
The reasoning leading to the identity \reff{eq.EGF.star0}
generalizes without any change whatsoever to $f_{n,k}(y,z)$:
the point is that each set $V_i$ is order-isomorphic to $[n_i]$
(by labeling the vertices in increasing order),
so that the meaning of ``proper edge'' is unaltered.
Therefore
\be
   f_{n,k}(y,z)
   \;=\;
   {1 \over k!}
   \!\!\!
   \sum\limits_{\begin{scarray}
                   n_1,\ldots,n_k \ge 1 \\
                   n_1 +\ldots+ n_k = n
                \end{scarray}}
   \!\!\!\!\!
   \binom{n}{n_1,\ldots,n_k} \, f_{n_1,1}(y,z) \,\cdots\, f_{n_k,1}(y,z)
   \;.
 \label{eq.EGF.star0bis}
\ee
In terms of the column exponential generating functions
\be
   \scrf_k(t;y,z)  \;\eqdef\;  \sum_{n=0}^\infty f_{n,k}(y,z) \, {t^n \over n!}
   \;,
\ee
we have
\be
   \scrf_k(t;y,z)  \;=\;  {\scrf_1(t;y,z)^k \over k!}
   \;.
 \label{eq.EGF.star1bis}
\ee
Therefore, the matrix $(f_{n,k}(y,z))_{n,k \ge 0}$
is the exponential Riordan array $\scrr[F,G]$
with $F(t) = 1$ and $G(t) = \scrf_1(t;y,z)$.
This fact will play a key role in the remainder of the paper.

Of course, it still remains to calculate the exponential generating function
$\scrf_1(t;y,z)$.  This calculation is not at all trivial,
but it was done a quarter-century ago by
Dumont and Ramamonjisoa \cite{Dumont_96};
we need only translate their results to our notation.

Let $\scrt^\bullet_n = \biguplus\limits_{i=1}^n \scrt^{[i]}_n$
denote the set of rooted trees on the vertex set $[n]$,
where $\scrt^{[i]}_n$ is the subset for which the root vertex is $i$.
Let $\scrt^{\<i\>}_n$ denote the subset of $\scrt^\bullet_n$
in which the vertex $i$ is a leaf (i.e.\ has no children).
Given a tree $T \in \scrt^\bullet_n$,
we write $\imp(T)$ for the number of improper edges of $T$.
Now define the generating polynomials
\begin{eqnarray}
   R_n(y,z)  \;=\;  f_{n,1}(y,z)
   & = &
   \sum_{T \in \scrt^\bullet_n} y^{\imp(T)} z^{n-1-\imp(T)}
   \qquad\qquad
        \label{def.Rn}  \\[2mm]
   S_n(y,z)
   & = &
   \sum_{T \in \scrt^{[1]}_{n+1}} y^{\imp(T)} z^{n-\imp(T)}
        \label{def.Sn}  \\[2mm]
   A_n(y,z)
   & = &
   \sum_{T \in \scrt^{\<1\>}_{n+1}} y^{\imp(T)} z^{n-\imp(T)}
        \label{def.An}
\end{eqnarray}
in which each improper (resp.\ proper) edge gets a weight $y$ (resp.~$z$),
and the corresponding exponential generating functions
\begin{eqnarray}
   \scrr(t;y,z)  \;=\;  \scrf_1(t;y,z)
   & = &  \sum\limits_{n=1}^\infty R_n(y,z) \: {t^n \over n!}
   \qquad\qquad\qquad
        \label{def.Rn.EGF}  \\[2mm]
   \scrs(t;y,z)  & = &  \sum\limits_{n=0}^\infty S_n(y,z) \: {t^n \over n!}
        \label{def.Sn.EGF}  \\[2mm]
   \scra(t;y,z)  & = &  \sum\limits_{n=0}^\infty A_n(y,z) \: {t^n \over n!}
        \label{def.An.EGF}
\end{eqnarray}
We then have the following key result \cite[Proposition~7]{Dumont_96}:

\begin{proposition}[Dumont--Ramamonjisoa]
   \label{prop.dumont}
The series $\scrr$, $\scrs$ and $\scra$ satisfy the following identities:
\begin{itemize}
   \item[(a)]  $\scrs(t;y,z) \:=\: \exp\big[ z \, \scrr(t;y,z) \bigr]$
   \item[(b)]  $\scra(t;y,z) \:=\: \displaystyle{1 \over 1 - y \scrr(t;y,z)}$
   \item[(c)]  $\displaystyle {d \over dt} \scrr(t;y,z)
                             \:=\: \scra(t;y,z) \, \scrs(t;y,z)$
\end{itemize}
and hence
\begin{itemize}
   \item[(d)] $\displaystyle {d \over dt} \scrr(t;y,z)
                             \:=\: 
               {\exp\big[ z \, \scrr(t;y,z) \bigr]
                \over
                1 - y \scrr(t;y,z)}$
\end{itemize}
\end{proposition}

Solving the differential equation of Proposition~\ref{prop.dumont}(d)
with the initial condition
$\scrr(0;y,z) = 0$, we obtain:

\begin{corollary}
   \label{cor.dumont}
The series $\scrr(t;y,z)$ satisfies the functional equation
\be
   y - z + yz \scrr  \;=\;  (y - z + z^2 t) \, e^{z\scrr}
 \label{eq.cor.dumont.1}
\ee
and hence has the solution
\be
   \scrr(t;y,z)
   \;=\;
   {1 \over z}
   \biggl[ T\Big( \Big(1 - {z \over y} + {z^2 \over y} t \Big) \:
                  e^{- \, \displaystyle \big(1 - {z \over y} \big)}
            \Big)
            \:-\: \Big(1 - {z \over y} \Big)
   \biggr]
 \label{eq.cor.dumont.2}
\ee
where $T(t)$ is the tree function \reff{def.treefn}.
\end{corollary}

For completeness, let us outline briefly the elegant proof
of Proposition~\ref{prop.dumont}, due to Jiang Zeng,
that was presented in \cite[section~7]{Dumont_96}:

\sketchofproofof{Proposition~\ref{prop.dumont}}
(a) Consider a tree $T \in \scrt^{[1]}_{n+1}$,
and suppose that the root vertex~1 has $k$ ($\ge 0$) children.
All $k$ edges emanating from the root vertex are proper.
Deleting these edges and the vertex~1, one obtains a partition
of $\{2,\ldots,n+1\}$ into blocks $B_1,\ldots,B_k$
and a rooted tree $T_j$ on each block $B_j$.
Standard enumerative arguments then yield the relation~(a)
for the exponential generating functions.

(b) Consider a tree $T \in \scrt^{\<1\>}_{n+1}$ with root $r$,
and let $r_1, r_2, \ldots, r_l, 1$ ($l \ge 0$) be the path in $T$
from the root $r_1 = r$ to the leaf vertex~1.\footnote{
   Here $l = 0$ corresponds to the case in which the vertex~1
   is {\em both}\/ a leaf and the root
   (and hence the tree consists of just this one vertex).
}
All $l$ edges of this path are improper.
Deleting these edges and the vertex~1,
one obtains an {\em ordered}\/ partition of $\{2,\ldots,n+1\}$
into blocks $B_1,\ldots,B_l$
and a rooted tree $(T_j,r_j)$ on each block.
Standard enumerative arguments then yield the relation~(b)
for the exponential generating functions.

(c) In a tree $T \in \scrt_n$, focus on the vertex~1
(which might be the root, a leaf, both or neither).
Let $T'$ be the subtree rooted at~1,
and let $T''$ be the tree obtained from $T$
by deleting all the vertices of $T'$ except the vertex~1
(it thus has the vertex~1 as a leaf).
The vertex set $[n]$ is then partitioned as $\{1\} \cup V' \cup V''$,
where $\{1\} \cup V'$ is the vertex set of $T'$
and $\{1\} \cup V''$ is the vertex set of $T''$;
and $T$ is obtained by joining $T'$ and $T''$ at the common vertex~1.
Standard enumerative arguments then yield the relation~(c)
for the exponential generating functions.
\qed

\medskip

{\bf Remarks.}
1. Dumont and Ramamonjisoa also gave \cite[sections~2--5]{Dumont_96}
a second (and very interesting) proof of Proposition~\ref{prop.dumont},
based on a context-free grammar \cite{Chen_93}
and its associated differential operator.

2. We leave it as an open problem to find a direct combinatorial proof of the
functional equation \reff{eq.cor.dumont.1},
without using the differential equation of Proposition~\ref{prop.dumont}(d).

3. The polynomials $R_n$ also arise \cite{Josuat-Verges_15}
as derivative polynomials for the tree function:
in the notation of \cite{Josuat-Verges_15} we have $R_n(y,1) = G_n(y-1)$.
The formula \reff{eq.cor.dumont.2} is then equivalent to
\cite[Theorem~4.2, equation for~$G_n$]{Josuat-Verges_15}.
\myendremark

\subsection[The matrix $(f_{n,k}(y,\bphi))_{n,k \ge 0}$]{The matrix $\bm{(f_{n,k}(y,\bphi))_{n,k \ge 0}}$}
    \label{subsec.EGF.3}

Recall that $f_{n,k}(y,\bphi)$ is defined combinatorially
as the generating polynomial for $k$-component forests of rooted trees
on the vertex set $[n]$,
in which each improper edge gets a weight $y$
and each vertex with $m$~proper children
gets a weight $\phihat_m \eqdef m! \, \phi_m$.
The reasoning leading to the identity \reff{eq.EGF.star0}
again generalizes verbatim to $f_{n,k}(y,\bphi)$, so that
\be
   f_{n,k}(y,\bphi)
   \;=\;
   {1 \over k!}
   \!\!\!
   \sum\limits_{\begin{scarray}
                   n_1,\ldots,n_k \ge 1 \\
                   n_1 +\ldots+ n_k = n
                \end{scarray}}
   \!\!\!\!\!
   \binom{n}{n_1,\ldots,n_k} \,
       f_{n_1,1}(y,\bphi) \,\cdots\, f_{n_k,1}(y,\bphi)
   \;.
 \label{eq.EGF.star0bis.phi}
\ee
In terms of the column exponential generating functions
\be
   \scrf_k(t;y,\bphi)
   \;\eqdef\;
   \sum_{n=0}^\infty f_{n,k}(y,\bphi) \, {t^n \over n!}
   \;,
\ee
we have
\be
   \scrf_k(t;y,\bphi)  \;=\;  {\scrf_1(t;y,\bphi)^k \over k!}
   \;.
 \label{eq.EGF.star1bis.phi}
\ee
Therefore, the matrix $(f_{n,k}(y,\bphi))_{n,k \ge 0}$
is the exponential Riordan array $\scrr[F,G]$
with $F(t) = 1$ and $G(t) = \scrf_1(t;y,\bphi)$.

We now show how Proposition~\ref{prop.dumont} can be generalized
to incorporate the additional indeterminates $\bphi = (\phi_m)_{m \ge 0}$.
For a rooted tree $T$ on a totally ordered vertex set,
we define $\pc_m(T)$ to be the number of vertices of $T$
with $m$~proper children.
We define $\scrt^\bullet_n$, $\scrt^{[i]}_n$ and $\scrt^{\<i\>}_n$
as before, and then define the obvious generalizations of
\reff{def.Rn}--\reff{def.An.EGF}:
\begin{eqnarray}
   R_n(y,\bphi)  \;=\;  f_{n,1}(y,\bphi)
   & = &
   \sum_{T \in \scrt^\bullet_n} y^{\imp(T)}
      \prod_{m=0}^\infty (m! \, \phi_m)^{\pc_m(T)}
   \qquad\qquad
        \label{def.Rn.phi}  \\[2mm]
   S_n(y,\bphi)
   & = &
   \sum_{T \in \scrt^{[1]}_{n+1}} y^{\imp(T)}
      \prod_{m=0}^\infty (m! \, \phi_m)^{\pc_m(T)}
        \label{def.Sn.phi}  \\[2mm]
   A_n(y,\bphi)
   & = &
   \sum_{T \in \scrt^{\<1\>}_{n+1}} y^{\imp(T)}
      \prod_{m=0}^\infty (m! \, \phi_m)^{\pc_m(T)}
        \label{def.An.phi}
\end{eqnarray}
and
\begin{eqnarray}
   \scrr(t;y,\bphi)  \;=\;  \scrf_1(t;y,\bphi)
   & = &  \sum\limits_{n=1}^\infty R_n(y,\bphi) \: {t^n \over n!}
   \qquad\qquad\qquad
        \label{def.Rn.EGF.phi}  \\[2mm]
   \scrs(t;y,\bphi)  & = &  \sum\limits_{n=0}^\infty S_n(y,\bphi) \: {t^n \over n!}
        \label{def.Sn.EGF.phi}  \\[2mm]
   \scra(t;y,\bphi)  & = &  \sum\limits_{n=0}^\infty A_n(y,\bphi) \: {t^n \over n!}
        \label{def.An.EGF.phi}
\end{eqnarray}
Let us also define the generating function
\be
   \Phi(u)
   \;\eqdef\;
   \sum_{m=0}^\infty \phi_m \, u^m
   \;=\;
   \sum_{m=0}^\infty \phihat_m \, {u^m \over m!}
   \;.
 \label{def.Phi}
\ee
We then have:

\begin{proposition}
   \label{prop.dumont.phi}
The series $\scrr$, $\scrs$ and $\scra$
defined in \reff{def.Rn.EGF.phi}--\reff{def.An.EGF.phi}
satisfy the following identities:
\begin{itemize}
   \item[(a)]  $\scrs(t;y,\bphi) \:=\:
                  \Phi\big( \scrr(t;y,\bphi) \bigr)$
   \item[(b)]  $\scra(t;y,\bphi) \:=\:
                  \displaystyle{\phi_0 \over 1 - y \scrr(t;y,\bphi)}$
   \item[(c)]  $\displaystyle {d \over dt} \scrr(t;y,\bphi)  \:=\:
                  {1 \over \phi_0} \, \scra(t;y,\bphi) \, \scrs(t;y,\bphi)$
\end{itemize}
and hence
\begin{itemize}
   \item[(d)] $\displaystyle {d \over dt} \scrr(t;y,\bphi)
                             \:=\: 
               {\Phi\big( \scrr(t;y,\bphi) \bigr)
                \over
                1 - y \scrr(t;y,\bphi)}$
\end{itemize}
\end{proposition}

\proof
The proof is identical to that of Proposition~\ref{prop.dumont},
with the following modifications:

(a) Consider a tree $T \in \scrt^{[1]}_{n+1}$
in which the root vertex~1 has $k$ children.
Since all $k$ edges emanating from the root vertex are proper,
we get an additional factor $\phihat_k$
over and above what was seen in Proposition~\ref{prop.dumont}.
Therefore, the exponential function in Proposition~\ref{prop.dumont}
is replaced here by the generating function $\Phi$.

(b) Consider a tree $T \in \scrt^{\<1\>}_{n+1}$ with root $r$,
where $r_1, r_2, \ldots, r_l, 1$ is the path in $T$
from the root $r_1 = r$ to the leaf vertex~1.
Since all $l$ edges of this path are improper,
the weights associated to the vertices $r_1, r_2, \ldots, r_l$ in $T$
are identical to those associated to these vertices in the trees $(T_j,r_j)$;
therefore no modification is required.
However, the tree $T$ contains a leaf vertex~1
that is not present in any of the trees $(T_j,r_j)$,
so we get an additional factor $\phihat_0 = \phi_0$.

(c) In a tree $T \in \scrt_n$, focus on the vertex~1
and define $T'$ and $T''$ as before.
Since $T''$ has the vertex~1 as a leaf but $T$ does not,
a factor of $\phi_0$ needs to be removed from the right-hand side.
\qed

Let us give a name to the function appearing on the right-hand side
of the differential equation in Proposition~\ref{prop.dumont.phi}(d):
\be
   \Psi(s;y,\bphi)
   \;\eqdef\;
   {\Phi(s) \over 1 - ys}
   \;\eqdef\;
   \sum_{m=0}^\infty (\bphi * y^\N)_m \, s^m
 \label{def.Psi}
\ee
where $\bphi * y^\N$ is the convolution
\be
   (\bphi * y^\N)_m
   \;\eqdef\;
   \sum_{r=0}^m \phi_r \, y^{m-r}
   \;=\;
   \sum_{r=0}^m {\phihat_r \, y^{m-r}  \over r!}
   \;.
 \label{def.psim}
\ee
It follows from Proposition~\ref{prop.dumont.phi}(d) that
the generating function $\scrr(t;y,\bphi)$,
and hence the generic rooted-forest polynomials $f_{n,k}(y,\bphi)$,
depends on the indeterminates $y,\bphi$ only via the combination
$\bphi * y^\N$.
Otherwise put, if $\bphi',y'$  and $\bphi'',y''$
are two specializations of $y,\bphi$ to values in a commutative ring $R$
that satisfy $\bphi' * (y')^\N = \bphi'' * (y'')^\N$,
then $f_{n,k}(y',\bphi') = f_{n,k}(y'',\bphi'')$
for all $n,k \ge 0$.
We leave it as an open problem to find a bijective proof of this fact ---
possibly by bijection to a ``canonical'' specialization such as $y=0$,
i.e.\ a bijective proof of
\be
   f_{n,k}(y,\bphi)  \;=\;  f_{n,k}(0,\bphi * y^\N)
 \label{eq.fnk.phipsi}
\ee
(see also Section~\ref{sec.lah} below).

\bigskip

{\bf Remark.}
One might hope to generalize Proposition~\ref{prop.dumont.phi}
--- and thus also Theorem~\ref{thm1.4} ---
by refining the counting of improper edges, as follows:
Let $\bphi = (\phi_m)_{m \ge 0}$ and $\bxi = (\xi_\ell)_{\ell \ge 0}$
be indeterminates,
and let $f_{n,k}(\bxi,\bphi)$ be the generating polynomial for
$k$-component forests of rooted trees on the vertex set $[n]$
with a weight $m! \, \phi_m \, \xi_\ell$
for each vertex that has $m$ proper children and $\ell$ improper children.
Our polynomials $f_{n,k}(y,\bphi)$ thus correspond to the special case
$\xi_\ell = y^\ell$.
One might then hope that Proposition~\ref{prop.dumont.phi}
could be generalized to this case,
with $1/(1 - y\scrr)$ replaced by $\Xi(\scrr)$,
where $\Xi(u) = \sum_{\ell=0}^\infty \xi_\ell \, u^\ell$.
Indeed, Proposition~\ref{prop.dumont.phi}(a,c) do extend to this situation;
but Proposition~\ref{prop.dumont.phi}(b) does not,
because the ``global'' counting of improper edges implicit in the proof
does not correspond to the ``local'' counting of improper edges
(assigning them all to the parent vertex)
adopted in this definition of $f_{n,k}(\bxi,\bphi)$.
And in fact, the resulting polynomials are different:
the differential equation $\scrr'(t) = \Phi(\scrr) \, \Xi(\scrr)$
leads to
\be
   R_3(\bxi,\bphi)
   \;=\;
   \phi_1^2 \,+\, 4 \phi_1 \xi_1 \,+\, \xi_1^2
            \,+\, 2 \phi_2 \,+\, 2 \xi_2
   \;,
\ee
while the counting of the nine 3-vertex trees
with the specified weights yields
\be
   f_{3,1}(\bxi,\bphi)
   \;=\;
   \phi_1^2 \,+\, 4 \phi_1 \xi_1 \,+\, 2\xi_1^2
            \,+\, 2 \phi_2 \,+\, \xi_2
   \;.
\ee
The terms corresponding to trees with two improper edges are thus different:
$\xi_1^2 + 2 \xi_2$ from the differential equation,
and $2\xi_1^2 + \xi_2$ from the counting.

I leave it as an open problem to find a different way of ``localizing''
the improper edges that would provide a combinatorial interpretation
for the polynomials defined by
the differential equation $\scrr'(t) = \Phi(\scrr) \, \Xi(\scrr)$.
\myendremark

\section{Proof of Theorems~\ref{thm1.1}--\ref{thm1.4}}
     \label{sec.proofs}

We will prove Theorems~\ref{thm1.1}--\ref{thm1.4}
by explicitly exhibiting the production matrices
for $F$, $F(x)$, $F(x,y,z)$ and $F(x,y,\bphi)$
and then proving that these production matrices are
coefficientwise totally positive.
By Theorems~\ref{thm.iteration.homo} and \ref{thm.iteration2bis},
this will prove the claimed results.

It suffices of course to prove Theorem~\ref{thm1.4},
since Theorems~\ref{thm1.1}--\ref{thm1.3} are contained in it as special cases:
take $\phi_m = z^m/m!$ to get Theorem~\ref{thm1.3};
then take $y=z=1$ to get Theorem~\ref{thm1.2};
and finally take $x=0$ to get Theorem~\ref{thm1.1}.
However, we shall find it convenient to work our way up,
starting with Theorem~\ref{thm1.1} and then gradually adding extra parameters.

\subsection[The matrix $(f_{n,k})_{n,k \ge 0}$ and its production matrix]{The matrix $\bm{(f_{n,k})_{n,k \ge 0}}$ and its production matrix}

Let $F = (f_{n,k})_{n,k \ge 0}$ be the
unit-lower-triangular matrix defined by \reff{def.fnk}.
Straightforward computation gives for the first few rows
of its production matrix
\be
   P  \;\eqdef\;  F^{-1} \Delta F
   \;=\;
\Scale[0.95]{
   \begin{bmatrix*}[r]
 0  &   1  &      &      &      &      &      &       &  \\
 0  &   2  &   1  &      &      &      &      &       &  \\
 0  &   5  &   4  &   1  &      &      &      &       &  \\
 0  &   16  &   15  &   6  &   1  &      &      &       &  \\
 0  &   65  &   64  &   30  &   8  &   1  &      &       &  \\
 0  &   326  &   325  &   160  &   50  &   10  &   1  &       &  \\
 0  &   1957  &   1956  &   975  &   320  &   75  &   12  &   1   &  \\
 0  &   13700  &   13699  &   6846  &   2275  &   560  &   105  &   14   & \ddots  \\
 \vdots & \vdots & \vdots & \vdots & \vdots & \vdots & \vdots & \vdots &  \ddots
    \end{bmatrix*}
    \,.
 \label{eq.P.empirical}
}
\ee
Empirically this matrix seems to be \cite[A073107]{OEIS}
augmented by a column of zeros at the left.
Taking the explicit formula from \cite[A073107]{OEIS}
and inserting the extra column of zeros leads to the conjecture:

\begin{proposition}[Production matrix for $F$]
   \label{prop.prodmat.F}
Let $F = (f_{n,k})_{n,k \ge 0}$ be the
unit-lower-triangular matrix defined by \reff{def.fnk}.
Then its production matrix $P = (p_{jk})_{j,k \ge 0} = F^{-1} \Delta F$
has matrix elements
\be
   p_{jk}  \;=\;  \sum_{m=0}^j {j! \over m!} \, \binom{m}{k-1}
           \;=\;  {j! \over (k-1)!} \sum_{\ell=0}^{j+1-k} {1 \over \ell!}
 \label{eq.prop.prodmat.F}
\ee
(where $p_{j0} = 0$).
\end{proposition}

We will give two proofs of Proposition~\ref{prop.prodmat.F}:
a first proof using the theory of exponential Riordan arrays,
and a second proof by direct computation using Abel-type identities.

\firstproofof{Proposition~\ref{prop.prodmat.F}}
It was shown in Section~\ref{subsec.EGF.1}
that the matrix $(f_{n,k})_{n,k \ge 0}$
is the exponential Riordan array with $F(t) = 1$
and $G(t) = $ the tree function $T(t) = \sum_{n=1}^\infty n^{n-1} \, t^n/n!$.
Differentiation of the functional equation $T(t) = t \, e^{T(t)}$ gives
\be
   T'(t)  \;=\;  {e^{T(t)} \over 1 - T(t)}
   \;.
 \label{eq.star3}
\ee
Applying Theorem~\ref{thm.riordan.exponential.production},
we see by comparing \reff{eq.prop.riordan.exponential.production.1}
with \reff{eq.star3} that $Z = 0$ and $A(s) = e^s / (1-s)$,
which implies $z_n = 0$ and
\be
   a_n  \;=\;  \sum_{\ell=0}^n {1 \over \ell!}
   \;.
\ee
Inserting this into \reff{eq.thm.riordan.exponential.production}
yields \reff{eq.prop.prodmat.F}.
\qed

\secondproofof{Proposition~\ref{prop.prodmat.F}}
The production matrix $P = F^{-1} \Delta F$ satisfies the recurrence
[cf.\ \reff{def.iteration.bis}]
\be
   f_{n+1,k}  \;=\;  \sum_{j=0}^n  f_{n,j} \, p_{jk}
\ee
for all $k \ge 0$, or in other words
\be
   \binom{n+1}{k} \, k \, (n+1)^{n-k}
   \;=\;
   \sum_{j=0}^n \binom{n}{j} \, j \, n^{n-j-1}   \: p_{jk}
   \;.
 \label{eq.proof.prop.prodmat.F.1}
\ee
Now an Abel inverse relation
\cite[p.~95, eq.~(3) with $x=0$]{Riordan_68}
\cite[p.~154, Example~12 with $y=0$]{Roman_84}
states that
\be
   a_n  \;=\;  \sum_{j=0}^n \binom{n}{j} \, (-j)^{n-j} \: b_j
   \qquad\Longleftrightarrow\qquad
   b_n  \;=\;  \sum_{j=0}^n \binom{n}{j} \, j \, n^{n-j-1} \: a_j
   \;.
\ee
Applying this with $a_j = p_{jk}$
and $b_n = \binom{n+1}{k} \, k \, (n+1)^{n-k}$
at fixed $k \ge 0$,
we see that \reff{eq.proof.prop.prodmat.F.1} is equivalent to
\be
   p_{nk}
   \;=\;
   \sum_{j=0}^n \binom{n}{j} \, (-j)^{n-j} \: 
        \binom{j+1}{k} \, k \, (j+1)^{j-k}
   \;.
 \label{eq.proof.prop.prodmat.F.2}
\ee
A bit of algebra shows that the right-hand side of
\reff{eq.proof.prop.prodmat.F.2} can be rewritten as
\begin{subeqnarray}
   \hbox{RHS}
   & = &
   \binom{n}{k-1}
     \sum_{j=k-1}^n \binom{n+1-k}{j+1-k} \, (-j)^{n-j} \, (j+1)^{j+1-k}
          \\[2mm]
   & = &
   \binom{n}{k-1}
     \sum_{\ell=0}^N \binom{N}{\ell} \, (1-k-\ell)^{N-\ell} \, (k+\ell)^\ell
 \slabel{eq.proof.prop.prodmat.F.3.b}
\end{subeqnarray}
where $\ell = j+1-k$ and $N = n+1-k$.
But Cauchy's formula \cite[p.~21]{Riordan_68} implies that the
right-hand side of \reff{eq.proof.prop.prodmat.F.3.b} equals
\be
   \binom{n}{k-1}
   \sum_{\ell=0}^N \binom{N}{\ell} \, \ell! \: 1^{N-\ell}
   \;=\;
   {n! \over (k-1)!} \sum_{\ell=0}^{n+1-k} {1 \over (n+1-k-\ell)!}
   \;,
\ee
which equals $p_{nk}$ as defined in \reff{eq.prop.prodmat.F}.
\qed
   
\begin{corollary}[Production matrix for $F'$]
   \label{cor.prodmat.F}
Let $F' = (f_{n+1,k+1})_{n,k \ge 0} = \Delta F \Delta^{\rm T}$
be the unit-lower-triangular matrix obtained from $F$
by deleting its zeroth row and column.
Then its production matrix $P' = (p'_{jk})_{j,k \ge 0} = (F')^{-1} \Delta F'$
is obtained from $P$ by deleting its zeroth row and column,
i.e.\ $P' = \Delta P \Delta^{\rm T}$, and hence has matrix elements
\be
   p'_{jk}  \;=\;
   p_{j+1,k+1}  \;=\; {(j+1)! \over k!} \, \sum_{\ell=0}^{j+1-k} {1 \over \ell!}
   \;.
 \label{eq.cor.prodmat.F}
\ee
\end{corollary}

\proof
Apply Lemma~\ref{lemma.down-shifted} to Proposition~\ref{prop.prodmat.F},
with the matrix \reff{eq.P.empirical}/\reff{eq.prop.prodmat.F}
playing the role of $Q$.
\qed

We remark that the elements of $F'$ are
$f_{n+1,k+1} = \displaystyle \binom{n}{k} (n+1)^{n-k}$.

\bigskip

Let us now introduce the sequence $\bpsi = (\psi_m)_{m \ge 0}$
of positive rational numbers given by
\be
   \psi_m  \;=\;  \sum_{\ell=0}^m {1 \over \ell!}
   \;,
\ee
and the corresponding lower-triangular Toeplitz matrix $T_\infty(\bpsi)$:
\be
   T_\infty(\bpsi)_{ij}   \;=\;  \psi_{i-j}
\ee
with the convention $\psi_m \eqdef 0$ for $m < 0$.
Then the production matrix \reff{eq.prop.prodmat.F} can be written as
\be
    P  \;=\;   D T_\infty(\bpsi) D^{-1} \, \Delta
 \label{eq.prodmat.F.toeplitz}
\ee
where $D = \diag\big( (i!)_{i \ge 0} \big)$.
Moreover, this production matrix has a nice factorization into simpler matrices:

\begin{proposition}[Factorization of the production matrix]
   \label{prop.P.factorization}
The matrix $P = (p_{jk})_{j,k \ge 0}$ defined by \reff{eq.prop.prodmat.F}
has the factorization
\be
   P  \;=\;  B_1 \, D T_1 D^{-1} \, \Delta
 \label{eq.prop.P.factorization}
\ee
where $B_1$ is the binomial matrix [cf.\ \reff{def.Bx}],
$T_1$ is the lower-triangular matrix of all ones [cf.\ \reff{def.Tx}],
and $D = \diag\big( (i!)_{i \ge 0} \big)$.
\end{proposition}

\proof
We have $\bpsi = \ba * \bb$ where $a_n = 1/n!$ and $b_n = 1$,
and hence
\begin{subeqnarray}
   P
   & = &
   D T_\infty(\bpsi) D^{-1} \, \Delta  \\[1mm]
   & = &
   D T_\infty(\ba) T_\infty(\bb) D^{-1} \, \Delta  \\[1mm]
   & = &
   (D T_\infty(\ba) D^{-1}) \, (D T_\infty(\bb) D^{-1}) \, \Delta  \\[1mm]
   & = &
   B_1 \, D T_1 D^{-1} \, \Delta
   \;.
 \label{eq.P.factorization}
\end{subeqnarray}
\qed

%

{\bf Remarks.}
1.  Since $P' = \Delta P \Delta^{\rm T}$, this also implies
\be
   P' \;=\; \Delta \, B_1 \, D T_1 D^{-1}
   \;.
 \label{eq.prop.Pprime.factorization}
\ee

2. It follows from \reff{eq.prop.P.factorization}
that the augmented production matrix
$\displaystyle   
   \widetilde{P}
   \eqdef
   \left[
    \begin{array}{c}
         1 \; 0 \; 0 \; 0 \; \cdots \; \vphantom{\Sigma} \\
         \hline
         P
    \end{array}
    \right]
$
is given here by
\be
   \widetilde{P}
   \;=\;
   \left[
    \begin{array}{c|c}
         1 & \bzero  \\
         \hline
                &   \\[-4mm]
         \bzero & B_1 \, D T_1 D^{-1}
    \end{array}
    \right]
    \;.
 \vspace*{-3mm}
\ee
\myendremark

\bigskip

The sequence $\bpsi$ has the ordinary generating function
\be
   \Psi(s)  \;\eqdef\; \sum_{m=0}^\infty \psi_m s^m
       \;=\;  {e^s \over 1-s}
   \;.
 \label{eq.Psi}
\ee
Since this generating function is of the form \reff{eq.thm.aissen},
it follows that the sequence $\bpsi$ is Toeplitz-totally positive.
(This can equivalently be seen by observing that $\bpsi = \ba * \bb$,
 where $a_n = 1/n!$ and $b_n = 1$ are both Toeplitz-totally positive.)
In view of \reff{eq.prodmat.F.toeplitz}, this proves:

\begin{proposition}[Total positivity of the production matrix for $F$]
   \label{prop.prodmat.F.TP}
The matrix $P = (p_{jk})_{j,k \ge 0}$ defined by \reff{eq.prop.prodmat.F}
is totally positive (in $\Z$).
\end{proposition}

\begin{corollary}[Total positivity of the production matrix for $F'$]
   \label{cor.prodmat.F.TP}
The matrix $P' = \Delta P \Delta^{\rm T}$
defined by \reff{eq.cor.prodmat.F} is totally positive (in $\Z$).
\end{corollary}

Equivalently, we can observe that the total positivity of $P$ and $P'$
follows from the factorizations
\reff{eq.prop.P.factorization}/\reff{eq.prop.Pprime.factorization}
together with
Lemmas~\ref{lemma.toeplitz.power} and \ref{lemma.binomialmatrix.TP}.

\proofof{Theorem~\ref{thm1.1}}
Applying Theorem~\ref{thm.iteration.homo} to the matrix $P$
and using Propositions~\ref{prop.prodmat.F} and \ref{prop.prodmat.F.TP},
we deduce Theorem~\ref{thm1.1}(a).
Similarly, applying Theorem~\ref{thm.iteration2bis} to the matrix $P'$
and using Corollaries~\ref{cor.prodmat.F} and \ref{cor.prodmat.F.TP},
we deduce Theorem~\ref{thm1.1}(b).
\qed

\subsection[The matrix $(F_{n,k}(x))_{n,k \ge 0}$ and its production matrix]{The matrix $\bm{(F_{n,k}(x))_{n,k \ge 0}}$ and its production matrix}

We now turn our attention to the matrix $F(x) = (F_{n,k}(x))_{n,k \ge 0}$
of binomial partial row-generating polynomials defined by \reff{def.Fnk}.
The matrix factorization $F(x) = F B_x$  [cf.\ \reff{eq.FBx}]
implies, by Lemma~\ref{lemma.production.AB},
that the production matrix of $F(x)$ is $B_x^{-1} P B_x$,
where $P$ is the production matrix of $F$
as determined in the preceding subsection [cf.\ \reff{eq.prop.prodmat.F}]
and $B_x$ is the $x$-binomial matrix [cf.\ \reff{def.Bx}].
But Lemma~\ref{lemma.BxinvPBx} shows that
$B_x^{-1} P B_x = P (I + x \Delta^{\rm T})$.
This, together with Proposition~\ref{prop.prodmat.F.TP},
immediately implies:

\begin{proposition}[Total positivity of the production matrix for $F(x)$]
   \label{prop.prodmat.Fx.TP}
The matrix $B_x^{-1} P B_x$ defined by \reff{eq.prop.prodmat.F}
and \reff{def.Bx} is totally positive in the ring $\Z[x]$
equipped with the coefficientwise order.
\end{proposition}

\proofof{Theorem~\ref{thm1.2}}
Applying Theorem~\ref{thm.iteration.homo} to the matrix $B_x^{-1} P B_x$
and using Propositions~\ref{prop.prodmat.F} and \ref{prop.prodmat.Fx.TP},
we deduce Theorem~\ref{thm1.2}(a).
Similarly, applying Theorem~\ref{thm.iteration2bis}
to the matrix $B_x^{-1} P B_x$
and using Propositions~\ref{prop.prodmat.F} and \ref{prop.prodmat.Fx.TP},
we deduce Theorem~\ref{thm1.2}(b).
\qed

An equivalent way of stating this proof of Theorem~\ref{thm1.2}(b)
is that we have applied Theorem~\ref{thm.hankelTP.expriordan}
to the matrices $F$ and $P$.

\subsection[The matrices $(f_{n,k}(y,z))_{n,k \ge 0}$ and
      $(F_{n,k}(x,y,z))_{n,k \ge 0}$ and their production matrices]{The matrices $\bm{(f_{n,k}(y,z))_{n,k \ge 0}}$ and
      $\bm{(F_{n,k}(x,y,z))_{n,k \ge 0}}$ and their production matrices}
   \label{sec.proof.thm1.3}

We now generalize the results of the preceding two subsections
to include the indeterminates $y$ and $z$.
The key result is the following:

\begin{proposition}[Production matrix for $F(y,z)$]
   \label{prop.prodmat.Fyz}
Let $F(y,z) = (f_{n,k}(y,z))_{n,k \ge 0}$ be the
unit-lower-triangular matrix defined by \reff{def.fnkyz}.
Then its production matrix
$P(y,z) = (p_{nk}(y,z))_{n,k \ge 0} = F(y,z)^{-1} \Delta F(y,z)$
has matrix elements
\be
   p_{nk}(y,z)
   \;=\;
   {n! \over (k-1)!} \sum_{\ell=0}^{n+1-k} {y^{n-\ell} z^\ell \over \ell!}
   \;.
 \label{eq.prop.prodmat.Fyz}
\ee
\end{proposition}

This time we have only a proof using exponential Riordan arrays:

\proofof{Proposition~\ref{prop.prodmat.Fyz}}
It was shown in Section~\ref{subsec.EGF.2}
that the matrix
\linebreak
$(f_{n,k}(y,z))_{n,k \ge 0}$
is the exponential Riordan array with $F(t) = 1$ and $G(t) = \scrr(t;y,z)$,
where $\scrr(t;y,z)$ solves the differential equation
of Proposition~\ref{prop.dumont}(d) with initial condition $\scrr(0;y,z) = 0$.
Applying Theorem~\ref{thm.riordan.exponential.production}
and comparing this differential equation with
\reff{eq.prop.riordan.exponential.production.1},
we see that $Z(s) = 0$ and $A(s) = e^{zs} / (1- ys)$,
which implies $z_n = 0$ and
\be
   a_n  \;=\;  \sum_{\ell=0}^n {y^{n-\ell} z^\ell \over \ell!}
   \;.
\ee
Inserting this into \reff{eq.thm.riordan.exponential.production}
yields \reff{eq.prop.prodmat.Fyz}.
\qed

Let us now introduce the sequence $\bpsi(y,z) = (\psi_m(y,z))_{m \ge 0}$
of polynomials with nonnegative rational coefficients given by
\be
   \psi_m(y,z)  \;=\;  \sum_{\ell=0}^m {y^{n-\ell} z^\ell \over \ell!}
   \;,
\ee
and the corresponding lower-triangular Toeplitz matrix $T_\infty(\bpsi(y,z))$.
Then \reff{eq.prop.prodmat.Fyz} can be written as
\be
   P(y,z)  \;=\;  D T_\infty(\bpsi(y,z)) D^{-1} \, \Delta
 \label{eq.prodmat.Fyz.toeplitz}
\ee
where $D = \diag\big( (i!)_{i \ge 0} \big)$;
the elements of these matrices lie in the ring $\Q[y,z]$.
Moreover, this production matrix has a nice factorization into simpler matrices:

\begin{proposition}[Factorization of the production matrix]
   \label{prop.P.factorization.Fyz}
The matrix $P = (p_{jk})_{j,k \ge 0}$ defined by \reff{eq.prop.prodmat.Fyz}
has the factorization
\be
   P(y,z)  \;=\;  B_z \, D T_y D^{-1} \, \Delta
 \label{eq.prop.P.factorization.Fyz}
\ee
where $B_z$ is the weighted binomial matrix \reff{def.Bx},
$T_y$ is the Toeplitz matrix of powers \reff{def.Tx},
and $D = \diag\big( (i!)_{i \ge 0} \big)$.
\end{proposition}

\proof
We have $\bpsi(y,z) = \ba * \bb$ where $a_n = z^n/n!$ and $b_n = y^n$,
and hence
\begin{subeqnarray}
   P(y,z)
   & = &
   D T_\infty(\bpsi(y,z)) D^{-1} \, \Delta  \\[1mm]
   & = &
   D T_\infty(\ba) T_\infty(\bb) D^{-1} \, \Delta  \\[1mm]
   & = &
   (D T_\infty(\ba) D^{-1}) \, (D T_\infty(\bb) D^{-1}) \, \Delta  \\[1mm]
   & = &
   B_z \, D T_y D^{-1} \, \Delta
   \;.
 \label{eq.P.factorization.Fyz}
\end{subeqnarray}
\qed

{\bf Remark.}
It follows from \reff{eq.prop.P.factorization.Fyz}
that the augmented production matrix
$\displaystyle   
   \widetilde{P}(y,z)
   =
   \left[
    \begin{array}{c}
         1 \; 0 \; 0 \; 0 \; \cdots \; \vphantom{\Sigma} \\
         \hline
         P(y,z)
    \end{array}
    \right]
$
is given by
\be
   \widetilde{P}(y,z)
   \;=\;
   \left[
    \begin{array}{c|c}
         1 & \bzero  \\
         \hline
                &   \\[-4mm]
         \bzero & B_z \, D T_y D^{-1}
    \end{array}
    \right]
    \;.
   \label{eq.augmented.Pyz}
\ee
Two interpretations of \reff{eq.augmented.Pyz}/\reff{eq.prodmat.u.iterated}
in terms of digraphs are given by Gilmore \cite{Gilmore_inprep}.
\myendremark

\bigskip

The sequence $\bpsi(y,z)$ has the ordinary generating function
\be
   \Psi(s;y,z)  \;\eqdef\; \sum_{m=0}^\infty \psi_m(y,z) \, s^m
       \;=\;  {e^{zs} \over 1-ys}
   \;.
 \label{eq.Psi.yz}
\ee
Since this generating function is of the form \reff{eq.thm.toeplitz.1a},
Lemma~\ref{lemma.toeplitz.1a} implies that
the sequence $\bpsi$ is coefficientwise Toeplitz-totally positive.
(This can equivalently be seen by observing that $\bpsi(y,z) = \ba * \bb$,
 where $a_n = z^n/n!$ and $b_n = y^n$ are both
 coefficientwise Toeplitz-totally positive.)
In other words, the Toeplitz matrix $T_\infty(\bpsi(y,z))$
is totally positive in the ring $\Q[y,z]$
equipped with the coefficientwise order.
It follows from \reff{eq.prodmat.Fyz.toeplitz} that the same goes for $P(y,z)$.
But the elements of $P(y,z)$ actually lie in the ring
$\Z[y,z] \subseteq \Q[y,z]$.  We have therefore proven:

\begin{proposition}[Total positivity of the production matrix for $F(y,z)$]
   \label{prop.prodmat.Fyz.TP}
The matrix $P(y,z) = (p_{jk}(y,z))_{j,k \ge 0}$
defined by \reff{eq.prop.prodmat.Fyz}
is totally positive in the ring $\Z[y,z]$
equipped with the coefficientwise order.
\end{proposition}

Equivalently, the total positivity of $P(y,z)$
follows from the factorization \reff{eq.prop.P.factorization.Fyz}
together with
Lemmas~\ref{lemma.toeplitz.power} and \ref{lemma.binomialmatrix.TP}.

\bigskip

We now consider the matrix $F(x,y,z) = (F_{n,k}(x,y,z))_{n,k \ge 0}$
of binomial partial row-generating polynomials defined by \reff{def.Fnkyz}.
The matrix factorization $F(x,y,z) = F(y,z) B_x$
implies, by Lemma~\ref{lemma.production.AB},
that the production matrix of $F(x,y,z)$ is
\linebreak
$B_x^{-1} P(y,z) B_x$,
where $P(y,z)$ is the production matrix of $F(y,z)$
[cf.\ \reff{eq.prop.prodmat.Fyz}].
But Lemma~\ref{lemma.BxinvPBx} shows that
$B_x^{-1} P(y,z) B_x = P(y,z) (I + x \Delta^{\rm T})$.
This, together with Proposition~\ref{prop.prodmat.Fyz.TP},
immediately implies:

\begin{proposition}[Total positivity of the production matrix for $F(x,y,z)$]
   \label{prop.prodmat.Fxyz.TP}
The matrix $B_x^{-1} P(y,z) B_x$ defined by \reff{eq.prop.prodmat.Fyz}
and \reff{def.Bx} is totally positive in the ring $\Z[x,y,z]$
equipped with the coefficientwise order.
\end{proposition}

\proofof{Theorem~\ref{thm1.3}}
Applying Theorem~\ref{thm.iteration.homo} to the matrix $B_x^{-1} P(y,z) B_x$
and using Propositions~\ref{prop.prodmat.Fyz} and \ref{prop.prodmat.Fxyz.TP},
we deduce Theorem~\ref{thm1.3}(a).

Similarly, applying Theorem~\ref{thm.iteration2bis}
to the matrix $B_x^{-1} P(y,z) B_x$
and using Propositions~\ref{prop.prodmat.Fyz} and \ref{prop.prodmat.Fxyz.TP},
we deduce Theorem~\ref{thm1.3}(b).

Theorem~\ref{thm1.3}(c) follows from Theorem~\ref{thm1.3}(b)
by noting that
\be
   f_{n+1,1}(y,z)
   \;=\;
   \left. {F_{n+1}(x,y,z) \over x} \right|_{x=0}
   \;.
\ee 
\qed

Once again, an equivalent way of stating this proof of Theorem~\ref{thm1.3}(b)
is that we have applied Theorem~\ref{thm.hankelTP.expriordan}
to the matrices $F(y,z)$ and $P(y,z)$.

%
%

\subsection[The matrices $(f_{n,k}(y,\bphi))_{n,k \ge 0}$ and
      $(F_{n,k}(x,y,\bphi))_{n,k \ge 0}$ and their production matrices]{The matrices $\bm{(f_{n,k}(y,\bphi))_{n,k \ge 0}}$ and
      $\bm{(F_{n,k}(x,y,\bphi))_{n,k \ge 0}}$ and their production matrices}

We now generalize the results of the preceding three subsections
to include the indeterminates $\bphi$.

\begin{proposition}[Production matrix for $F(y,\bphi)$]
   \label{prop.prodmat.Fyz.phi}
The lower-triangular matrix
$F(y,\bphi) = (f_{n,k}(y,\bphi))_{n,k \ge 0}$
has production matrix
$P(y,\bphi) = (p_{nk}(y,\bphi))_{n,k \ge 0}$
\linebreak
 $= F(y,\bphi)^{-1} \Delta F(y,\bphi)$
given by
\be
   p_{nk}(y,\bphi)
   \;=\;
   {n! \over (k-1)!} \: (\bphi * y^\N)_{n-k+1}
   \;,
 \label{eq.prop.prodmat.Fyz.phi}
\ee
where
\be
   (\bphi * y^\N)_m
   \;\eqdef\;
   \sum_{\ell=0}^m \phi_\ell \, y^{m-\ell}
   \;=\;
   \sum_{\ell=0}^m {\phihat_\ell \, y^{m-\ell}  \over \ell!}
   \;.
 \label{def.psim.bis}
\ee
\end{proposition}

\proof
It was shown in Section~\ref{subsec.EGF.3}
that the matrix $(f_{n,k}(y,\bphi))_{n,k \ge 0}$
is the exponential Riordan array with $F(t) = 1$
and $G(t) = \scrr(t;y,\bphi)$,
where $\scrr(t;y,\bphi)$ solves the differential equation
of Proposition~\ref{prop.dumont.phi}(d)
with initial condition $\scrr(0;y,\bphi) = 0$.
Applying Theorem~\ref{thm.riordan.exponential.production}
and comparing this differential equation with
\reff{eq.prop.riordan.exponential.production.1},
we see that $Z(s) = 0$ and $A(s) = \Psi(s;y,\bphi)$
as defined in \reff{def.Psi},
which implies that $z_n = 0$ and $a_n = (\bphi * y^\N)_n$.
Inserting this into \reff{eq.thm.riordan.exponential.production}
yields \reff{eq.prop.prodmat.Fyz.phi}.
\qed

Now suppose that $\bphi$ is specialized to be a sequence,
with values in a partially ordered commutative ring $R$,
that is Toeplitz-totally positive of order $r$.
Then the sequence $\bphi$ is obviously Toeplitz-TP${}_r$ in the ring $R[y]$
equipped with the coefficientwise order.
And  by Lemma~\ref{lemma.toeplitz.power},
the sequence $y^\N \eqdef (y^n)_{n \ge 0}$ is Toeplitz-TP
in the ring $R[y]$ equipped with the coefficientwise order.
It follows that their convolution $\bphi * y^\N$
is Toeplitz-TP${}_r$ in the ring $R[y]$
equipped with the coefficientwise order.
On the other hand, \reff{eq.prop.prodmat.Fyz.phi} can be written as
\be
   P(y,\bphi) \;=\;  T_\infty(\bphi * y^\N)^\sharp \,  \Delta
   \;,
\ee
where the operation ${}^\sharp$ is defined in
Section~\ref{subsec.diagonal_scaling}.
Lemma~\ref{lemma.diagmult.TP} then implies that
the matrix $P(y,\bphi)$ is TP${}_r$
in the ring $R[y]$ equipped with the coefficientwise order.
We have therefore proven:

\begin{proposition}[Total positivity of the production matrix for $F(y,\bphi)$]
   \label{prop.prodmat.Fyz.TP.phi}
\quad\hfill\break\hspace*{-1.7mm}
Fix $1 \le r \le \infty$.
Let $R$ be a partially ordered commutative ring,
and let $\bphi = (\phi_m)_{m \ge 0}$ be a sequence in $R$
that is Toeplitz-totally positive of order~$r$.
Then the matrix $P(y,\bphi) = (p_{jk}(y,\bphi))_{j,k \ge 0}$
defined by \reff{eq.prop.prodmat.Fyz.phi}
is totally positive of order $r$ in the ring $R[y]$
equipped with the coefficientwise order.
\end{proposition}

\begin{remark}
\rm
If the ring $R$ contains the rationals, then we have the factorization
\begin{subeqnarray}
   P(y,\bphi)
   & = &
   D \, T_\infty(\bphi * y^\N) \, D^{-1} \,  \Delta
       \\[1mm]
   & = &
   D \, T_\infty(\bphi) \,
        T_\infty(y^\N) \, D^{-1} \, \Delta
       \\[1mm]
   & = &
   D \, T_\infty(\bphi) \, D^{-1} \:
   D \, T_y \, D^{-1} \: \Delta
\end{subeqnarray}
where $D = \diag\big( (i!)_{i \ge 0} \big)$,
by analogy with Propositions~\ref{prop.P.factorization}
and \ref{prop.P.factorization.Fyz}.
\myendremark
\end{remark}

\bigskip

We now consider the matrix $F(x,y,\bphi) = (F_{n,k}(x,y,\bphi))_{n,k \ge 0}$
of binomial partial row-generating polynomials defined by \reff{def.Fnkyz.phi}.
The matrix factorization
\linebreak
$F(x,y,\bphi) = F(y,\bphi) B_x$
implies, by Lemma~\ref{lemma.production.AB},
that the production matrix of $F(x,y,\bphi)$ is $B_x^{-1} P(y,\bphi) B_x$,
where $P(y,\bphi)$ is the production matrix of $F(y,\bphi)$
[cf.\ \reff{eq.prop.prodmat.Fyz.phi}].
But Lemma~\ref{lemma.BxinvPBx} shows that
$B_x^{-1} P(y,\bphi) B_x =$
${P(y,\bphi) (I + x \Delta^{\rm T})}$.
This, together with Proposition~\ref{prop.prodmat.Fyz.TP.phi},
immediately implies:

\begin{proposition}[Total positivity of the production matrix for $F(x,y,\bphi)$]
   \label{prop.prodmat.Fxyz.TP.phi}
\quad\hfill\break
Fix $1 \le r \le \infty$.
Let $R$ be a partially ordered commutative ring,
and let $\bphi = (\phi_m)_{m \ge 0}$ be a sequence in $R$
that is Toeplitz-totally positive of order~$r$.
Then the matrix $B_x^{-1} P(y,\bphi) B_x$
defined by \reff{eq.prop.prodmat.Fyz.phi} and \reff{def.Bx}
is totally positive of order $r$ in the ring $R[x,y]$
equipped with the coefficientwise order.
\end{proposition}

\proofof{Theorem~\ref{thm1.4}}
Applying Theorem~\ref{thm.iteration.homo}
to the matrix $B_x^{-1} P(y,\bphi) B_x$
and using Propositions~\ref{prop.prodmat.Fyz.phi}
and \ref{prop.prodmat.Fxyz.TP.phi},
we deduce Theorem~\ref{thm1.4}(a).

Similarly, applying Theorem~\ref{thm.iteration2bis}
to the matrix $B_x^{-1} P(y,\bphi) B_x$
and using Propositions~\ref{prop.prodmat.Fyz.phi}
and \ref{prop.prodmat.Fxyz.TP.phi},
we deduce Theorem~\ref{thm1.4}(b).

Then Theorem~\ref{thm1.4}(c) follows from Theorem~\ref{thm1.4}(b)
by noting that
\be
   f_{n+1,1}(y,\bphi)
   \;=\;
   \left. {F_{n+1}(x,y,\bphi) \over x} \right|_{x=0}
   \;.
 \vspace*{-3mm}
\ee
\qed

Once again, an equivalent way of stating this proof of Theorem~\ref{thm1.4}(b)
is that we have applied Theorem~\ref{thm.hankelTP.expriordan}
to the matrices $F(y,\bphi)$ and $P(y,\bphi)$.

\section{Connection with the generic Lah polynomials}  \label{sec.lah}

In a recent paper \cite{latpath_lah} we introduced the
{\em generic Lah polynomials}\/, which are defined as follows:

Recall first \cite[pp.~294--295]{Stanley_86}
that an {\em ordered tree}\/ (also called {\em plane tree}\/)
is a rooted tree in which the children of each vertex are linearly ordered.
An {\em unordered forest of ordered trees}\/
is an unordered collection of ordered trees.
An {\em increasing ordered tree}\/ is an ordered tree
in which the vertices carry distinct labels from a linearly ordered set
(usually some set of integers) in such a way that
the label of each child is greater than the label of its parent;
otherwise put, the labels increase along every path downwards from the root.
An {\em unordered forest of increasing ordered trees}\/
is an unordered forest of ordered trees with the same type of labeling.

Now let $\bphi = (\phi_m)_{m \ge 0}$ be indeterminates,
and let $L_{n,k}(\bphi)$ be the generating polynomial for
unordered forests of increasing ordered trees on the vertex set $[n]$,
having $k$ components (i.e.\ $k$ trees),
in which each vertex with $m$ children gets a weight $\phi_m$.
Clearly $L_{n,k}(\bphi)$ is a homogeneous polynomial of degree $n$
with nonnegative integer coefficients;
it is also quasi-homogeneous of degree $n-k$
when $\phi_m$ is assigned weight~$m$.
The first few polynomials $L_{n,k}(\bphi)$
[specialized for simplicity to $\phi_0 = 1$] are
%
%
%
\begin{table}[H]
\centering
\footnotesize
\begin{tabular}{c|cccccc}
$n \setminus k$ & 0 & 1 & 2 & 3 & 4 & 5 \\
\hline
0 & 1 &  &  &  &  &  \\
1 & 0 & 1 &  &  &  &  \\
2 & 0 & $\phi_1$ & 1 &  &  &  \\
3 & 0 & $\phi_1^2 + 2 \phi_2$ & $3 \phi_1$ & 1 &  &  \\
4 & 0 & $\phi_1^3 + 8 \phi_1 \phi_2 + 6 \phi_3$ & $7 \phi_1^2 + 8 \phi_2$ & $6 \phi_1$ & 1 &  \\
5 & 0 & $\phi_1^4 + 22 \phi_1^2 \phi_2 + 16 \phi_2^2 + 42 \phi_1 \phi_3 + 24 \phi_4$ & $15 \phi_1^3 + 60 \phi_1 \phi_2 + 30 \phi_3$ & $25 \phi_1^2 + 20 \phi_2$ & $10 \phi_1$ & 1 \\
\end{tabular}
\end{table}
\vspace*{-5mm}
\noindent
Now let $x$ be an additional indeterminate,
and define the row-generating polynomials
$L_n(\bphi,x) = \sum_{k=0}^n L_{n,k}(\bphi) \: x^k$.
Then $L_n(\bphi,x)$ is quasi-homogeneous of degree~$n$
when $\phi_i$ is assigned weight~$i$ and $x$ is assigned weight~1.
We call $L_{n,k}(\bphi)$ and $L_n(\bphi,x)$
the \textbfit{generic Lah polynomials},
and we call the lower-triangular matrix $\sfL = (L_{n,k}(\bphi))_{n,k \ge 0}$
the \textbfit{generic Lah triangle}.
Here $\bphi = (\phi_i)_{i \ge 0}$ are in the first instance indeterminates,
so that $L_{n,k}(\bphi) \in \Z[\bphi]$ and $L_n(\bphi,x) \in \Z[\bphi,x]$; 
but we can then, if we wish, substitute specific values for $\bphi$
in any commutative ring $R$,
leading to values $L_{n,k}(\bphi) \in R$ and $L_n(\bphi,x) \in R[x]$.

We can relate the generic Lah polynomials $L_{n,k}(\bphi)$
to the generic rooted-forest polynomials $f_{n,k}(y,\bphi)$, as follows:
First of all, the fact that we chose to define the generic Lah polynomials
in terms of {\em ordered}\/ trees is unimportant.
Since the vertices of our trees are labeled,
the children of each vertex are distinguishable;
therefore, for each unordered labeled tree and each vertex with $m$ children,
there are $m!$ possible orderings of those children.
It follows that the generic Lah polynomials,
defined initially as a sum over
unordered forests of increasing {\em ordered}\/ trees
with a weight $\phi_m$ for each vertex with $m$ children,
can equivalently be defined as a sum over
unordered forests of increasing {\em unordered}\/ trees
with a weight $\phihat_m = m! \, \phi_m$ for each vertex with $m$ children.
(This is why we inserted the factors $m!$ into our definition of
the generic rooted-forest polynomials.)
We shall henceforth reinterpret the generic Lah polynomials in this manner,
as a sum over unordered forests of increasing {\em unordered}\/ trees.

Now, the generic Lah polynomials are defined as a sum over
forests of {\em increasing}\/ trees on the vertex set $[n]$,
while the generic rooted-forest polynomials are defined as a sum over
forests of {\em arbitrary}\/ trees on the vertex set $[n]$
with a weight $y$ for each improper edge.
Furthermore, the generic Lah polynomials are defined as giving a weight
$\phihat_m$ for each vertex with $m$ children,
while the generic rooted-forest polynomials are defined as giving a weight
$\phihat_m$ for each vertex with $m$ {\em proper}\/ children.
But a tree is increasing if and only if all its edges are proper!
Therefore, by setting $y=0$ in the generic rooted-forest polynomials,
we ensure that the sum runs precisely over forests of {\em increasing}\/ trees,
and we also ensure that all the children at each vertex are proper.
It follows that the generic Lah polynomials $L_{n,k}(\bphi)$
are equal to the generic rooted-forest polynomials $f_{n,k}(y,\bphi)$
specialized to $y=0$:

\begin{proposition}[Generic Lah polynomials as specialization of generic rooted-forest polynomials]
   \label{prop.generic_lah}
We have $L_{n,k}(\bphi) = f_{n,k}(\bphi,0)$.
\end{proposition}

In \cite[Proposition~1.4]{latpath_lah} we showed that the
production matrix $P = (p_{ij})_{i,j \ge 0}$ for the generic Lah triangle
$\sfL = (L_{n,k}(\bphi))_{n,k \ge 0}$ is
\be
   p_{ij}
   \;=\;
   \begin{cases}
      0   &  \textrm{if $j=0$ or $j > i+1$}  \\[1mm]
      {\displaystyle {i! \over (j-1)!}} \, \phi_{i-j+1}
          &  \textrm{if $1 \le j \le i+1$}
   \end{cases}
 \label{eq.prop.prodmat.lah}
\ee
This is precisely Proposition~\ref{prop.prodmat.Fyz.phi}
of the present paper specialized to $y=0$.

So the generic rooted-forest polynomials are a generalization of
the generic Lah polynomials, to which they reduce when $y=0$.
On the other hand, the generic rooted-forest polynomials
are {\em also}\/ a {\em specialization}\/ of the generic Lah polynomials,
since \reff{eq.fnk.phipsi} and Proposition~\ref{prop.generic_lah}
immediately imply:

\begin{proposition}[Generic rooted-forest polynomials as specialization of generic Lah polynomials]
   \label{prop.generic_lah.2}
We have $f_{n,k}(y,\bphi) = L_{n,k}(\bphi * y^\N)$.
\end{proposition}

\noindent
We leave it as an open problem to find a direct (ideally bijective)
proof of Proposition~\ref{prop.generic_lah.2}.

Proposition~\ref{prop.generic_lah.2} can also be interpreted
in the language of exponential Riordan arrays.
As remarked in \cite[Section~8]{latpath_lah},
the generic Lah triangle $\sfL = (L_{n,k}(\bphi))_{n,k \ge 0}$
is in fact the {\em general}\/ exponential Riordan array $\scrr[F,G]$
of the ``associated subgroup'' $F=1$,
expressed in terms of its $A$-sequence $\ba = \bphi$
(cf.\ Theorem~\ref{thm.riordan.exponential.production}).
That is, the theory of the generic Lah triangle
is {\em equivalent}\/ to the theory of
exponential Riordan arrays of the ``associated subgroup'' $\scrr[1,G]$.
So, since the generic rooted-forest triangle
is indeed an exponential Riordan array of the associated subgroup
(Section~\ref{subsec.EGF.3}),
it {\em must}\/ be a specialization of the generic Lah triangle.

Let us remark, finally, that \cite[Section~3.1]{latpath_lah}
introduced a generalization of the generic Lah triangle
--- called the \textbfit{refined generic Lah triangle} ---
in which the weight for a vertex with $m$ children
now depends also on a quantity called its ``level'' $L$
\cite[Definition~3.1]{latpath_lah}.
The production matrix of the refined generic Lah triangle
was determined in \cite[Proposition~3.2]{latpath_lah};
the proof employed a bijection from
ordered forests of increasing ordered trees
to a set of labeled reversed partial \L{}ukasiewicz paths.
It would be interesting to know whether that construction
can be generalized to $y \neq 0$,
i.e.\ to forests of trees that are not necessarily increasing.

\section{Open problems}  \label{sec.open}

We conclude by proposing some open problems,
which are variants or generalizations of the results found here.

\subsection{Schl\"afli--Gessel--Seo polynomials}

In 1847, Schl\"afli \cite{Schlafli_1847} introduced the polynomials
\begin{subeqnarray}
   P_0(x;a,b)   & = &  1  \\
   P_n(x;a,b)   & = &  x \, \prod\limits_{i=1}^{n-1} [x + ia + (n-i)b]
     \qquad\hbox{for $n \ge 1$}
 \label{def.gessel-seo}
\end{subeqnarray}
and showed, using Lagrange inversion,
that their exponential generating function
\be
   \bfF(t;x,a,b)
   \;\eqdef\;
   \sum_{n=0}^\infty P_n(x;a,b) \, {t^n \over n!}
\ee
satisfies the functional equation
\be
   {\bfF^{-a/x} \,-\, \bfF^{-b/x}  \over  b-a}
   \;=\;
   t
   \;.
 \label{eq.gessel-seo.functeqn}
\ee
This immediately implies that
\be
   \bfF(t;x,a,b)  \;=\;  \bfF(t;1,a,b)^x
\ee
and hence \cite[eqns.~(4.1)/(4.2)]{Gould_74}
that the polynomials $P_n(x;a,b)$ form a {\em sequence of binomial type}\/
\cite{Mullin_70,Rota_73,Garsia_73,Fillmore_73,Labelle_81,Roman_78,Roman_84}:
that is,
\be
   P_n(x+y;a,b)
   \;=\;
   \sum_{k=0}^n \binom{n}{k} P_k(x;a,b) \, P_{n-k}(y;a,b)
   \;.
 \label{eq.binomial_type}
\ee
An equivalent statement is that, if we define
\be
   P_{n,k}(a,b)  \;=\;  [x^k] \, P_n(x;a,b)
     \quad\hbox{for $n,k \ge 0$}
   \;,
 \label{def.gessel-seo.Pnk}
\ee
the unit-lower-triangular matrix
$\bigl( P_{n,k}(a,b) \bigr)_{n,k \ge 0}$
is an exponential Riordan array $\scrr[F,G]$ with $F(t) = 1$.
The identity \reff{eq.binomial_type} goes back in fact
(in a slightly different notation)
to Rothe \cite{Rothe_1793} in 1793 and Pfaff \cite{Pfaff_1795} in 1795.
Rothe's identity is usually expressed in terms of the polynomials
\begin{subeqnarray}
   R_0(x;h,w)   & = &  1  \\
   R_n(x;h,w)   & = &  x \, \prod\limits_{i=1}^{n-1} (x + ih + nw)
     \qquad\hbox{for $n \ge 1$}
 \label{def.rothe}
\end{subeqnarray}
which obviously satisfy
\be
   R_n(x;h,w)  \;=\;  P_n(x;h+w,w)
   \;.
\ee
(But Schl\"afli's more symmetric formulation has a cleaner
combinatorial interpretation, as we shall see.)
The Rothe--Pfaff--Schl\"afli identity \reff{eq.binomial_type}
is thus a two-parameter identity that includes
as special cases the binomial theorem ($a=b=0$),
a variant of the Chu--Vandermonde identity ($a=0$ or $b=0$),
and a variant of Abel's \cite{Abel_1826} 1826 generalization
of the binomial theorem ($a=b$).
See \cite{Gould_56,Gould_57,Gould_66,Gould_74,Riordan_68,Strehl_92,Johnson_07,%
Scott-Sokal_expidentities,Sokal_chu-vandermonde}
for further discussion.

It is worth observing that the polynomials $P_n$ and $P_{n,k}$
are symmetric in ${a \leftrightarrow b}$;
that $P_n$ is homogeneous of degree $n$ in $x,a,b$;
and that $P_{n,k}$ is homogeneous of degree $n-k$ in $a,b$.
Note also that $P_{n,1}(a,b) = \prod\limits_{i=1}^{n-1} [ia + (n-i)b] \:$
\cite[A067948]{OEIS} \cite{Khidr_84}
and $P_{n,n-1}(a,b) = \binom{n}{2} (a+b)$.
The triangular array $\bigl( P_{n,k}(a,b) \bigr)_{n,k \ge 0}$ begins
\be
\Scale[0.7]{
\begin{bmatrix}
 1 &   &   &   &   &      &   \\
 0 & 1 &   &   &   &      &   \\
 0 & a+b & 1 &   &   &      &   \\
 0 & 2 a^2+5 a b+2 b^2 & 3 a+3 b & 1 &   &      &   \\
 0 & 6 a^3+26 a^2 b+26 a b^2+6 b^3 & 11 a^2+26 a b+11 b^2 & 6 a+6 b & 1 &      &   \\
 0 & 24 a^4+154 a^3 b+269 a^2 b^2+154 a b^3+24 b^4 & 50 a^3+200 a^2 b+200 a b^2+50 b^3 & 35 a^2+80 a b+35 b^2 & 10 a+10 b & 1    &   \\
 \vdots & \vdots & \vdots & \vdots & \vdots & \vdots &  \ddots
\end{bmatrix}
}
\ee

Furthermore --- and most importantly ---
we see from \reff{eq.Fn.explicit}/\reff{def.gessel-seo}
that $P_n(x;1,1) = x(x+n)^{n-1} = F_n(x)$
and hence that $P_{n,k}(1,1) = f_{n,k}$ [cf.~\reff{def.fnk}].
It follows that the polynomials $P_{n,k}(a,b)$
enumerate forests of rooted trees on the vertex set $[n]$ with $k$ components
according to some bivariate statistic.

%

Some years after Rothe, Pfaff and Schl\"afli
--- in 2006, to be precise ---
Gessel and Seo \cite{Gessel-Seo_06}, in a very interesting paper,
reintroduced the polynomials \reff{def.gessel-seo}
as enumerators of forests of rooted trees
and gave {\em two}\/ versions of this bivariate statistic,
as follows\footnote{
   I have altered their notation slightly: my $x$ is their $c$.
   I have also introduced explicitly the polynomials $P_{n,k}$;
   they introduced only $P_{n,1}$, calling it $Q_n$.
}:
Recall first that an edge $e = ij$ in a forest $\scrf$,
ordered so that $j$ is a child of~$i$,
is called a {\em proper edge}\/ if all the descendants of~$j$,
including $j$ itself, are higher-numbered than~$i$;
and in this case we say that $j$ is a {\em proper child}\/ of~$i$.
These were the key concepts in the present paper.
We now define a related but different concept:
we say that a vertex~$i$ is a \textbfit{proper vertex}
if all the descendants of~$i$, other than~$i$ itself,
are higher-numbered than~$i$.
(Equivalently, a vertex is proper in case
 {\em all}\/ of its children are proper children.)
 Note that every leaf is proper,
 and that the smallest-numbered vertex in each tree is proper.
Let us write $\propv(\scrf)$ for the number of proper vertices
in the forest $\scrf$.
Writing $\sfF_{n,k}$ for the set of forests of rooted trees
on the vertex set $[n]$ with $k$ components,
Gessel--Seo's first combinatorial interpretation is
\cite[Theorem~6.1]{Gessel-Seo_06}
\be
   P_{n,k}(a,b)
   \;=\;
   \sum_{\scrf \in \sfF_{n,k}}  a^{n - \propv(\scrf)} \, b^{\propv(\scrf) - k}
   \;.
 \label{eq.gessel-seo.first}
\ee
Gessel and Seo \cite{Gessel-Seo_06} gave two proofs of
\reff{eq.gessel-seo.first}:  one using exponential generating functions,
the other partly combinatorial.
A fully bijective proof was given by Seo and Shin \cite{Seo_07}.
Note that the symmetry $a \leftrightarrow b$ is far from obvious in
\reff{eq.gessel-seo.first};
a combinatorial explanation was recently given by Hou \cite{Hou_16}.

Yet another related concept is as follows:
We say that an edge $e = ij$ in a forest $\scrf$,
ordered so that $j$ is a child of~$i$,
is an \textbfit{ascent} if $i < j$ and a \textbfit{descent} if $i > j$.
Let us write $\asc(\scrf)$ [resp.\ $\des(\scrf)$]
for the number of ascents (resp.\ descents) in the forest $\scrf$.
Gessel--Seo's second combinatorial interpretation
\cite[Theorem~9.1]{Gessel-Seo_06}
--- a special case of a result found earlier by
 E\u{g}ecio\u{g}lu and Remmel \cite{Egecioglu_86} --- is
\be
   P_{n,k}(a,b)
   \;=\;
   \sum_{\scrf \in \sfF_{n,k}}  a^{\des(\scrf)} \, b^{\asc(\scrf)}
   \;.
 \label{eq.gessel-seo.second}
\ee
Note that the symmetry $a \leftrightarrow b$ is manifest in
\reff{eq.gessel-seo.second}:
it suffices to relabel the vertices $i \mapsto n+1-i$.

See also Drake \cite[Example~1.7.2]{Drake_08}
for another combinatorial interpretation of the polynomials $P_{n,1}(a,b)$.

Finally, it follows from either \reff{eq.gessel-seo.first} or
\reff{eq.gessel-seo.second}, using arguments identical to those used in
Sections~\ref{subsec.EGF.1} and \ref{subsec.EGF.2},
that the unit-lower-triangular matrix
$\bigl( P_{n,k}(a,b) \bigr)_{n,k \ge 0}$
is an exponential Riordan array $\scrr[F,G]$ with $F(t) = 1$.
As noted earlier, this implies (and is in fact equivalent to)
the Rothe--Pfaff--Schl\"afli identity \reff{eq.binomial_type}.

By analogy with Theorem~\ref{thm1.3},
I conjecture the following:

\begin{samepage}
\begin{conjecture}[Total positivities for the Schl\"afli--Gessel--Seo polynomials]
   \label{conj.gessel-seo}
\quad\hfill\vspace*{-1mm}
\begin{itemize}
   \item[(a)]  The unit-lower-triangular polynomial matrix
      $P(a,b) = \bigl( P_{n,k}(a,b) \bigr)_{n,k \ge 0}$
      is coefficientwise totally positive (jointly in $a,b$).
   \item[(b)]  The polynomial sequence
       $\bP = \bigl( P_{n}(x;a,b) \bigr)_{n \ge 0}$
       is coefficientwise Hankel-totally positive (jointly in $x,a,b$).
   \item[(c)]  The polynomial sequence
       $\bP^\triangle = \bigl( P_{n+1,1}(a,b) \bigr)_{n \ge 0}$
       is coefficientwise Hankel-totally positive (jointly in $a,b$).
\end{itemize}
\end{conjecture}
\end{samepage}

\noindent
I have verified part~(a) up to $15 \times 15$,
and part~(b) up to $11 \times 11$;
part~(c) is an immediate consequence of part~(b).
Conjecture~\ref{conj.gessel-seo} of course
implies analogous statements for the Rothe polynomials \reff{def.rothe},
but not conversely.

In view of the approach used in Section~\ref{sec.proofs}
to prove Theorems~\ref{thm1.1}--\ref{thm1.4},
it is natural to try to employ the same production-matrix method
to prove Conjecture~\ref{conj.gessel-seo}.
Alas, this does not work.
Straightforward computation gives for the first few rows
of the production matrix
\begin{eqnarray}
   & &
   \hspace*{-5mm}
   \Pi
   \;\eqdef\;
   P(a,b)^{-1} \Delta P(a,b)
   \;=\;
          \nonumber \\[1mm]
   & &
   \hspace*{-1cm}
\Scale[0.6]{
   \begin{bmatrix}
 0 & 1 &   &   &   &     &   \\
 0 & a+b & 1 &   &   &     &   \\
 0 & a^2+3 a b+b^2 & 2 a+2 b & 1 &   &     &   \\
 0 & a^3+7 a^2 b+7 a b^2+b^3 & 3 a^2+9 a b+3 b^2 & 3 a+3 b & 1 &     &   \\
 0 & a^4+15 a^3 b+33 a^2 b^2+15 a b^3+b^4 & 4 a^3+28 a^2 b+28 a b^2+4 b^3 & 6 a^2+18 a b+6 b^2 & 4 a+4 b & 1   &   \\
 0 & a^5+31 a^4 b+131 a^3 b^2+131 a^2 b^3+31 a b^4+b^5 & 5 a^4+75 a^3 b+165 a^2 b^2+75 a b^3+5 b^4 & 10 a^3+70 a^2 b+70 a b^2+10 b^3 & 10 a^2+30 a b+10 b^2 & 5 a+5 b   & \ddots   \\
 \vdots & \vdots & \vdots & \vdots & \vdots & \vdots &  \ddots
    \end{bmatrix}
}
          \nonumber \\
 \label{eq.prodmat.gessel-seo}
\end{eqnarray}
(this is \reff{eq.thm.riordan.exponential.production}
with $\bz = 0$ and $\ba$ given by \cite[A046802]{OEIS})
and we have
\begin{eqnarray}
   & &
   \!\!\!\!
   \pi_{41} \pi_{52} - \pi_{51} \pi_{42}
         \nonumber \\[2mm]
   & &
   \quad
   =\;
   (a^4+15 a^3 b+33 a^2 b^2+15 a b^3+b^4)(5 a^4+75 a^3 b+165 a^2 b^2+75 a b^3+5 b^4)    \nonumber \\
    & & \qquad - \, (4 a^3+28 a^2 b+28 a b^2+4 b^3) (a^5+31 a^4 b+131 a^3 b^2+131 a^2 b^3+31 a b^4+b^5)
        \nonumber \\[2mm]
   & &
   \quad
   =\;
   a^8-2 a^7 b+35 a^6 b^2+36 a^5 b^3+121 a^4 b^4+36 a^3 b^5+35 a^2 b^6-2 a b^7+b^8
        \nonumber \\[1mm]
   & &
   \quad
   \not\myge\;  0
   \;.
\end{eqnarray}
So the production matrix is not even coefficientwise TP${}_2$!
Clearly, new techniques will be needed to prove
Conjecture~\ref{conj.gessel-seo}, if indeed it is true.

We can also take Conjecture~\ref{conj.gessel-seo}(c) one step further.
Note first that $P_n(1;1,1) = f_n = (n+1)^{n-1}$
is a Stieltjes moment sequence:
it is the product of the Stieltjes moment sequences $(n+1)^n$
(see footnote~\ref{footnote.stieltjes.fn} above)
and $1/(n+1)$.
This known fact is a specialization of the claim
in Conjecture~\ref{conj.gessel-seo}(b)
that the sequence of polynomials $P_n(x;a,b)$ is
coefficientwise Hankel-totally positive.
On the other hand, we also know a stronger fact:
not only is $(n+1)^{n-1}$ a Stieltjes moment sequence,
but so is ${(n+1)^{n-1}/n!}$,
since it is the product of the Stieltjes moment sequences ${(n+1)^n/n!}$
(see again footnote~\ref{footnote.stieltjes.fn})
and ${1/(n+1)}$.
(This latter fact is stronger, because multiplication by $n!$ preserves the
 Stieltjes moment property.)
This suggests to ask whether the sequence of polynomials $P_n(x;a,b)/n!$
is coefficientwise Hankel-totally positive.
The answer is negative;
indeed, this sequence is not even coefficientwise log-convex, since
\be
   {P_0 \over 0!} \, {P_2 \over 2!} \:-\: \Bigl( {P_1 \over 1!} \Bigr)^2
   \;=\;
   {x(x+a+b) \over 2} \,-\, x^2
   \;=\;
   \half (ax + bx - x^2)
   \;\not\myge\; 0  \;.
\ee
So Conjecture~\ref{conj.gessel-seo}(b) does {\em not}\/ have an analogue
involving division by $n!$.
But Conjecture~\ref{conj.gessel-seo}(c) may:

\begin{conjecture}[Hankel-TP for the Schl\"afli--Gessel--Seo polynomials, bis]
   \label{conj.gessel-seo.bis}
\quad\hfill\vspace*{-1mm}
\begin{itemize}
   \item[(a)]  The polynomial sequence
$\bigl( P_{n+1,1}(a,b)/n! \bigr)_{n \ge 0}$
is coefficientwise Hankel-totally positive (jointly in $a,b$).
   \item[(b)]  The polynomial sequence
$\bigl( P_{n+1,1}(a,b)/(n+1)! \bigr)_{n \ge 0}$
is coefficientwise Hankel-totally positive (jointly in $a,b$).
\end{itemize}
\end{conjecture}

\noindent
I have verified parts~(a) and (b) up to $11\times 11$.

\bigskip

{\bf Remark.}
Among the Conjectures~\ref{conj.gessel-seo}(c),
\ref{conj.gessel-seo.bis}(a) and \ref{conj.gessel-seo.bis}(b),
Conjecture~\ref{conj.gessel-seo.bis}(a) is ``morally'' the strongest,
because one would expect that
multiplication by the Stieltjes moment sequences $n!$ or $1/(n+1)$
would preserve coefficientwise Hankel-total positivity
(as it does for Stieltjes moment sequences of real numbers).
But, rather suprisingly, it turns out \cite{Sokal_totalpos}
that this is {\em not}\/ a general property:
there exist coefficientwise Hankel-TP sequences $(p_n(x))_{n \ge 0}$
in the polynomial ring $\R[x]$
for which $(n! \, p_n(x))_{n \ge 0}$ and $(p_n(x)/(n+1))_{n \ge 0}$
are {\em not}\/ coefficientwise Hankel-TP.
So Conjectures~\ref{conj.gessel-seo}(c),
\ref{conj.gessel-seo.bis}(a) and \ref{conj.gessel-seo.bis}(b)
need to be considered separately.
\myendremark

\subsection[$q$-generalizations of the forest numbers]{$\bm{q}$-generalizations of the forest numbers}

It is natural to seek $q$-generalizations of the results and conjectures
in this paper.
Recall the definition of the $q$-integers
\be
   [n]_q
   \;\eqdef\; {1 - q^n \over 1-q}
   \;=\;  \begin{cases}
               0  & \textrm{if $n=0$}  \\
               1+q+q^2+\ldots+q^{n-1}  & \textrm{if $n \ge 1$}
          \end{cases}
\ee
and the $q$-factorials
\be
   [n]_q !
   \;\eqdef\;
   \prod_{i=1}^n [i]_q
   \;.
\ee
We treat $q$ as an indeterminate.
Then the $q$-binomial coefficients
\be
   \qbinom{n}{k}{q}
   \;\eqdef\;
   {[n]_q !  \over  [k]_q ! \: [n-k]_q !}
\ee
are polynomials in $q$ with nonnegative integer coefficients
\cite[Theorem~3.1]{Andrews_76}.

The simplest $q$-generalization of the forest numbers \reff{def.fnk}
simply replaces $n$ and $k$ (in the base but not in the exponent)
by $q$-integers, and the binomial coefficient by a $q$-binomial coefficient:
\be
   f_{n,k}(q)
   \;\eqdef\;
   \begin{cases}
       \delta_{k0}   & \textrm{if $n=0$} \\[2mm]
       \displaystyle
       \qbinom{n}{k}{q} \, [k]_q \, ([n]_q)^{n-k-1}
       \;=\;
       \qbinom{n-1}{k-1}{q} \, ([n]_q)^{n-k}
                     & \textrm{if $n \ge 1$}
   \end{cases}
 \label{def.fnk.q}
\ee
The triangular array $\bigl( f_{n,k}(q) \bigr)_{n,k \ge 0}$
of $q$-forest numbers begins
\be
\Scale[0.6]{
\begin{bmatrix}
 1  &     &     &     &     &       &   \\
 0  &  1  &     &     &     &       &   \\
 0  &  1 + q  &  1  &     &     &       &   \\
 0  &  1 + 2 q + 3 q^2 + 2 q^3 + q^4  &  1 + 2 q + 2 q^2 + q^3  &  1  &     &       &   &  \\
 0  &  1 + 3 q + 6 q^2 + 10 q^3 + 12 q^4 + 12 q^5 + 10 q^6 + 6 q^7 + 
   3 q^8 + q^9  &  
  1 + 3 q + 6 q^2 + 9 q^3 + 10 q^4 + 9 q^5 + 6 q^6 + 3 q^7 + q^8  &  
  1 + 2 q + 3 q^2 + 3 q^3 + 2 q^4 + q^5  &  1   &  \\
\vdots & \vdots & \vdots & \vdots & \vdots & \ddots\hspace*{-7mm}
\end{bmatrix}
   \;.
}
\ee
It follows easily from \reff{def.fnk.q}
that $f_{n,k}(q)$ is a monic self-reciprocal polynomial
of degree~$(n-1)^2 - (k-1)^2$.

Very recently, Gilmore \cite{Gilmore_inprep}
has generalized Theorem~\ref{thm1.1}(a) to the $q$-forest numbers:

\begin{theorem}[Gilmore \cite{Gilmore_inprep}]
   \label{thm.gilmore}
The unit-lower-triangular polynomial matrix $F(q) = (f_{n,k}(q))_{n,k \ge 0}$
is coefficientwise totally positive.
\end{theorem}

\noindent
Gilmore's \cite{Gilmore_inprep} method is very different from the one used here:
he uses planar networks and the Lindstr\"om--Gessel--Viennot lemma
(along the lines of \cite{Brenti_95}), not production matrices.
Indeed, the production matrix
\begin{eqnarray}
   & &
   \hspace*{-5mm}
   P
   \;\eqdef\;
   F(q)^{-1} \Delta F(q)
   \;=\;
          \nonumber \\[1mm]
   & &
   \hspace*{-6mm}
\Scale[0.8]{
   \begin{bmatrix}
0   &    1   &        &      &   \\
0   &    1+q   &    1   &      &   \\
0   &    2 q^2+2 q^3+q^4   &    q+2 q^2+q^3   &    1 &   \\
0   &    -q^2-q^3+3 q^5+6 q^6+5 q^7+3 q^8+q^9   &    -q^2+2 q^4+5 q^5+5 q^6+3 q^7+q^8   &    q^2+2 q^3+2 q^4+q^5 &  \ddots  \\
 \vdots & \vdots & \vdots & \vdots &  \ddots
    \end{bmatrix}
}
          \nonumber \\
 \label{eq.prodmat.Fnk.q}
\end{eqnarray}
is not even coefficientwise TP${}_1$;
and numerical tests strongly suggest that it is pointwise TP
(that is, for a real number $q$)
only when $q=0$ or $q=1$.

Finally, the production-matrix method {\em cannot}\/ work here
because the generalization of Theorem~\ref{thm1.1}(b) is {\em false}\/:
the row-generating polynomials
\be
   F_n(x,q)  \;=\;  \sum_{k=0}^n f_{n,k}(q) \, x^k
\ee
are {\em not}\/ coefficientwise Hankel-TP.
Indeed, the $3 \times 3$ Hankel minor
\begin{eqnarray}
   & &
   \hspace{-1cm}
   \begin{vmatrix}
      F_0(x,q)  &  F_1(x,q)  &  F_2(x,q)  \\
      F_1(x,q)  &  F_2(x,q)  &  F_3(x,q)  \\
      F_2(x,q)  &  F_3(x,q)  &  F_4(x,q)
   \end{vmatrix}
   \;\;=\;
        \nonumber \\[2mm]
   & &
   \qquad
   (-q^2+3 q^4+8 q^5+12 q^6+12 q^7+8 q^8+4 q^9+q^{10}) \, x^2
          \nonumber \\
   & &
   \qquad\quad +\:
     (-q-2 q^2-2 q^3-q^4+3 q^5+7 q^6+7 q^7+4 q^8+q^9) \, x^3
          \nonumber \\
   & &
   \qquad\quad +\:
     (-1-q+2 q^2+2 q^3-q^4-q^5) \, x^4
\end{eqnarray}
is coefficientwise nonnegative in $x$ (for real $q$) only when $q=1$.

But the Hankel-TP of the row-generating polynomials can possibly
be restored by inserting a simple additional factor: let us define
\be
   f_{n,k}^\star(q) \;\eqdef\; q^{k(k-1)/2} \, f_{n,k}(q)
 \label{def.fnk.q.star}
\ee
and then
\be
   F_n^\star(x,q)  \;=\; \sum_{k=0}^n f_{n,k}^\star(q) \, x^k
   \;.
 \label{def.Fn.q.star}
\ee
Note that this $k$-dependent factor does not change the
total positivity of the lower-triangular matrix
(it corresponds to right-multiplication by a diagonal matrix of monomials),
but it does change the row-generating polynomials.
It is not difficult to show,
using the $q$-binomial theorem \cite[Theorem~3.3]{Andrews_76},
that
\be
   F_n^\star(x,q)  \;=\; x \, \prod_{i=1}^{n-1} (q^i x + [n]_q)
   \quad\hbox{for $n \ge 1$}
   \;.
 \label{eq.q-abel}
\ee
This latter formula --- revealing $F_n^\star(x,q)$ as a kind of
``$q$-Abel polynomial'' \cite{Johnson_96,Cigler_08} ---
suggests that the numbers $f_{n,k}^\star(q)$, and not $f_{n,k}(q)$,
may be the most natural $q$-generalization of the forest numbers.

Since
\begin{subeqnarray}
   & &
   F_0^\star(x,q) \, F_2^\star(x,q) \:-\: F_1^\star(x,q)^2
   \;=\;
   (q+1) x \,+\, (q-1) x^2
          \\[2mm]
   & &
   F_1^\star(x,q) \, F_3^\star(x,q) \:-\: F_2^\star(x,q)^2
   \;=\;
          \nonumber \\[0.5mm]
   & &
   \qquad
   q^2 (q^2 + 2 q + 2) x^2  \,+\, q (q^3 + 2 q^2 - 1) x^3
                              \,+\, q^2 (q-1) x^4
   \qquad\qquad
\end{subeqnarray}
we certainly need $q \ge 1$ in order to have coefficientwise Hankel-TP${}_2$
(or even pointwise Hankel-TP${}_2$ for large positive $x$),
even if we restrict to the subsequence with ${n \ge 1}$.
But computations by Tomack Gilmore and myself
suggest that we might have coefficientwise Hankel-TP
(in $x$) of all orders whenever $q \ge 1$,
and that this might even hold coefficientwise (in the two variables)
after a change of variables $q = 1+r$.
That is, we conjecture:

\begin{conjecture}[with Tomack Gilmore]
   \label{conj.Fn.q.star}
\quad
The polynomial sequence $\bF^\star =$
\break
$\bigl( F_n^\star(x,1+r) \bigr)_{n \ge 0}$
is coefficientwise Hankel-totally positive (jointly in $x,r$).
\end{conjecture}

\noindent
I have verified this conjecture up to $10 \times 10$.

\subsection[$q$-generalizations of the Schl\"afli--Gessel--Seo polynomials]{$\bm{q}$-generalizations of the Schl\"afli--Gessel--Seo polynomials}

We can go farther and introduce $q$-generalizations
of the Schl\"afli--Gessel--Seo polynomials \reff{def.gessel-seo}.
There are several ways in which this can be done;
the best-behaved seems to be\footnote{
   Since $[i]_q + q^i [n-i]_q = [n]_q$,
   this is also a rewriting of the ``$q$-Rothe'' polynomials $a_n(x;b,h,w,q)$
   defined in \cite[Section~4]{Johnson_96}.
}
\begin{subeqnarray}
   P_0(x;y,a,b,q)   & = &  1  \\
   P_n(x;y,a,b,q)   & = &
     x \, \prod\limits_{i=1}^{n-1}
               (q^i x \,+\, y \,+\, [i]_q \, a \,+\, q^i [n-i]_q \, b)
     \quad\hbox{for $n \ge 1$} \qquad
 \label{def.gessel-seo.q}
\end{subeqnarray}
Let us define also
\be
   P_{n,k}(y,a,b,q)  \;=\;  [x^k] \, P_n(x;y,a,b,q)
     \quad\hbox{for $n,k \ge 0$}
   \;.
 \label{def.gessel-seo.Pnk.q}
\ee
Since $[i]_q + q^i [n-i]_q = [n]_q$,
it follows from \reff{eq.q-abel}/\reff{def.gessel-seo.q} that
\be
   P_{n,k}(0,1,1,q)
   \;=\;
   f^\star_{n,k}(q)
   \;.
 \label{eq.qgessel5.fnk}
\ee
So the $q$-Schl\"afli--Gessel--Seo polynomials \reff{def.gessel-seo.Pnk.q}
are a refinement of the {\em modified}\/ $q$-forest numbers
\reff{def.fnk.q.star}.
I now make the following generalization of
Conjectures~\ref{conj.gessel-seo} and \ref{conj.Fn.q.star}:

\begin{samepage}
\begin{conjecture}[Total positivities for the $q$-Schl\"afli--Gessel--Seo polynomials]
   \label{conj.gessel-seo.q}
\quad\hfill\vspace*{-1mm}
\begin{itemize}
   \item[(a)]  The unit-lower-triangular polynomial matrix
      $P(y,a,b,q) = \bigl( P_{n,k}(y,a,b,q) \bigr)_{n,k \ge 0}$
      is coefficientwise totally positive (jointly in $y,a,b,q$).
   \item[(b)]  The polynomial sequence
       $\bP = \bigl( P_{n}(x;0,a,b,1+r) \bigr)_{n \ge 0}$
       is coefficientwise Hankel-totally positive (jointly in $x,a,b,r$).
\end{itemize}
\end{conjecture}
\end{samepage}

\noindent
I have verified part~(a) up to $12 \times 12$,
and up to $14 \times 14$ when specialized to $y=0$.
I~have also verified part~(b) up to $8 \times 8$.

Please note also that in part~(b) it is crucial that we set $y=0$.
Indeed, even when $x=a=0$, the $3 \times 3$ Hankel minor
\begin{subeqnarray}
   \Delta_3^{(0)}(y,b,q)
   & \eqdef &
   \begin{vmatrix}
      P_0(0;y,0,b,q)  &  P_1(0;y,0,b,q)  &  P_2(0;y,0,b,q)  \\
      P_1(0;y,0,b,q)  &  P_2(0;y,0,b,q)  &  P_3(0;y,0,b,q)  \\
      P_2(0;y,0,b,q)  &  P_3(0;y,0,b,q)  &  P_4(0;y,0,b,q)
   \end{vmatrix}
        \\[3mm]
   & = &
   (\ldots) y^2 \:+\: (\ldots) y^3 \:-\: q^2 (q-1)^2 b^2 \, y^4
\end{subeqnarray}
can be coefficientwise nonnegative in $y$
(or even pointwise nonnegative for large positive $y$) 
only when $q=0$ or $q=1$ or $b=0$.


\begin{samepage}
\begin{problem}[Properties of the $q$-Schl\"afli--Gessel--Seo polynomials]
   \label{problem.gessel-seo.q}
\quad\hfill\vspace*{-1mm}
\begin{itemize}
   \item[(a)]  Find combinatorial interpretations of the polynomials
       \reff{def.gessel-seo.q}/\reff{def.gessel-seo.Pnk.q}
       or variants thereof,
       generalizing \reff{eq.gessel-seo.first} and \reff{eq.gessel-seo.second}.
       (Here the work of E\u{g}ecio\u{g}lu and Remmel \cite{Egecioglu_86}
        may be relevant.)
   \item[(b)]  Find ``$q$-binomial'' identities satisfied by the
       polynomials \reff{def.gessel-seo.q} or variants thereof,
       generalizing \reff{eq.binomial_type}.
       (Some partial results have been obtained by Johnson \cite{Johnson_96}.)
       In particular, do the polynomials \reff{def.gessel-seo.Pnk.q}
       form a $q$-exponential Riordan array in the sense of
       Cheon, Jung and Lim \cite{Cheon_13}?
\end{itemize}
\end{problem}
\end{samepage}

\subsection{Ordered forests of rooted trees, and functional digraphs}

An \textbfit{ordered forest of rooted trees} is simply a forest of rooted trees
in which we have specified a linear ordering of the trees.
Obviously the number of ordered forests of rooted trees
on the vertex set $[n]$ with $k$ components is
\be
   f^{\rm ord}_{n,k} \;=\;
   k! \, f_{n,k}  \;=\;
   {(n-1)! \over (n-k)!} \, k \, n^{n-k}  \;=\;
   {n! \over (n-k)!} \, k \, n^{n-k-1}
 \label{def.fnk.ord}
\ee
(to be interpreted as $\delta_{k0}$ when $n=0$).
The total number of ordered forests of rooted trees on the vertex set $[n]$ is
\be
   f^{\rm ord}_n  \,\;\eqdef\;\, \sum_{k=0}^n f^{\rm ord}_{n,k}  \,\;=\;\, n^n
 \label{def.fn.ord}
\ee
(see below for a proof).
The first few $f^{\rm ord}_{n,k}$ and $f^{\rm ord}_n$ are
%
%
\vspace*{-5mm}
\begin{table}[H]
\footnotesize
\hspace*{-5mm}
\begin{tabular}{c|rrrrrrrrr|r}
$n \setminus k$ & 0 & 1 & 2 & 3 & 4 & 5 & 6 & 7 & 8 & $n^n$ \\
\hline
0 & 1 &  &  &  &  &  &  &  &  & 1  \\
1 & 0 & 1 &  &  &  &  &  &  &  & 1  \\
2 & 0 & 2 & 2 &  &  &  &  &  &  & 4  \\
3 & 0 & 9 & 12 & 6 &  &  &  &  &  & 27  \\
4 & 0 & 64 & 96 & 72 & 24 &  &  &  &  & 256  \\
5 & 0 & 625 & 1000 & 900 & 480 & 120 &  &  &  & 3125  \\
6 & 0 & 7776 & 12960 & 12960 & 8640 & 3600 & 720 &  &  & 46656  \\
7 & 0 & 117649 & 201684 & 216090 & 164640 & 88200 & 30240 & 5040 &  & 823543  \\
8 & 0 & 2097152 & 3670016 & 4128768 & 3440640 & 2150400 & 967680 & 282240 & 40320 & 16777216  \\
\end{tabular}
\end{table}
\vspace*{-5mm}

\noindent
\!\!\cite[A066324 and A000312]{OEIS}.
The lower-triangular matrix
$F^{\rm ord} = (f^{\rm ord}_{n,k})_{n,k \ge 0}$
has the exponential generating function
\be
   \sum_{n=0}^\infty \sum_{k=0}^n f^{\rm ord}_{n,k} \, {t^n \over n!} \, x^k
   \;=\;
   {1 \over 1 \,-\, x T(t)}
   \;,
 \label{eq.fnk.ord.egf}
\ee
where $T(t)$ is the tree function \reff{def.treefn}
[compare \reff{eq.fnk.egf}].

The ordered forest numbers have another combinatorial interpretation.
Please recall that a \textbfit{functional digraph}
is a directed graph $G = (V,\vec{E})$ in which every vertex has out-degree 1;
the terminology comes from the fact that such digraphs
are in obvious bijection with functions $f \colon\, V \to V$
[namely, $\vec{ij} \in \vec{E}$ if and only~if $f(i) = j$].
Note that each weakly connected component of a functional digraph
consists of a directed cycle (possibly of length 1)
together with a collection of (possibly trivial) directed trees
rooted at the vertices of the cycle (with edges pointing towards the root).
We say that a vertex of a functional digraph is {\em cyclic}\/
if it lies on one of the cycles
(or equivalently, is the root of one of the underlying trees).
A functional digraph on the vertex set $V$ with $k$ cyclic vertices
can obviously be constructed by taking a (unordered) forest of rooted trees
on $V$ with $k$ components and then connecting the roots of those trees
into cycles --- that is, by choosing a permutation of those $k$ roots.
It follows that the number of functional digraphs on the vertex set $[n]$
with $k$ cyclic vertices is $f^{\rm ord}_{n,k}$.
This also proves that $\sum_{k=0}^n f^{\rm ord}_{n,k} = n^n$.

Since the ordered forest triangle $F^{\rm ord}$
is simply the forest triangle $F$ right-multiplied
by the diagonal matrix $\diag(k!)$,
its total positivity is an immediate consequence of
(and in fact equivalent to) Theorem~\ref{thm1.1}(a);
and since its $k=1$ column is identical to that of the forest triangle,
its Hankel-total positivity is equivalent to Theorem~\ref{thm1.1}(b).
However, no such trivial relation connects its row-generating polynomials
\be
   F^{\rm ord}_n(x)
   \;\eqdef\;
   \sum_{k=0}^n f^{\rm ord}_{n,k} \, x^k
   \;=\;
   \sum_{k=0}^n f_{n,k} \, k! \, x^k
 \label{def.Fn.ord}
\ee
to their unordered counterpart.
Nevertheless, the analogue of Theorem~\ref{thm1.2}(b) appears to be true:

\begin{conjecture}
   \label{conj.orderedforest}
The polynomial sequence
$\bF^{\rm ord} = \bigl( F^{\rm ord}_{n}(x) \bigr)_{n \ge 0}$
is coefficientwise Hankel-totally positive.
\end{conjecture}

\noindent
I have verified this conjecture up to $12 \times 12$.
Indeed, this result even appears to hold after division by $n!$:

\begin{conjecture}
   \label{conj.orderedforest.2}
The polynomial sequence
$\widetilde{\bF}^{\rm ord} = \bigl( F^{\rm ord}_{n}(x)/n! \bigr)_{n \ge 0}$
is coefficientwise Hankel-totally positive.
\end{conjecture}

\noindent
I have verified this conjecture up to $13 \times 13$.

\subsection{Functional digraphs by number of components}

Let $\psi_{n,k}$ be the number of functional digraphs
on the vertex set $[n]$ with $k$ (weakly connected) components;
obviously $\sum_{k=0}^n \psi_{n,k} = n^n$.
The first few $\psi_{n,k}$ are
%
%
\vspace*{-5mm}
\begin{table}[H]
\centering
\footnotesize
\begin{tabular}{c|rrrrrrrrr|r}
$n \setminus k$ & 0 & 1 & 2 & 3 & 4 & 5 & 6 & 7 & 8 & $n^n$ \\
\hline
0 & 1 &  &  &  &  &  &  &  &  & 1  \\
1 & 0 & 1 &  &  &  &  &  &  &  & 1  \\
2 & 0 & 3 & 1 &  &  &  &  &  &  & 4  \\
3 & 0 & 17 & 9 & 1 &  &  &  &  &  & 27  \\
4 & 0 & 142 & 95 & 18 & 1 &  &  &  &  & 256  \\
5 & 0 & 1569 & 1220 & 305 & 30 & 1 &  &  &  & 3125  \\
6 & 0 & 21576 & 18694 & 5595 & 745 & 45 & 1 &  &  & 46656  \\
7 & 0 & 355081 & 334369 & 113974 & 18515 & 1540 & 63 & 1 &  & 823543  \\
8 & 0 & 6805296 & 6852460 & 2581964 & 484729 & 49840 & 2842 & 84 & 1 & 16777216  \\
\end{tabular}
\end{table}
\vspace*{-5mm}

\noindent
\!\!\cite[A060281 and A000312]{OEIS}.
The unit-lower-triangular matrix $\Psi = (\psi_{n,k})_{n,k \ge 0}$
has the exponential generating function \cite{Knuth_89}
\be
   \sum_{n=0}^\infty \sum_{k=0}^n \psi_{n,k} \, {t^n \over n!} \, y^k
   \;=\;
   [1 - T(t)]^{-y}
   \;,
 \label{eq.hnk.egf}
\ee
where $T(t)$ is the tree function \reff{def.treefn}.
An equivalent statement is that
the unit-lower-triangular matrix $(\psi_{n,k})_{n,k \ge 0}$
is the exponential Riordan array
$\scrr[F,G]$ with $F(t) = 1$ and $G(t) = -\log[1 - T(t)]$ \cite[A001865]{OEIS}.

Let us introduce the row-generating polynomials
\be
   \Psi_n(y)  \;=\;  \sum_{k=0}^n \psi_{n,k} \, y^k
   \;.
\ee
We refer to the $\Psi_n(y)$ as the
\textbfit{functional-digraph polynomials}.\footnote{
   Knuth and Pittel \cite{Knuth_89} call them the {\em tree polynomials}\/
   because of the link \reff{eq.hnk.egf} with the tree function.
   But it seems to me that this name is potentially misleading,
   because these polynomials count functional digraphs, not trees.
}
We then have:

\begin{conjecture}[Total positivities for the functional-digraph polynomials]
   \label{conj.functional_digraph}
\quad\hfill\vspace*{-1mm}
\begin{itemize}
   \item[(a)]  The unit-lower-triangular matrix
      $\Psi = (\psi_{n,k})_{n,k \ge 0}$
      is totally positive.
   \item[(b)]  The polynomial sequence
       $\bm{\Psi} = \bigl( \Psi_n(y) \bigr)_{n \ge 0}$
       is coefficientwise Hankel-totally positive (in $y$).
   \item[(c)]  The integer sequence
       $\bm{\Psi}^\triangle = ( \psi_{n+1,1} )_{n \ge 0}$
       is Hankel-totally positive (i.e.\ is a Stieltjes moment sequence).
   \item[(d)]  The rational-number sequence
       $( \psi_{n+1,1}/(n+1)! )_{n \ge 0}$
       is Hankel-totally positive (i.e.\ is a Stieltjes moment sequence).
\end{itemize}
\end{conjecture}

\noindent
I have verified part~(a) up to $17 \times 17$,
part~(b) up to $13 \times 13$,
and parts~(c) and (d) up to $500 \times 500$
(by computing the classical S-fraction);
of course part~(c) is also an immediate consequence of either (b) or (d).

The production matrix $P = \Psi^{-1} \Delta \Psi$ is not totally positive:
its $9 \times 9$ leading principal submatrix is
\be
   P_{9 \times 9}
   \;=\;
   \begin{bmatrix}
       0 & 1 & 0 & 0 & 0 & 0 & 0 & 0 & 0 \\
       0 & 3 & 1 & 0 & 0 & 0 & 0 & 0 & 0 \\
       0 & 8 & 6 & 1 & 0 & 0 & 0 & 0 & 0 \\
       0 & 19 & 24 & 9 & 1 & 0 & 0 & 0 & 0 \\
       0 & 41 & 76 & 48 & 12 & 1 & 0 & 0 & 0 \\
       0 & 84 & 205 & 190 & 80 & 15 & 1 & 0 & 0 \\
       0 & 171 & 504 & 615 & 380 & 120 & 18 & 1 & 0 \\
       0 & 347 & 1197 & 1764 & 1435 & 665 & 168 & 21 & 1 \\
       0 & 690 & 2776 & 4788 & 4704 & 2870 & 1064 & 224 & 24 \\
   \end{bmatrix}
\ee
(this array is not in \cite{OEIS})
and the bottom-left $4 \times 4$ minor is negative:
\be
   \begin{vmatrix}
       84 & 205 & 190 & 80 \\
       171 & 504 & 615 & 380 \\
       347 & 1197 & 1764 & 1435 \\
       690 & 2776 & 4788 & 4704 \\
   \end{vmatrix}
   \;=\;
   -36570734
   \;<\; 0
   \;.
\ee

\subsection{Functional digraphs by number of cyclic vertices and number of components}

We can now combine the polynomials of the two preceding subsections
into a single bivariate polynomial:
let $\Psi_n(x,y)$ be the generating polynomial
for functional digraphs on the vertex set $[n]$
with a weight $x$ for each cyclic vertex
and a weight $y$ for each component.
Thus, $\Psi_n(x,1) = F_n^{\rm ord}(x)$ and $\Psi_n(1,y) = \Psi_n(y)$.
We refer to the $\Psi_n(x,y)$ as the
\textbfit{bivariate functional-digraph polynomials}.
They have the exponential generating function
\be
   \sum_{n=0}^\infty \Psi_n(x,y) \, {t^n \over n!}
   \;=\;
   [1 - xT(t)]^{-y}
   \;.
 \label{eq.Psixy.egf}
\ee
{}From these bivariate polynomials we can form two different
lower-triangular matrices:
\begin{subeqnarray}
   \Psi^{\rm X}  \;=\;  \big( \psi^{\rm X}_{n,k}(y) \big)_{n,k \ge 0}
   & \hbox{where} &
      \psi^{\rm X}_{n,k}(y) \;\eqdef\; [x^k] \, \Psi_n(x,y)
          \\[2mm]
   \Psi^{\rm Y}  \;=\;  \big( \psi^{\rm Y}_{n,k}(x) \big)_{n,k \ge 0}
   & \hbox{where} &
      \psi^{\rm Y}_{n,k}(x) \;\eqdef\; [y^k] \, \Psi_n(x,y)
\end{subeqnarray}
The matrix $\Psi^{\rm X}$ is not an exponential Riordan array;
but $\Psi^{\rm Y}$ is the exponential Riordan array
$\scrr[F,G]$ with $F(t) = 1$ and $G(t) = -\log[1 - xT(t)]$.
On the other hand, $\Psi^{\rm X}$
is obtained from the forest triangle $F = (f_{n,k})_{n,k \ge 0}$
by right-multiplication by the diagonal matrix
$\diag\big( (y^{\overline{k}})_{k \ge 0} \big)$,
where $y^{\overline{k}} \eqdef y (y+1) \cdots (y+k-1)$:
\be
   \psi^{\rm X}_{n,k}(y)  \;=\;  f_{n,k} \: y^{\overline{k}}
   \;.
\ee
The coefficientwise total positivity (in $y$) of the matrix $\Psi^{\rm X}$
is thus an immediate consequence of Theorem~\ref{thm1.1}(a);
and since its $k=1$ column is $y$ times that of the forest triangle,
its coefficientwise Hankel-total positivity
is equivalent to Theorem~\ref{thm1.1}(b).
For the rest, we have the following conjectures:

\begin{conjecture}[Total positivities for the bivariate functional-digraph polynomials]
   \label{conj.bivariate_functional_digraph}
\quad\hfill\vspace*{-1mm}
\begin{itemize}
   \item[(a)]  The unit-lower-triangular matrix
      $\Psi^{\rm Y} = \big( \psi^{\rm Y}_{n,k}(x) \big)_{n,k \ge 0}$
      is coefficientwise totally positive (in $x$).
   \item[(b)]  The polynomial sequence
       $\bm{\Psi} = \bigl( \Psi_n(x,y) \bigr)_{n \ge 0}$
       is coefficientwise Hankel-totally positive (in $x,y$).
   \item[(c)]  The polynomial sequence
       $\bm{\Psi}^{{\rm Y}\triangle} =
        \bigl( \psi^{\rm Y}_{n+1,1} \bigr)_{n \ge 0}$
       is coefficientwise Hankel-totally positive (in $x$).
\end{itemize}
\end{conjecture}

\noindent
I have verified part~(a) up to $16 \times 16$,
part~(b) up to $11 \times 11$,
and part~(c) up to $12 \times 12$;
of course, part~(c) is an immediate consequence of part~(b).

Here $\psi^{\rm Y}_{n,k}(x)$ enumerates functional digraphs
on the vertex set $[n]$ with $k$ components,
with a weight $x$ for each cyclic vertex.
In particular, the polynomials in the $k=1$ column,
which have exponential generating function $G(t) = -\log[1 - xT(t)]$, 
enumerate {\em connected}\/ functional digraphs
on the vertex set $[n]$ with a weight $x$ for each cyclic vertex;
equivalently, they enumerate {\em cyclically ordered}\/ forests of rooted trees
on the vertex set $[n]$ with a weight $x$ for each tree:
\be
   \psi^{\rm Y}_{n,1}(x)
   \;=\;
   \sum_{k=1}^n (k-1)! \, f_{n,k} \, x^k
   \quad\hbox{for $n \ge 1$}
   \;.
\ee
The coefficient matrix for these polynomials is \cite[A201685]{OEIS}.

\subsection{Forests by number of components and number of root descents}

The forest matrix $F = (f_{n,k})_{n,k \ge 0}$
is the exponential Riordan array $\scrr[1,G]$
with $G(t) = T(t) = \sum\limits_{n=1}^\infty n^{n-1} (t^n/n!)$.
Let us now generalize this by considering the exponential Riordan array
$\scrr[1,G] \eqdef F^\sharp = \big( f^\sharp_{n,k}(w) \big)_{n,k \ge 0}$ with
\be
   G(t)
   \;=\;
   {e^{w T(t)} - 1  \over  w}
   \;=\;
   \sum_{n=1}^\infty (w+n)^{n-1} \, {t^n \over n!}
\ee
where $w$ is an indeterminate;
this reduces to the forest matrix $F$ when specialized to $w=0$.
The triangular array $F^\sharp$ begins
\be
\Scale[0.79]{
\begin{bmatrix}
1  &  0  &  0  &  0  &  0  &  0  \\
0  &  1  &  0  &  0  &  0  &  0  \\
0  &  2+w  &  1  &  0  &  0  &  0  \\
0  &  9+6 w+w^2  &  6+3 w  &  1  &  0  &  0  \\
0  &  64+48 w+12 w^2+w^3  &  48+36 w+7 w^2  &  12+6 w  &  1  &  0  \\
0  &  625+500 w+150 w^2+20 w^3+w^4  &  500+450 w+140 w^2+15 w^3  &  150+120 w+25 w^2  &  20+10 w  &  1  \\
\vdots & \vdots & \vdots & \vdots & \vdots & \ddots
\end{bmatrix}
}
\ee

The polynomial $f^\sharp_{n,1}(w) = (w+n)^{n-1}$
has at least two combinatorial interpretations: it is
\begin{itemize}
   \item[(a)] the generating polynomial for rooted trees
      on the vertex set $[n]$ with a weight $1+w$ for each
      {\em root descent}\/
      (i.e.\ child of the root that is lower-numbered than the root)
      \cite{Chauve_99,Chauve_00,Sokal_trees_enumeration};
\end{itemize}
and it is also
\begin{itemize}
   \item[(b)] the generating polynomial for unrooted trees
      on the vertex set $[n+1]$ with a weight $w$ for each neighbor
      of the vertex~1 except one
      [see \reff{eq.Fn.explicit} and the sentence preceding \reff{eq.fnk.egf}].
\end{itemize}
Since in case~(a) the trees have size $n$,
it follows that $f^\sharp_{n,k}(w)$ counts $k$-component forests
of rooted trees on the vertex set $[n]$ with a weight $1+w$ for each
root descent.

Defining, as usual, the row-generating polynomials
\be
   F^\sharp_n(x,w)  \;=\;  \sum_{k=0}^n f^\sharp_{n,k}(w) \, x^k
   \;,
\ee
I conjecture:

\begin{conjecture}[Total positivities for the bivariate forest polynomials]
   \label{conj.bivariate_forest}
\quad\hfill\vspace*{-1mm}
\begin{itemize}
   \item[(a)]  The unit-lower-triangular matrix
      $F^\sharp = \big( f^\sharp_{n,k}(w) \big)_{n,k \ge 0}$
      is coefficientwise totally positive (in $w$).
   \item[(b)]  The polynomial sequence
       $\bm{F}^\sharp = \bigl( F^\sharp_n(x,w) \bigr)_{n \ge 0}$
       is coefficientwise Hankel-totally positive (in $x,w$).
   \item[(c)]  The polynomial sequence
       $\bm{F}^{\sharp\triangle} =
        \bigl( f^\sharp_{n+1,1}(w) \bigr)_{n \ge 0}$
       is coefficientwise Hankel-totally positive (in $w$).
\end{itemize}
\end{conjecture}

\noindent
I have verified part~(a) up to $16 \times 16$,
and part~(b) up to $11 \times 11$;
of course, part~(c) is an immediate consequence of part~(b).

In fact, more seems to be true.
Suppose that we make the shift $w = -1 + w'$,
so that the matrix begins
\be
\begin{bmatrix}
1  &  0  &  0  &  0  \\
0  &  1  &  0  &  0  \\
0  &  1+w'  &  1  &  0  \\
0  &  4+4 w'+ w^{\prime 2}  &  3+3 w'  &  1  \\
\vdots & \vdots & \vdots & \ddots
\end{bmatrix}
\ee
This matrix is {\em not}\/ coefficientwise TP${}_2$ in the variable $w'$,
or even pointwise TP${}_2$ at $w'=0$, since
\be
   (1+w')(3+3w') - 1 (4+4 w'+ w^{\prime 2})
   \;=\;
   -1 \,+\, 2w' \,+\, 2 w^{\prime 2}
   \;.
\ee
So part~(a) of Conjecture~\ref{conj.bivariate_forest}
does {\em not}\/ extend to the shifted matrix.
But parts~(b) and (c) do appear to extend:

\begin{conjecture}[Total positivities for the bivariate forest polynomials, bis]
   \label{conj.bivariate_forest.bis}
\quad\hfill\vspace*{-1mm}
\begin{itemize}
   \item[(b)]  The polynomial sequence
       $\bm{F}^\sharp = \bigl( F^\sharp_n(x,-1+w') \bigr)_{n \ge 0}$
       is coefficientwise Hankel-totally positive (in $x,w'$).
   \item[(c)]  The polynomial sequence
       $\bm{F}^{\sharp\triangle} =
        \bigl( f^\sharp_{n+1,1}(-1+w') \bigr)_{n \ge 0}$
       is coefficientwise Hankel-totally positive (in $w'$).
\end{itemize}
\end{conjecture}

\noindent
I have verified part~(b) up to $11 \times 11$;
of course, part~(c) is an immediate consequence of part~(b).

In fact, Xi Chen and I have recently {\em proven}\/
Conjecture~\ref{conj.bivariate_forest.bis}(c)
[and hence also the weaker Conjecture~\ref{conj.bivariate_forest}(c)]:
that is,

\begin{theorem}[Chen and Sokal \cite{Chen-Sokal_trees_totalpos}]
The polynomial sequence $\big( (w' +n)^n \big)_{n \ge 0}$
is coefficientwise Hankel-totally positive (in $w'$).
\end{theorem}

\noindent
The proof, which uses production-matrix methods similar to those used
in the present paper,
but for exponential Riordan arrays $\scrr[F,G]$ with $F \neq 1$,
will appear elsewhere \cite{Chen-Sokal_trees_totalpos}.
This proof does not, however, seem to extend to
Conjecture~\ref{conj.bivariate_forest.bis}(b).

Using \reff{eq.prop.riordan.exponential.production.1} and \reff{eq.star3},
it is straightforward to show that the $A$-series
for the exponential Riordan array $\scrr[1,G]$ is
\be
   A(s)
   \;=\;
   {(1+ws)^{1+1/w} \over 1 - {1 \over w} \log(1+ws)}
   \;.
\ee
Since for $w \in \R$ this is of the form \reff{eq.thm.aissen} only when $w=0$,
it follows that the production matrix $P = (F^\sharp)^{-1} \Delta F^\sharp$
is totally positive only when $w=0$.
Furthermore, since
\be
   A(s)
   \;=\;
   1 \,+\, (2+w)s \,+\, (5+2w) {s^2 \over 2!} \,+\, (16-2w^2) {s^3 \over 3!}
     \,+\, \ldots
   \;,
\ee
the production matrix $P$ is not coefficientwise (in $w$) even TP${}_1$.
So the production-matrix method used in the present paper
does not extend to proving Conjecture~\ref{conj.bivariate_forest};
other approaches will have to be devised, if indeed this conjecture is true.

\subsection{Some refinements of the Ramanujan and rooted-forest polynomials}

Let us now return to the generalized Ramanujan polynomials
$f_{n,k}(y,z)$ defined in \reff{def.fnkyz},
which enumerate forests of rooted trees
according to the number of improper and proper edges,
and the corresponding matrix $F(y,z) = (f_{n,k}(y,z))_{n,k \ge 0}$.
We saw in Propositions~\ref{prop.prodmat.Fyz}
and \ref{prop.P.factorization.Fyz}
that the production matrix of $F(y,z)$ is
\be
   P(y,z)  \;=\;  B_z \, D T_y D^{-1} \, \Delta
   \;,
\ee
which is totally positive by
Lemmas~\ref{lemma.toeplitz.power} and \ref{lemma.binomialmatrix.TP}.

Observe now that the Toeplitz matrix of powers $T_y$
is simply a special case of the inverse bidiagonal matrix
\be
   T(y_1,y_2,\ldots)
   \;=\;
   \begin{bmatrix}
      1   &     &     &     &         \\
      y_1   &  1  &     &     &            \\
      y_1 y_2 &  y_2  &  1  &     &            \\
      y_1 y_2 y_3 &  y_2 y_3  &  y_3  &  1    &            \\
      \vdots  & \vdots   & \vdots  & \vdots   & \ddots
   \end{bmatrix}
\ee
with entries $T(\bfy)_{ij} = y_{j+1} y_{j+2} \cdots y_i$ for $0 \le j \le i$.
(We call it the ``inverse bidiagonal matrix'' because it is the inverse of
 the lower-bidiagonal matrix that has 1 on the diagonal
 and $-y_1, -y_2, \ldots$ on the subdiagonal.)
It is easy to prove \cite{Sokal_totalpos}
that the inverse bidiagonal matrix $T(\bfy)$ is totally positive,
coefficientwise in the indeterminates $\bfy = (y_i)_{i \ge 1}$.
So let $F(\bfy,z) = (f_{n,k}(\bfy,z))_{n,k \ge 0}$
be the output matrix corresponding to the production matrix
\be
   P(\bfy,z)  \;=\;  B_z \, D T(\bfy) D^{-1} \, \Delta
   \;.
\ee
We then have the following generalization of
Theorem~\ref{thm1.3}(a,c):\footnote{
   For simplicity we refrain from including the variables $x$
   in part~(a), but this can of course be done:
   since the matrix $F(\bfy,z)$ is totally positive,
   so is $F(\bfy,z) \, B_x$.
}

\begin{theorem}[Total positivity of the refined Ramanujan polynomials]
   \label{thm1.3.MULTI}
\quad\hfill\vspace*{-1mm}
\begin{itemize}
   \item[(a)]  The unit-lower-triangular polynomial matrix
      $F(\bfy,z) = \bigl( f_{n,k}(\bfy,z) \bigr)_{n,k \ge 0}$
      is coefficientwise totally positive (jointly in $\bfy,z$).
   \item[(c)] The polynomial sequence
       $\bigl( f_{n+1,1}(\bfy,z) \bigr)_{n \ge 0}$
       is coefficientwise Hankel-totally positive (jointly in $\bfy,z$).
\end{itemize}
\end{theorem}

\proof
Since the production matrix $P(\bfy,z)$ is totally positive,
part~(a) follows immediately from Theorem~\ref{thm.iteration.homo}.
Moreover, Lemma~\ref{lemma.down-shifted} implies that
\be
   P'(\bfy,z) \;\eqdef\; \Delta P(\bfy,z) \Delta^{\rm T}
              \;=\;  \Delta \, B_z \, D T(\bfy) D^{-1}
\ee
is the production matrix for $F'(\bfy,z) = \Delta F(\bfy,z) \Delta^{\rm T}$.
Since $P'(\bfy,z)$ is totally positive,
part~(c) follows from Theorem~\ref{thm.iteration2bis},
because the zeroth column of $F'(\bfy,z)$
is $\bigl( f_{n+1,1}(\bfy,z) \bigr)_{n \ge 0}$.
\qed


By contrast, the analogue of Theorem~\ref{thm1.3}(b)
does {\em not}\/ hold in this generality:
that~is, the row-generating polynomials of $F(\bfy,z)$
are {\em not}\/ coefficientwise Hankel-totally positive.
Indeed, I~have been unable to find any interesting specializations
of the $\bfy$ (other than $y_i = y$ for all $i$)
in which this coefficientwise Hankel-totally positivity holds.
For example, if we take $y_i = q^i$,
then the $3 \times 3$ Hankel determinant is a degree-4 polynomial in $x$
whose coefficient of $x^4$ is
\be
   - q^2 (q-1)^2 \,+\, 3q (q-1)^2 z
   \;,
\ee
which is not coefficientwise nonnegative in $z$
for any real number $q \neq 0,1$.

Even so, Theorem~\ref{thm1.3.MULTI} shows that the polynomials
$f_{n,k}(\bfy,z)$ are of some interest.  What do they count?
Obviously they are enumerating forests of rooted trees
according to the number of proper edges
together with some refinement of the improper edges
into classes $1,2,3,\ldots$ with weights $y_1,y_2,y_3,\ldots\;$.
What are these classes?

\begin{problem}[Interpretation of the refined Ramanujan polynomials]
   \label{problem.ramanujan.MULTI}
Find a combinatorial interpretation of the
refined Ramanujan polynomials $f_{n,k}(\bfy,z)$.
\end{problem}

\bigskip

Alternatively, in the production matrix
$P(y,z) = B_z \, D T_y D^{-1} \, \Delta$
or
$P(y,\bphi) = (D T_\infty(\bphi) D^{-1}) (D T_y D^{-1}) \,\Delta$,
we could replace $T_y$ by a more general Toeplitz matrix $T_\infty(\bxi)$:
\be
   P(\bxi,\bphi)
   \;\eqdef\;
   (D T_\infty(\bphi) D^{-1}) (D T_\infty(\bxi) D^{-1}) \,\Delta
   \;.
\ee
Since $P(\bxi,\bphi) = P(\bphi * \bxi, 0)$,
it is immediate that the polynomials $f_{n,k}(\bxi,\bphi)$
generated by the production matrix $P(\bxi,\bphi)$
possess all the properties asserted in Theorem~\ref{thm1.4}
whenever both $\bphi$ and $\bxi$ are Toeplitz-totally positive.
Furthermore, since $f_{n,k}(\bxi,\bphi) = L_{n,k}(\bphi * \bxi)$
by Proposition~\ref{prop.generic_lah},
these polynomials have a trivial interpretation
as enumerating forests of increasing rooted trees
with a weight $m! \, (\bphi * \bxi)_m$ for each vertex with $m$ children.
But we would like, rather, an interpretation
in terms of forests of {\em general}\/ (not necessarily increasing)
rooted trees and that reduces to our original definition
of the rooted-forest polynomials $f_{n,k}(y,\bphi)$ when $\xi_\ell = y^\ell$:

\begin{problem}[Interpretation of the refined rooted-forest polynomials]
   \label{problem.refined.rooted-forest}
Find a combinatorial interpretation of the polynomials $f_{n,k}(\bxi,\bphi)$
in which each vertex with $m$ proper children gets a weight $m! \, \phi_m$
and weights $\bxi$ are somehow assigned to the improper children/edges.
\end{problem}

\noindent
In this context, see the Remark at the end of Section~\ref{subsec.EGF.3}.

\section*{Acknowledgments}

I wish to thank Xi Chen, Bishal Deb, Tomack Gilmore and Mathias P\'etr\'eolle
for helpful conversations.

This work was immeasurably facilitated by
the On-Line Encyclopedia of Integer Sequences \cite{OEIS}.
I warmly thank Neil Sloane for founding this indispensable resource,
and the hundreds of volunteers for helping to maintain and expand it.

This research was supported in part by
Engineering and Physical Sciences Research Council grant EP/N025636/1.

\end{document}